\documentclass[11pt]{amsart}
\usepackage{amssymb}
\usepackage{amsmath}
\usepackage{amsthm}
\usepackage{times}
\usepackage{changepage}
\usepackage{caption}
\usepackage{geometry}
\usepackage[curve,all]{xy}
\xyoption{curve}
\xyoption{line}
\xyoption{arc}
\xyoption{color}
\usepackage{array}
\usepackage{enumitem}
\usepackage[dvipdfmx]{graphicx}
\usepackage[dvipdfmx]{color}
\usepackage{color}
\usepackage{amsthm}
\newtheorem{theorem}{Theorem}[section]
\newtheorem{lemma}[theorem]{Lemma}
\newtheorem{corollary}[theorem]{Corollary}
\newtheorem{proposition}[theorem]{Proposition}
\newtheorem{prop}[theorem]{Proposition}

\usepackage{subfig}
 
\theoremstyle{definition}

\theoremstyle{remark}
\newtheorem{remark}[theorem]{Remark}

\newtheorem{assumption}[theorem]{Assumption}
\numberwithin{equation}{section}

\newcommand{\calC}{\mathcal{C}}

\newcommand{\calS}{\mathcal{S}}

\newcommand{\bfA}{\mathbf{A}}
\newcommand{\bfC}{\mathbf{C}}
\newcommand{\bfF}{\mathbf{F}}
\newcommand{\bfP}{\mathbf{P}}
\newcommand{\bfQ}{\mathbf{Q}}
\newcommand{\bfR}{\mathbf{R}}
\newcommand{\bfZ}{\mathbf{Z}}
\newcommand{\la}{\langle}
\newcommand{\ra}{\rangle}

\def\Aut{{\text{Aut}}}

\def\Pic{{\text{Pic}}}
\def\Ker{{\text{Ker}}}
\def\Im{{\text{Im}}}

\def\ss{{\text{ss}}}

\def\rank{{\text{rank}}}
\def\det{{\text{det}}}

\def\rank{{\text{rank}}}

\def\tor{{\text{tor}}}
\def\Hom{{\text{Hom}}}
\def\End{{\text{End}}}

\def\NS{{\rm{NS}}}
\def\Triv{{\rm{Triv}}}

\def\Km{{\text{Km}}}

\def\I{{\text{I}}}
\def\II{{\text{II}}}
\def\III{{\text{III}}}

\def\Pic{{\text{Pic}}}
\newcommand{\frakS}{\mathfrak{S}}

\setenumerate{
label=(\arabic*),
font=\normalfont
}

\begin{document}
\title [Automorphism groups of Kummer surfaces]{The automorphism groups of Kummer surfaces in characteristic two
and their complex analogues}

\author{Shigeyuki Kond\=o}
\address{Graduate School of Mathematics, Nagoya University, Nagoya 464-8602, Japan}
\email{kondo@math.nagoya-u.ac.jp}

\author{Shigeru Mukai}
\address{Research Institute for Mathematical Sciences, Kyoto University, Kyoto 606-8502, Japan}
\email{mukai@kurims.kyoto-u.ac.jp}

\thanks{This work was partly supported by JSPS Grant-in-Aid for Scientific Research (A) No.20H00112, (B) No.25K00906, and by the Research Institute for Mathematical Sciences, an International Joint Usage/Research Center located in Kyoto University.}

\begin{abstract}
We calculate the automorphism group of the Kummer surface associated with a curve of genus 2 or the product of two elliptic curves in characteristic two under the assumption that the Kummer surface is a $K3$ surface.
Moreover we discuss the complex $K3$ surfaces with the same Picard lattice as these Kummer surfaces. 
The paper has two appendices. 
\end{abstract}

\maketitle

\begin{center}
APPENDIX B BY SHIGERU MUKAI
\end{center}

\bigskip
\bigskip

\begin{center}
{\it Dedicated to Professor Yukihiko Namikawa on the occasion of his 80th birthday}
\end{center}

\bigskip
\bigskip

\section{Introduction}\label{sec1}
Let $C$ be a complex non-singular curve of genus 2, $J(C)$ the Jacobian of $C$ and $\iota$ the inversion
of $J(C)$.  Then the quotient surface $S=J(C)/\la \iota \ra$ has sixteen nodes and sixteen tropes (= double conics) and is called a Kummer quartic surface (via $|2\Theta|$ map).  
Its minimal non-singular model $X$ is a $K3$ surface called the Kummer surface associated with $C$.
There are two sets of mutually disjoint sixteen $(-2)$-curves (= non-singular rational curves) forming a Kummer $(16_6)$-configuration, that is, each member of one set meets exactly six members of the other.
One set consists of exceptional curves over sixteen nodes and the other of the proper transforms of sixteen 
tropes.

It is classically known that many birational automorphisms of $S$ exist: 
$(\bfZ/2\bfZ)^5$ generated by involutions induced from 
translations of $J(C)$ by sixteen 2-torsion points of $J(C)$ and 
the dual map $S\to S^\vee$ $(\cong S)$ called a
switch; sixteen projections from a node and sixteen correlations which are the duals of 
the projections (Klein \cite{Kl1}); 60 Hutchinson--G\"opel involutions associated with G\"opel tetrads
(Hutchinson \cite{Hut2}).  Later, in his paper \cite{Keum}, Keum found 192 automorphisms of infinite order.  Then 
the first author \cite{Kondo2} proved that the automorphism group $\Aut(X)$ for $X$ with Picard number 17  (called Picard general) 
is generated by these automorphisms, based on a theory of hyperbolic reflection groups due to Conway \cite{Co}
and Borcherds \cite{Bor1}, by showing that there exists a finite polyhedron $\calC(X)$ in the hyperbolic space 
$\NS(X)_\bfR=\NS(X)\otimes \bfR$ with $316 (= 32+ 32 + 60+192)$-faces which is a fundamental domain of 
$\Aut(X)$.  The faces correspond to
$32$ $(-2)$-curves, $32$ projections and correlations, $60$ Hutchinson--G\"opel involutions, $192$ 
Keum's automorphisms.  Later Ohashi \cite{Ohashi} pointed out the existence of 
192 Hutchinson--Weber involutions associated with Weber hexads (Hutchinson \cite{Hut1}) 
which work well instead of Keum's automorphisms. 
We remark that the existence of Keum's automorphism depends on
the Torelli type theorem for $K3$ surfaces, but other arguments are characteristic free.  Thus the above results hold for a Picard general $X$ defined over an algebraically closed 
field in positive characteristic $p\ne 2$.  
For completeness, we give its proof in Appendix \ref{OddChar}.

Later, Keum and the first author \cite{KeKo} gave similar results for complex 
Kummer surfaces associated with the product of two elliptic curves.  
We summarize the results as in Table \ref{table0}.
{\small
\begin{table}[ht]
 \centering
  \begin{tabular}{|c|c|c|c|}
   \hline  
Abelian surface  & $\NS(X)$ & $R$ & Faces of $\calC(X)$ \\ \hline 
$J(C)$   &$U\oplus D_7\oplus 2D_4$& $A_3\oplus 6A_1$ & $32 (-2)+ 32(-1)+60(-1)+192({-3\over 4})$ \\
\hline
$E\times F$ $(E\nsim F)$ &$U\oplus 2D_8$& $2D_4$ & $24(-2)+24(-1)$ \\
  \hline
$E\times E$ &$U\oplus D_8\oplus D_9$& $D_4\oplus A_3$ & $28(-2)+3(-1)+36(-1)+8({-5\over 4})+ 72({-1\over 4})$\\  
\hline
$E_\omega\times E_\omega$ &$U\oplus D_{12}\oplus E_6$& $D_4\oplus A_2$ & $32(-2)+72(-1)+8({-4\over 3})+ 96({-1\over 3})$ \\
\hline
$E_{\sqrt{-1}}\times E_{\sqrt{-1}}$ &$U\oplus E_8\oplus 2D_5$& $2A_3$ & $40 (-2)+ 40(-1)+64({-5\over 4})+160({-1\over 2})+ 320({-1\over 4})$ \\
\hline
\end{tabular}
\caption{Complex Kummer surfaces}\label{table0}
\end{table}
}

\noindent
In Table \ref{table0}, $R$ is a root sublattice of an even unimodular lattice $\II_{1,25}$ 
of signature $(1,25)$ whose orthogonal complement is isomorphic to $\NS(X)$, and $\calC(X)$ is 
a finite polyhedron in $\NS(X)_\bfR$ which is the restriction of Conway's fundamental domain $\calC$ 
of the reflection group of $\II_{1,25}$ and is
a fundamental domain of $\Aut(X)$ up to finite groups.
For example, $32 (-2)$ in Table \ref{table0} means that there are 32 faces each of them is defined by a hyperplane perpendicular to a vector with norm $-2$.  We denote by $E\nsim F$ when two elliptic curves $E$ and $F$ are non-isogenous.
An elliptic curve $E_\omega$ (resp. $E_{\sqrt{-1}}$) is the one with a primitive third root of unity (a primitive fourth root of unity) as its period.

On the other hand, in case $p=2$, the situation is entirely different.  The Kummer quartic surface has four rational double points of type $D_4$ for $J(C)$ ordinary (i.e. $p$-rank 2 case), two rational double points of type $D_8$ for
$J(C)$ of $p$-rank 1 and an elliptic singularity for $J(C)$ supersingular (i.e. $p$-rank 0).  

The main purpose of this paper is to give a set of generators of the automorphism groups of several types
of Kummer surfaces in characteristic 2 as in the case of complex Kummer surfaces
(Theorems \ref{MainCurve}, \ref{MainCurvep-rank1}, \ref{MainProd1}, \ref{MainProd2}, \ref{MainProd3}).
We summarize our calculation in Table \ref{table1}.  
{\small
\begin{table}[ht]
 \centering
  \begin{tabular}{|c|c|c|c|}
   \hline  
Abelian surface  & $\NS(X)$ & $R$ & Faces of $\calC(X)$\\ \hline 
$J(C)$,\ ordinary  &$U\oplus E_{8}\oplus D_4\oplus A_3$& $D_4\oplus D_5$ & $20 (-2)+ 4(-1)+6(-1)+8({-3\over 4})$ \\ \hline
$J(C)$,\ $p$-rank 1&$U\oplus E_8\oplus D_7$& $D_9$ & $18(-2)+2(-1)$ \\
\hline
$E\times F$ $(E\nsim F)$, ordinary &$U\oplus E_8\oplus D_{8}$& $D_8$ & $20(-2)+2(-1)$ \\
\hline
$E\times E$, ordinary &$U\oplus E_8\oplus D_9$& $D_7$ & $22(-2)+3(-1)+12({-1\over 4})$ \\
\hline
$E\times F$, $p$-rank 1 &$U\oplus E_8\oplus E_8$& $E_8$ & $19 (-2)$ \\
\hline
\end{tabular}
\caption{Kummer surfaces in characteristic 2}\label{table1}
\end{table}
}

\noindent
Here we assume that all cases are Picard general, that is, the N\'eron--Severi group of a Kummer surface is as minimal as possible.
We exclude the case of the algebraic closure $\overline\bfF_2$ of  $\bfF_2$.  
In this case, the Picard number is even (e.g. see the proof of \cite[Chap.17, Cor. 2.9]{Huy}) by the Tate conjecture \cite{ArtinSwinnerton}, \cite{Totaro}).
Also we exclude the cases that $C$ is a supersingular curve of genus 2 and that both $E$ and $F$ are supersingular elliptic curves because the associated Kummer surfaces are rational, not $K3$ surfaces.  

Note that only $(-2)$-faces appear in the last case of Table \ref{table1}. This implies that
$\calC(X)$ is nothing but the ample cone of $X$ and $\Aut(X)$ is finite in this case (Theorem \ref{MainProd3}).
In the remaining cases,
$\calC(X)$ is a proper subdomain of the ample cone and $\Aut(X)$ is infinite.

We also remark that in Table \ref{table1}, $J(C)$ (resp. $E \times F$) with $p$-rank 1 and of $p$-rank 2 have the same Picard number 1 (resp. 2).
This implies that the associated Kummer surfaces have the same Picard number.  However, their quotient singularities are different,
i.e., 2$D_8$-singularities and $4D_4$-singularities, which makes the difference of their Neron-Severi lattices $U\oplus E_8\oplus D_7$ and $U\oplus E_8\oplus D_4\oplus A_3$
(resp. $U\oplus E_8\oplus E_8$ and $U\oplus E_8\oplus D_8$).  

When $C$ is a supersingular curve of genus 2 and both $E$ and $F$ are supersingular elliptic curves in characteristic 2, 
the method as above can not be applicable to these cases.  
It would be interesting to study supersingular case in more detail.

Note that a complex $K3$ surface $X$ with the same N\'eron--Severi lattice as in Table \ref{table1} is not a Kummer surface.
However the method in this paper works well in these cases, too.  Thus almost the same result holds for these $K3$ surfaces, which is the secondary purpose of this paper.
The last case  $\NS(X)\cong U\oplus E_8\oplus E_8$  is well known  and studied in many papers under the name {\it Shioda-Inose} $K3$ surfaces.
The first and second cases are realized by quartic surfaces with the same singularities $4D_4, 2D_8$ and with the same number  ($= 4, 2$)  of tropes  over  $\bfC$, respectively
(see \S \ref{sec7}, \ref{Qw2D8}).
In the case of  $2D_8$, we give an explicit equation of a 3-dimensional 
family of quartic surfaces both in characteristic 0 and 2, 
which looks new (see the equations \eqref{lift of Duquesne 2} and \eqref{Gen'dDuquesne}).
In the case of characteristic 2, it
contains Duquesne's Kummer quartic surfaces with $p$-rank 1 in subsection \ref{Qw2D8-2}.
The third case is realized by the $K3$ cover of an Enriques surface studied by Barth-Peters~\cite{BP} (see \S \ref{DQ4D4}).
The fourth case  $\NS(X)\cong U\oplus E_8\oplus D_9$  is a special case of the third and also known as a {\it mirror quartic} (see \S \ref{Q6A3}).

We summarize the complex case in the following table.
Here the last column shows which abelian surfaces  $A$  have the same transcendental lattices with  $X$.
$U$-{\it polarized} (for $R = E_8$)  means   $A \simeq E \times F$, and 
 $U(2)$-{\it polarized} (for $R = D_8$)  means that $A$ is a double cover of  $E \times F$,
 where both  $E$  and $F$  are elliptic curves. 
$E \times (E/\la a \ra)$  with a 2-torsion $a \in E$  (for $R = D_7$)  is both $U$- and $U(2)$-polarized (see Remark~\ref{SI3}).

{\small
\begin{table}[ht]
 \centering
  \begin{tabular}{|c|c|c|c|}
   \hline  
 $R$ & Equation & Remark & Abelian surf. \\ \hline 
 $D_4\oplus D_5$ &
$A(x^2y^2+z^2w^2) + B(x^2z^2+y^2w^2)+ C(x^2w^2 + y^2z^2)$  & \eqref{Heisenbergwith E=0}& None \\
& $+ Dxyzw=0$ &&  \\ \hline
 $D_9$ & $\{q_1(x, z) + yw\}\{q_2(x, z) + yw\}   + x^2zw + xyz^2  = 0$ & \eqref{lift of Duquesne 2}& deg. 4 polarized \\
\hline
 $D_8$ & $w^2 =  \left(x + \frac{a}x\right) + b + \left(y + \frac{c}y\right)$ & \eqref{4D4-eq} & $U(2)$-polarized \\
\hline
 $D_7$ & $(x+y+z+w)^4 = \lambda xyzw$ & \eqref{mirror2} & $E \times (E/\la a \ra)$ \\
\hline
 $E_8$ & $y^2zw - 4x^3z + 3axzw^2 - \frac12(z^2w^2+w^4)+bzw^3 =0$ & Inose\cite{Ino} & $U$-polarized\\
\hline
\end{tabular}
\caption{Complex analogues}\label{table2}
\end{table}
}

The main method to calculate the automorphism group 
in characteristic 2 is the same as 
that over the complex numbers. However, 
in this paper, we use Curtis's MOG (miracle octad generator) to compute 
the polyhedron $\calC(X)$ (\S \ref{Leech}).  
Also we employ a unified way to describe a set of generators 
of automorphism group 
$\Aut(X)$ which might be applicable to other cases: for the vector $\delta$ with norm not equal to $-2$ defining a face of $\calC(X)$, 
we need
to find an involution of $X$ sending $\delta$ to $-\delta$.  
To do this, we will find a suitable elliptic pencil with a section 
whose inversion automorphism is the desired one
(\S\ref{Mordell}).  This method was used in Ujikawa \cite{Ujikawa}.

The following is the plan of this paper.
In Section \ref{sec2}, we recall Kummer surfaces in characteristic 2, Borcherds' \cite{Bor1} and Conway's \cite{Co} theorems on hyperbolic reflection groups, and Mordell-Weil groups of elliptic pencils with a section.
In Section \ref{sec3}, we consider the case of a Picard general Kummer surface $X$ 
associated with an ordinary curve of genus 2.  
In Section \ref{sec4}, we give a set of generators of the automorphism group of a Picard general Kummer
surface associated with a curve of genus 2 and $p$-rank 1.  
In Sections \ref{sec5} and \ref{sec6}, we study a Kummer surface associated with 
the product of two elliptic
curves $E, F$ where $E$ is ordinary (i.e. $p$-rank 1) and $F$ is ordinary or supersingular ($p$-rank 0).  
In Sections \ref{sec7}--\ref{Q6A3}, we discuss the complex $K3$ surfaces 
with the same N\'eron--Severi lattice as in Table \ref{table1}.

\medskip
{\bf Notations and convention.}
Throughout the paper, we assume that the base field is an algebraically closed field.  
However we exclude the algebraic closure $\overline\bfF_2$ of  $\bfF_2$.

A {\it lattice} is a free $\bfZ$-module $L$ of finite rank with 
a $\bfZ$-valued non-degenerate symmetric bilinear form $x\cdot y$ ($x, y \in L$).
We denote by $L^*$ the dual of $L$, that is, $L^*=\Hom(L,\bfZ)$.  
A lattice is called {\it unimodular} if $L \cong L^*$.
A lattice $L$ is called {\it even} if $x\cdot x$ is even for any $x\in L$.  A vector $x\in L$ with 
$x\cdot x=-2$ is called a $(-2)$-vector or a {\it root}.  A {\it root lattice} is 
a {\it negative} definite lattice generated by $(-2)$-vectors.  
We denote by $A_m, D_n$ or $E_k$ $(m\geq 1, n\geq 4, k=6,7,8)$ 
a root lattice defined by the Cartan matrix of the same type.  
Denote by $\tilde{A}_m, \tilde{D}_n$ or $\tilde{E}_k$ an extended root lattice which is also used
for a singular fiber of the same type of an elliptic pencil. 
We denote by $U$ the hyperbolic plane, that is, an even unimodular lattice of signature $(1,1)$.
For lattices $L, M$, we denote by $L\oplus M$
the orthogonal direct sum, and by $mL$ the orthogonal direct sum of $m$ copies of $L$.  
For an even lattice $L$, we define $A_L=L^*/L$ and $q_L: A_L \to {\bf Q}/2{\bf Z}$, 
$q_L(x+L)=x\cdot x + 2{\bf Z}$
which are called the discriminant group and the discriminant quadratic form of $L$, respectively.

A non-singular rational curve on a $K3$ surface is called a $(-2)$-curve for simplicity.  A dual graph of
some $(-2)$-curves is the following simplicial complex $\Gamma$:  (i) the set of vertices, denoted by
$\bullet$, consists of $(-2)$-curves; (ii) two vertices are joined by $m$-tuple line if the intersection number of the corresponding $(-2)$-curves is $m$ ($m=0,1,2$).  For example, see Figure \ref{20curves}.  We add vertices, denoted by $\circ$, corresponding to some divisors to a dual graph of $(-2)$-curves.
If $m \geq 3$ or $m < 0$, then two vertices are joined by
a dotted line with $m$.
For example, see Figure \ref{p-rank1dualgraph}.

For an Abelian variety  $A$,  $\Aut_0(A)$  denotes the automorphism group of  $A$  as group variety.
It is a subgroup of  the automorphism group $\Aut(A)$  as variety and also the unit group of the endomorphism ring $\End(A)$.

\medskip
{\bf Acknowledgement.}
This joint work was started when both of authors stayed in the Institute for Basic Sciences for Complex Geometry (IBS-CCG) in Daejeon, South Korea, in Spring 2024.
The authors are very grateful for the invitation and its hospitality.
The authors thank Toshiyuki Katsura and Hisanori Ohashi for useful comments.  
The authors also would like to thank the referee
who gave them many helpful comments and advice.

\section{Preliminaries}\label{sec2}

\subsection{Kummer surfaces associated with curves $C$ 
of genus 2 in characteristic 2 whose Jacobian have $p$-rank 2 }\label{Kummer}
First we recall the Kummer quartic surface $S$ associated with $C$.
The equation of Igusa's canonical model of $C$ is given by
\begin{equation}\label{IgusaCan}
y^2+ y = a x + b x^{-1} + c (x+1)^{-1}
\end{equation}
where $a, b, c$ are non-zero constants.  Thus their moduli space has dimension three.

\begin{prop}\label{KummerQuartic}{\rm (Laszlo--Pauly \cite{LP} \cite{LP2}, Katsura--Kond\=o \cite[Remarks 4.3, 6.2]{KK2}; See Appendix B also.)}
\begin{enumerate}
\item[(i)]\ The Kummer quartic surface $S$ associated with $C$ is given by the equation
\begin{equation}\label{KummerQuarticEq}
\begin{array}{l}
(\sqrt{c}(xy+zw) + \sqrt{b}(xz+yw)+ \sqrt{a}(xw + yz))^2+ xyzw=0.
\end{array}
\end{equation}
\item[(ii)]\ The quartic surface $S$ has exactly four rational double points 
\[p_1=(1, 0, 0, 0),\ p_{2}=(0, 1, 0, 0),\ p_{3}=(0, 0, 1, 0),\ p_{4}=(0, 0, 0, 1)\]
of type $D_4$
and contains four tropes $T_1, T_2, T_3, T_4$ defined by $x=0, y=0, z=0, w=0$ respectively.
\end{enumerate}
\end{prop}

\begin{remark}
Katsura--Kond\=o \cite{KK2} obtained the equation (\ref{KummerQuarticEq}) as the Kummer quartic surface 
associated with a quadratic line complex.
\end{remark}

The Cremona transformation
\begin{equation}\label{CremonaOrdinary}
\sigma : (x, y, z, w)\dashrightarrow \left(yzw, xzw, xyw, xyz\right)
\end{equation}
preserves the equation of the Kummer quartic surface and
interchanges $p_i$ and $T_i$.  
Also there are the following involutions of $S$ which generate $(\bfZ/2\bfZ)^2$:
\begin{equation}\label{translation}
\left\{ \begin{array}{l}
\varphi_1: (x, y, z, w)\to (y, x, w, z),\\
\varphi_2: (x, y, z, w)\to (z, w, x, y).\\
\end{array}\right.
\end{equation}

\begin{remark}
If $a = b$, $b=c$ or $c=a$, then the following projective transformation
induces an involution of $S$:
\begin{equation}\label{transformation}
\left\{ \begin{array}{l}
\psi_{a,b} : (x, y, z, w)\to (x, y, w,z)\quad (a=b),\\
\psi_{b,c} : (x, y, z, w)\to (x, z, y,w)\quad (b=c),\\
\psi_{c,a} : (x, y, z, w)\to (x, w, z,y)\quad (c=a).\\
\end{array}\right.
\end{equation}
These transformations are induced from automorphisms of the curve $C$: 
\begin{equation}\label{}
\left\{ \begin{array}{l}
(x,y)\to (x^{-1}, y +c_0) \quad (a = b,\ c_0^2 + c_0 = c),\\ 
(x,y)\to (x+1, y + a_0)\quad (b=c,\ a_0^2 + a_0 = a),\\
(x,y)\to (x(x+1)^{-1}, y + b_0)\quad (c = a, \ b_0^2 + b_0 = b)
\end{array} \right.
\end{equation}
(Ancochea \cite{An}, Igusa \cite{Ig}).
\end{remark}

\subsection{Kummer surfaces associated with curves $C_1$ of genus 2 whose Jacobian have $p$-rank 1}\label{Kummer2}

Next we recall the Kummer surface $S_1$ associated with $C_1$.
The equation of Igusa's canonical model of $C_1$ is given by
\begin{equation}\label{IgusaCan2}
y^2+ y = x^3 + a x + b x^{-1}
\end{equation}
where $a, b \in k$ with $b\ne 0$.
Thus their moduli space has dimension two.

\begin{prop}\label{KummerQuartic2}
\begin{enumerate}
\item[(i)]\ {\rm (Duquesne \cite[\S 3]{Du}; see Katsura--Kond\=o \cite[\S 7.4]{KK2})}\ The Kummer quartic surface $S_1$ is given by the equation
\begin{equation}\label{KummerEq1}
b^2x^4+a^2x^2z^2 + x^2zw + xyz^2 + y^2w^2 + z^4=0.
\end{equation}
\item[(ii)]\ {\rm (Shioda \cite{Shioda}, Katsura \cite{Ka})}\ 
The surface $S_1$ has exactly two singular points 
\[p_1=[0,0,0, 1],\ \ p_2=[0, 1, 0, 0]\]
of type $D_8$  and contains two tropes $T_1, T_2$ cut out by
the hyperplane sections $x=0$ and $z=0$, respectively.
Two tropes meet at $p_1$ and $p_2$.
\end{enumerate}
\end{prop}

There exist a Cremona involution
\begin{equation}\label{Cremonaprank1}
\sigma : [x,y,z,w]\dashrightarrow [xz^2, yz^2, b x^2z, b x^2w]
\end{equation}
and a projective linear transformation 
\begin{equation}\label{translprank1}
\varphi : [x, y, z, w] \mapsto [z, w, b x, b y].
\end{equation}
which act on $S_1$ (Dolgachev--Kondo \cite{DK2}).

\subsection{The Leech lattice and the reflection group of a hyperbolic lattice}\label{Leech}

Let $\Lambda$ be the Leech lattice, that is, an even unimodular negative definite lattice of rank 24 without $(-2)$-vectors.
Let ${\II}_{1,25}$ be an even unimodular lattice of signature $(1,25)$ which is unique up to isomorphisms.
Note that ${\II}_{1,25}\cong U\oplus \Lambda$.  We fix an orthogonal decomposition 
\[\II_{1,25} = U \oplus  \Lambda,\]
and denote each vector
$x\in {\II}_{1,25}$ by $(m,n,\lambda)$, where 
$\lambda\in \Lambda$, and $x = mf+ng +\lambda,$ with $f,g$ being 
generators of $U$ satisfying $f^2 = g^2 = 0$ and $f\cdot g= 1$.
Let $\rho=(1,0,0)$ which is called {\it the Weyl vector}.  A $(-2)$-vector $r$ is called {\it a Leech root} if $r\cdot \rho =1$.
Note that a Leech root $r$ can be written as $r = (-1-\frac{\langle \lambda,\lambda \rangle} {2}, 1, \lambda)$. 
We denote by $\Delta$ the
set of all Leech roots.   
Let $W(\II_{1,25})$ be the subgroup generated by
reflections in the orthogonal group $O(\II_{1,25})$ of $\II_{1,25}$. Let $P^+$ be a
connected component of
\[\{ x \in \bfP({\II_{1,25}}\otimes \bfR) : x^2 > 0 \}.\] 
Then $W(\II_{1,25})$ acts naturally on $P^+$. Conway \cite{Co} showed that
a fundamental domain with respect to this action of
$W(\II_{1,25})$ on $P^+$ is given by
\[\calC = \{ x \in P^+ : x\cdot r > 0, \ r\in \Delta \}.\]

Let $R$ be a root sublattice of $\II_{1,25}$ generated by some Leech roots.  Let $L$ be the orthogonal complement of $R$ in $\II_{1,25}$ which has signature $(1, 25 -\rank R)$.  Let $\calC(L)$ be the restriction
of $\calC$ under the embedding $L\subset \II_{1,25}$.  
We denote by $w$ the projection of $\rho$ into $L^*$.
Then, Borcherds proved the following:

\begin{proposition}\label{Bor}{\rm (Borcherds \cite[Lemmas 4.1, 4.2]{Bor1})}\ \  
$\calC(L)$ contains $w$, in particular it is non-empty, and has only a finite number of faces.
\end{proposition}

We now recall the details of Leech lattice.  Let $\Omega$ be a set of twenty-four 
points.  Let $\bfR^{24}$ be spanned by the orthonormal basis $\nu_i (i\in \Omega)$.
Define $\nu_S$ to be $\sum_{i\in S} \nu_i$ for a subset $S \subset \Omega$. 
A {\it Steiner system} $\calS(5,8,24)$ is a set consisting of 
eight-element subsets of $\Omega$ such that any five-element subset belongs to a unique element 
of $\calS(5,8,24)$. An eight-element subset in $\calS(5,8,24)$ is called an {\it octad}. Then
$\Lambda$ is defined as a lattice in $\bfR^{24}$, equipped with inner product $-\frac{x\cdot y}{8}$, 
generated by the vectors $\nu_\Omega-4\nu_\infty$ and $2\nu_K$, 
where $K$ belongs to the Steiner system $\calS(5,8,24)$.

In the following we will use MOG (miracle octad generator) for describing the elements of $\Lambda$ (see Conway--Sloane \cite[Chap.11, \S 11, Chap.23, \S 4]{CS} for the details).  Consider a box consisting 24 points.  We label three special 
points $\I, \II, \III$ called {\it Romans}.  The remaining 21 points are identified with 21 points on 
$\bfP^2(\bfF_4)$ the projective plane defined over the finite field $\bfF_4$.

\begin{figure}[ht]
 \centering
  \begin{tabular}{|c|c|c|}
   \hline  
$\infty_\infty$  $\infty_0$  & $\ \cdot \ $ \ \ $\ \cdot \ $ & $\ \cdot \ $ \ \ $\ \cdot \ $ \\ 
I \ \ \ \ \ $\infty_1$  & $\ \cdot \ $ \ \ $\ \cdot \ $ & $\ \cdot \ $ \ \ $\ \cdot \ $ \\     \hline 
II \ \ \ \ $\infty_\omega$  & $\ \cdot \ $ \ \ $\ \cdot \ $ & $\ \cdot \ $ \ \ $\ \cdot \ $ \\  
III \ \ \ $\infty_{\bar\omega}$  & $\ \cdot \ $ \ \ $\ \cdot \ $ & $\ \cdot \ $ \ \ $\ \cdot \ $  \\
\hline
\end{tabular}
\caption{MOG as a projective plane (\cite[Figure 11.22]{CS})}\label{MOG}
\end{figure}
The five points $\infty_s$ on the left denote the points at infinity with slope $y/x = s$  for  $s \in \bfP^1(\bfF_4)$.
The rest are the 16 points of the affine plane  $\bfA^2(\bfF_4)$.

\begin{proposition}\label{MOGoctads}{\rm (\cite[p.315]{CS})}
The octads in $\calS(5,8,24)$ are given by$:$
\begin{enumerate}
\item[(3+5)] All three Romans together with a line of the plane {\rm (21} octads{\rm );}
\item[(2+6)] Two Romans with an oval in the plane {\rm (168 --);}
\item[(1+7)] One Roman, together with a subplane defined over $\bfF_2$ {\rm (360 --);}
\item[(0+8)] No Romans, and the symmetric difference of two lines of the plane {\rm (210 --)}. 
\end{enumerate}
Here an oval means five points of $\bfP^2(\bfF_4)$ on a non-singular conic plus a unique point
not lying on chords of the conic.  
Thus there are $759 (=21+168+360+210)$ octads.  
\end{proposition}
 
The vectors in $D_n$ diagrams involve the following vectors:

$\emptyset$: the zero vector;

$[P]$: $P$ is a point in the plane or $\I$, $\II$ or $\III$.  The vector has $-3$ at the position $P$ and $1$'s elsewhere;

$[\widehat{P}]$: it has $5$ at $P$ and 1's elsewhere.

$[L]$: $L$ is a line in the plane.  The vector has 2 at the five points on $L$ and at $\I$, $\II$, $\III$, and $0$'s elsewhere.

$[Q]$: $Q$ is an oval in the plane.  The vector has 2 at the six points on $Q$ and at two points from $\I$, 
$\II$, $\III$, and $0$'s elsewhere.  For example, $[Q_0]$ has 2 at $\II$, $\III$.

Consider the following vectors in $\Lambda$:
\[\emptyset,\ [\widehat{\I}], \ [\widehat{\II}], \ [\widehat{\III}].\]
The corresponding Leech roots
\begin{equation}\label{xyz}
\alpha_0 = (-1, 1, \emptyset),\ 
\alpha_1 = (2, 1, [\widehat{\I}]), \ \alpha_2 = (2, 1, [\widehat{\II}]),\ \alpha_3 = (2, 1, [\widehat{\III}])
\end{equation}
generate a root lattice $R_0$ isomorphic to $D_{4}$ (see Figure \ref{D4}).   
\begin{figure}[h]
\begin{center}
\scalebox{1.0}{\xy 
(-60,20)*{};
(0,-10)*{};
@={(-15,0),(-30,0),(-45,0),(-30,15)}@@{*{\bullet}};
(-15,0)*{};(-30,0)*{}**\dir{-};(-30,0)*{};(-45,0)*{}**\dir{-};
(-30,0)*{};(-30,15)*{}**\dir{-};
(-45,-4)*{\alpha_1};(-30,-4)*{\alpha_0};(-15,-4)*{\alpha_3};(-34,15)*{\alpha_2};
\endxy}
\end{center}
\caption{Root lattice $R_0$ of type $D_4$}\label{D4}
\end{figure}

\begin{lemma}\label{l2} {\rm (\cite[Lemma 3.2]{DK0})}
There are exactly forty-two Leech roots
which are orthogonal to $R_0$.  They are given by 
\[(1,1,[L]),\ (1,1,[P])\]
where $L$ is a line and $P$ is a point in the plane.
\end{lemma}

Among many sub-diagrams of type $D_4$ perpendicular to $R_0$, the standard one $R_1$ consists of $(1,1,[L_\infty])$ with the line of infinity $L_\infty$ and three points at infinity $\infty_1, \infty_\omega, \infty_{\bar\omega}$.
The following is immediate from the above.

\begin{figure}[h]
\begin{center}
\scalebox{1.0}{\xy 
(-60,20)*{};
(0,-10)*{};
@={(-15,0),(-30,0),(-45,0),(-30,15)}@@{*{\bullet}};
(-15,0)*{};(-30,0)*{}**\dir{-};(-30,0)*{};(-45,0)*{}**\dir{-};
(-30,0)*{};(-30,15)*{}**\dir{-};
(-45,-4)*{\infty_1};(-30,-4)*{L_\infty};(-14,-4)*{\infty_{\bar\omega}};(-29,18)*{\infty_\omega};
\endxy}
\end{center}
\caption{Root lattice $R_1$ of type $D_4$}\label{D4+}
\end{figure}

\begin{lemma}\label{l3} {\rm (\cite[Lemma 3.2, Case I]{KeKo})}
There are exactly twenty four Leech roots
which are orthogonal to $R_0 \oplus R_1$.  They are given by 8 lines and 16 points
\[(1,1,[x=a]), (1,1,[y=b]), \ (1,1,[(a, b)]),\]
where $a, b \in \bfF_4$ and $(a, b)$ are the points in the affine plane $\bfA^2(\bfF_4)=\bfP^2(\bfF_4) \setminus L_\infty$.
\end{lemma}

Now we consider a Leech root  $r$  such that $R_0$ and $r$  generate $D_5$.

\begin{lemma}\label{l2-2} {\rm (\cite[(3.9)]{DK0})}
The number of such Leech roots  $r$  is 168.
They are given by $(1,1,[C])$ for an oval  $C \subset \bfP^2(\bfF_4)$.
\end{lemma}

Finally we consider a Leech root  $r$  such that $R_0, R_1$ and $r$  generate $D_4 \oplus D_5$.
By Proposition~\ref{MOGoctads}, we have the following:

\begin{lemma}\label{l3-2} {\rm (\cite[Lemma 3.6, Case I]{KeKo})}
Such Leech roots correspond bijectively to the permutations $\sigma$ of  $\bfF_4 =\{0, 1, \omega, \bar\omega\}$.
In particular their number is twenty four.
\begin{enumerate}
\item If a permutation $\sigma$  is even, then $r$  is  $(1,1,[L_\sigma])$, where  $L_\sigma$  is the unique line passing through the four points  $(a, \sigma(a))$, $a \in \bfF_4$.
\item Otherwise, $r$  is  $(1,1,[C_\sigma])$, where  $C_\sigma$  is the oval consisting of two points of infinity of slope $0, \infty$ and 
the four points  $(a, \sigma(a))$.
\end{enumerate}
\end{lemma}

\begin{remark}
(1) The 42 $(-2)$ Leech roots in Lemma~\ref{l2} are realized as rational curves on the inseparable double cover $z^2 = xy(x^3+y^3+1)$ of $\bfP^2$, the model of superspecial $K3$ surface studied in \cite{DK0}.

(2) The inseparable double cover $z^2 = (x^4-x)(y^4-y)$ is the model realizing 24 $(-2)$ rational curves in Lemma~\ref{l3}.
Note that this surface has an extra $D_4$ singular point at $(\infty, \infty)$ and is the quotient of the product of two cuspidal rational curves $z^2=x^4-x$ by the diagonal action of $\mu_2$ (\cite{KS}).
\end{remark}

\subsection{Leech roots in terms of a fixed oval}
Sometimes it is more convenient to list Leech roots in relation with a fixed oval than 
$(x, y)$-coordinates.
Here we handle the diagrams containing a root lattice $D_n$ for various $n$.
To do this we specify MOG as follows.
We take an oval $Q_0$ consisting of the points $\{\infty, 0, 1, 2, 3, 4\}$ and label 
the remaining 15 points with Sylvester's {\it synthemes} such as 
$\infty 0.14.23$ (we do not distinguish 
$\infty 0.14.23$, $14.23.\infty 0$ etc.).  See Figure \ref{MOG2}.
\begin{figure}[ht]
 \centering
  \begin{tabular}{|c|c|c|}
   \hline  
$\infty 0.14.23$ \ \ $\infty 0.13.24$  & $\infty 4.03.12$ \ \ $\infty 2.01.34$ & $\infty 3.02.14$ \ \ $\infty 1.04.23$ \\ 
\quad \quad \I \ \quad \quad $\infty 0.12.34$  & $\infty 1.02.34$ \ \ $\infty 3.04.12$  &  $\infty 2.03.14$ \ \ $\infty 4.01.23$ \\    \hline 
\II \ \quad \quad \quad $\infty$  & $\infty 2.04.13$ \ \ $\infty 4.02.13$ & $1$\quad \quad \quad $3$ \\ 
\III \ \ \quad \quad \quad  $0$  & $\infty 3.01.24$ \ \ $\infty 1.03.24$ & $4$ \quad \quad \quad $2$ \\
\hline
\end{tabular}
\caption{Labels for positions in MOG (\cite[Figure 23.18]{CS})}\label{MOG2}
\end{figure}

Among 21 lines in the plane, 15 of them meet $Q_0$ in two points.  In this case, we label such a line
as the Sylvester's {\it duad} of these two points. For example, the duad $\infty 0$ meets $Q_0$ at $\infty$ and $0$, and
also contains the points $\infty 0.12.34$, $\infty 0.13.24$, $\infty 0.14.23$.
 The remaining 6 lines do not meet $Q_0$ which are the axes of $Q_0$.  They are labeled with {\it totals}.  Here a total is the set of five synthemes containing all 15 duads.  A total $a|bcdef$ being an abbreviation for the set of five synthemes $ad.ce.bf$, $ae.bc.df$, $af.be.cd$, $ab.cf.de$, $ac.bd.ef$, which are the points on the axes (see \cite[Figure 23.19]{CS}).
 For Sylvester's synthemes, duads and totals, we refer the reader to \cite[\S 10.4]{DK1}.

\subsection{The odd unimodular hyperbolic lattices and their even parts}\label{odd unimodular}
In this subsection we consider the odd unimodular hyperbolic lattice $I_{1, n-1}$  
of rank  $n \ge 3$, which is unique up to isomorphisms,
and review the orthogonal group of of the  sublattice  $I_{1, n-1}^{ev}$  
of elements of even norm in $I_{1, n-1}$.

First we recall that
a hyperbolic lattice is {\it reflective} if its orthogonal group contains a subgroup of finite index generated by reflections.

\begin{proposition}\label{reflectivity inherits}{\rm (\cite[\S2]{V2})}
If  an odd lattice $L$  is reflective, then so is  the even part $L^{ev}$ with the same Coxeter diagram as $L$.
\end{proposition}
\begin{proof}
Let ${\rm Disc}$  be the discriminant group of  $L^{ev}$ and $a_{odd} \in {\rm Disc}$ the element of order 2 corresponding to the overlattice  $L$.
The image of the natural injection  $O(L) \hookrightarrow O(L^{ev})$  consists of the orthogonal transformations preserving the element  $a_{odd} \in {\rm Disc}$.
Hence it is of finite index in  $O(L^{ev})$.
Therefore, just by replacing  $(-k)$-reflections  $R_v$, $k$ being odd, of $L$ with $(-4k)$-reflections  $R_{2v}$ of $L^{ev}$, we have our assertion.
\end{proof}

The lattice $I_{1, n-1}$ is reflective for  $n \le 18$  by Vinberg~\cite{V}  and for  $n = 19, 20$ by Vinberg-Kaplinskaja~\cite{VK}.
Hence we have

\begin{corollary}
The even lattice $I_{1, n-1}^{ev}$ is reflective for  $n \le 20$.
\end{corollary}
 
\begin{remark}
The converse of Proposition~\ref{reflectivity inherits} does not hold. 
In fact, though   $I_{1, 21}$  is not, its even part $I_{1,21}^{ev}$ is reflective (\cite[\S8]{Bor1}).
\end{remark}

The Vinberg case of $n \le 18$  is divided into two sub-cases:

\begin{description}
\item[Sub-case (F)] $n = 4,...,10, 12, 13,14, 16$
where the Coxeter diagram (see \cite[Table 4]{V}) contains a unique $(-1)$-vector.
For example, when $n=16$, the diagram is Figure~\ref{polygon2legs} (left) consisting of 17 $(-2)$-roots $\bullet$ and a $(-1)$-root  (or  $(-4)$-root)  $\circ$.

\item[Sub-case (ID)] 
$n = 3, 11, 15, 17, 18$ where the Coxeter diagram contains a unique pair of $(-1)$-vectors $u, v$  with  $(u, v) = 1$ (or two $\circ$ joined by a thick line in our diagram).
For example, when $n=15$, the diagram is Figure~\ref{polygon2legs} (right) consisting of 15 $(-2)$-roots $\bullet$ and two $(-1)$-roots $\circ$.
\end{description}
 (Note: Vinberg's diagrams for $n =11, 14, 15$ in \cite[Table 4, p.345--347]{V} contain misprints.  For correct ones, see Vinberg \cite[Table 5]{V1972}).
 
\begin{figure}[h]
\begin{center}
\xy
(60,45)*{};
(-30,-5)*{};
@={(0,40),(10,40),(20,40),(30,40),
(0,30),(10,30),(40,30),
(0,20),(40,20),
(0,10),(30,10),(40,10),
(0,0),(10,0),(20,0),(30,0),(40,0)
}@@{*{\bullet}};
(40,0)*{};(40,30)*{}**\dir{-};
(0,40)*{};(30,40)*{}**\dir{-};
(0,0)*{};(0,40)*{}**\dir{-};
(0,0)*{};(40,0)*{}**\dir{-};
(0,40)*{};(10,30)*{}**\dir{-};
(30,10)*{};(40,0)*{}**\dir{-};
(30,10)*{};(40,0)*{}**\dir{-};
(30,40)*{};(40,40)*{\circ}**\dir{=};
(40,30)*{};(40,40)*{}**\dir{=};
@={(50,40),(60,40),(70,40),
(50,30),(60,30),
(50,20),(90,20),
(50,10),(80,10),(90,10),
(50,0),(60,0),(70,0),(80,0),(90,0)
}@@{*{\bullet}};
(90,0)*{};(90,20)*{}**\dir{-};
(50,40)*{};(70,40)*{}**\dir{-};
(50,0)*{};(50,40)*{}**\dir{-};
(50,0)*{};(90,0)*{}**\dir{-};
(50,40)*{};(60,30)*{}**\dir{-};
(80,10)*{};(90,0)*{}**\dir{-};
(80,10)*{};(90,0)*{}**\dir{-};
(80,40)*{{\circ}};(90,30)*{{\circ}};
(70,40)*{};(80,40)*{}**\dir{=};
(90,20)*{};(90,30)*{}**\dir{=};
(80.9,39.8)*{};(89.9,30.8)*{}**\dir{=};
(80.8,39.7)*{};(89.8,30.7)*{}**\dir{=};
(80.7,39.6)*{};(89.7,30.6)*{}**\dir{=};
(80.6,39.5)*{};(89.6,30.5)*{}**\dir{=};
(80.5,39.4)*{};(89.5,30.4)*{}**\dir{=};
\endxy
\end{center}
\caption{Hexadecagon and pentadecagon with two legs in $I_{1, 15}$ and $I_{1, 14}$}
\label{polygon2legs}
\end{figure}

\begin{remark}
The five values of $n$ in the subcase (ID) correspond to the five $n\leq 18$ for which a Niemeier lattice of type $(D_{27-n} \oplus *)$  exists.
More explicitly such (five) Niemeier lattices are $N(D_{24}), N(D_{16}\oplus E_8), N(2D_{12}), N(D_{10}\oplus 2E_7)$  and  $N(D_{9}\oplus A_{15})$.
See also \cite[Table on p.339]{V} and \cite[Chap.28]{CS}.
\end{remark}

Let  $X_n$  be a $K3$ surface over  $\bfC$  whose Picard lattice is  $I_{1, n-1}^{ev}$.
In the sub-case (F), the lattice $I_{1, n-1}^{ev}$  is $(-2)$-reflective and the $K3$ surface  $X_n$  has only {\it finite} automorphisms by the Torelli-type theorem.
In the sub-case (ID), $\Aut(X_n)$ contains the {\it infinite dihedral} group generated by two $(-1)$-reflections as a subgroup of finite index.

Later we study the case  $n= 17, 18, 19$ in more details with the Leech lattice technique of \cite[Chap.28]{CS} and \cite{Bor1}.

\subsection{Mordell--Weil lattices}\label{Mordell}

Let $\pi: S \to C$ be an elliptic surface with a zero section $(O)$.
The {\it trivial lattice} $\Triv(\pi)$ is the sublattice of the N\'eron-Severi lattice $\NS(S)$ 
generated by $(O)$ and components of fibers of $\pi$.  Let $E$ be the generic fiber over $K=k(C)$. Let
$E(K)$ be the group of $K$-rational points on $E$ and $E(K)_{\tor}$ the torsion subgroup of $E(K)$.  Then the following theorem holds (cf. Shioda--Sch\"utt \cite[Lemma 6.18]{SS}).

\begin{theorem}\label{M-W}  There exists a group monomorphism of abelian groups
\[ \varphi: E(K)/E(K)_{\tor} \hookrightarrow \NS(S)\otimes \bfQ\]
such that the image $\Im(\varphi)$ is perpendicular to $\Triv(\pi)$.
In particular, the inversion acts on $\Im(\varphi)$ as $-1$. 
\end{theorem}

In this paper, we use this theorem as in the following.  
Assume that $\rank\ E(K) =1$ and 
let $\delta$ be a generator of $E(K)\otimes \bfQ$.  
Then the inversion of $E(K)$ sends $\delta$ to $-\delta$, 
but in general, the inversion does not act on $\NS(S)\otimes \bfR$ as a reflection with respect to $\delta$, that is, acts on $\Triv(\pi)$ non trivially.  However it works like a reflection, that is, 
it sends a half space $\la x, \delta \ra \ge 0$ in $\NS(S)\otimes \bfR$
to the opposite half space $\la x, \delta \ra \le 0$ (see Ujikawa \cite[Proposition 2.5]{Ujikawa}).

By the theorem, the inversions often give a geometric realization of orthogonal transformations of the form
$g = \sigma(R)\sigma(R')g'$
which are used in Borcherds~\cite[\S3]{Bor1}, when  $\NS(S)$  is isomorphic to the orthogonal complement of a root lattice $R$  in  $\II_{1,25}$.
This fact is not used in the sequel but a theoretical/historical background of our method.

\section{The automorphism group of the Kummer surface associated with a curve of genus 2 whose Jacobian has $p$-rank 2}\label{sec3}

Let $C$ be a curve of genus 2 whose Jacobian has $p$-rank 2.  Let $X$ be the minimal resolution of the Kummer quartic surface $S$, 
given by the equation (\ref{KummerQuarticEq}), associated with $C$.
It contains twenty $(-2)$-curves which are the proper transforms of four tropes
and the sixteen exceptional curves over the four rational double points of $S$ of type $D_4$.
The dual graph $\Gamma_{20}$ of twenty $(-2)$-curves is given as in Figure \ref{20curves}.

\begin{figure}[htbp]
\begin{center}
\scalebox{0.8}{
\xy
(70,65)*{};
(-15,-10)*{};
(0,0)*{};(0,30)*{}**\dir{-};
(00,00)*{};(30,00)*{}**\dir{-};
(0,30)*{};(30,30)*{}**\dir{-};
(30,0)*{};(30,30)*{}**\dir{-};
(40,10)*{};(30,0)*{}**\dir{-};
(40,10)*{};(50,20)*{}**\dir{-};
(00,00)*{};(20,20)*{}**\dir{-};
(20,20)*{};(50,20)*{}**\dir{-};
(20,20)*{};(20,50)*{}**\dir{-};
(50,20)*{};(50,50)*{}**\dir{-};
(30,30)*{};(50,50)*{}**\dir{-};
(00,30)*{};(20,50)*{}**\dir{-};
(50,50)*{};(20,50)*{}**\dir{-};
@={(0,0),(0,15),(0,30),(15,0),(15,30),(30,0),(30,30),(30,15),(40,10),(50,20),
(40,40),(50,50),(50,35),(10,40),(20,50),(35,50),(10,10),(20,20),(35,20),(20,35)}@@{*{\bullet}};
(-3,-3)*{E_{2}};(-4,33)*{T_{1}};(55,20)*{E_{1}};(34,-3)*{T_{4}};
(53,53)*{T_{2}};(34,28)*{E_{3}};(24,23)*{T_{3}};(17,54)*{E_{4}};
(15,-4)*{E_{42}};(45,10)*{E_{41}};(26,13)*{E_{43}};(14,8)*{E_{32}};(-4,13)*{E_{12}};
(6,43)*{E_{14}};(35,54)*{E_{24}};(54,35)*{E_{21}};(35,41)*{E_{23}};(25,36)*{E_{34}};
(13,27)*{E_{13}};(37,16)*{E_{31}};
\endxy 
}
\caption{The dual graph $\Gamma_{20}$ of 20 $(-2)$-curves on $X$}
\label{20curves}
\end{center}
\end{figure}

Here $\{T_1, T_2, T_3, T_4 \}$ are the proper transforms of tropes and 
$\{E_{i}, E_{ji}\}_{j\ne i}$ are the exceptional curves over the singular point $p_i$ 
of type $D_4$ $(1\leq i \leq 4)$.

\begin{remark}\label{Mukai-Peters}
Over the complex numbers, Peters--Stienstra \cite{PS} studied a 1-dimensional family of $K3$ surfaces
containing 20 $(-2)$-curves forming the dual graph in Figure \ref{20curves} and 
Mukai--Ohashi \cite{MO} also studied quartic surfaces given in (\ref{KummerQuarticEq}).
\end{remark}

In this section, we study the automorphism of Kummer surfaces under the following assumption.

\begin{assumption}\label{assum1}  
$\NS(J(C))$ has rank 1 (hence $\NS(X)$ has rank 17).  
We call such $X$ {\it Picard general}.  
\end{assumption}

Since $\NS(J(C))$ is naturally embedded into the endomorphism ring ${\rm End}(J(C))$,
we have an example of a Jacobian  $J(C)$ with  N\'eron-Severi rank 1 (in positive characteristic) from the following result.

\begin{proposition}{\rm (\cite{Koi}, \cite{Mori})}\label{Koizumi-Mori}
Let  ${\mathcal M}_g$  be the {\rm (}coarse{\rm )} moduli space of curves of genus $g$.
The generic curve  $C$  of genus  $g$  over the algebraic closure of the function field of  ${\mathcal M}_g$ is endo-general, that is, the endomorphism ring of its Jacobian is isomorphic to  $\bfZ$.
\end{proposition}

In our case  $g=2$, the base field has transcendence degree 3 over $\bfF_2$.

\begin{lemma}\label{(1,16)}  Let $L_{20}$ be the sublattice of $\NS(X)$ 
generated by $20$ $(-2)$-curves.  
Then,  
$L_{20}$ has the
signature $(1,16)$ and the determinant $16$.  
\end{lemma}
\begin{proof}
For the proof, see Peters--Stienstra \cite[Proposition 1]{PS}.
Also we can prove the assertion as follows.  
Here we may assume that $X$ is Picard general.
Let $\pi$ be 
the elliptic fibration defined by the linear system
\[|E_{12}+E_{13}+ E_{31}+E_{32}+
2(T_1+E_{14}+E_4+E_{34}+T_{3})|.\]
It has a reducible singular fiber of type $\tilde{D}_{8}$ and 
two sections $E_{1}$, $E_{3}$.  Moreover $E_{21}$, $T_2$, $E_{23}$, $T_4$, 
$E_{41}$, $E_{42}$, $E_{43}$ are components of 
reducible singular fibers
forming a root lattice $A_3\oplus D_4$.  Since $E_1\cdot E_{21}=E_1\cdot E_{41}=1$, $A_3$ and $D_4$ are supported on different fibers.
Since $\rank (\NS(X))=17$, 
these reducible fibers are of type $\tilde{D}_4$ and of 
type $\tilde{A}_3$.  Thus the Mordell-Weil group of $\pi$ is torsion.  
The components of fibers and 
a section $E_1$ generate a lattice isomorphic to 
$U\oplus D_8\oplus D_4\oplus A_3$ whose determinant is equal to $2^6$.  
By adding one more section $E_3$, these 18 $(-2)$-curves generate a lattice with determinant $2^4$.
The equality
\[
|E_{2}+E_{12}+ T_{1}+E_{14}+ E_{4}+E_{34}+T_{3}+E_{32}|
= |T_{4}+E_{43}+ E_{3}+E_{23}+T_2+E_{21}+E_1+E_{41}|
\]
implies that the remaining $(-2)$-curve $E_2$ is represented by a linear combination of the above 18 $(-2)$-curves.  The curve $E_{24}$ is similar.  Thus we have proved the assertion.
\end{proof}

\begin{lemma}\label{NeronSeveri}
Let $C$ be any ordinary curve of genus two and let $X$ be the associated Kummer surface.  Then $L_{20}$ is primitive in $\NS(X)$, that is, $\NS(X)/L_{20}$ is torsion free.  In particular, if $X$ is Picard-general, then
$L_{20}= \NS(X)$.
\end{lemma}
\begin{proof}
It suffices to show the first assertion.  
Recall that $\NS(X)$ contains a sublattice 
$\la 4 \ra \oplus 4D_4$ where $\la 4\ra$ is generated by the total transform $H$ of the
hyperplane section of the Kummer quartic surface $S$ and
$D_4$ is generated by the exceptional divisors over a 
rational double point on $S$ (Propositions \ref{KummerQuartic}).
Consider the two elliptic fibrations with two singular fibers of type 
$\tilde{D}_{6}$ defined by the complete linear
systems
\[
|E_{13}+E_{23} + E_{41}+E_{42}+2(E_3+E_{43}+T_4)|=|E_{14}+E_{24}+E_{31} + E_{32}+2(E_{4}+E_{34}+T_{3})|,
\]
\[|E_{12}+E_{13} + E_{24}+E_{34}+2(E_4+E_{14}+T_1)|=|E_{21}+E_{31}+E_{42} + E_{43}+2(E_{1}+E_{41}+T_{4})|.\]
For each fibration, the sum of 8 simple components of two singular fibers of type $\tilde{D}_{6}$ is divided by 2 in $\NS(X)$, and hence we have two non-zero vectors in the discriminant group $(4D_4)^*/4D_4 \cong (\bfZ/2\bfZ)^8$.  
We can easily see that two vectors are linearly
independent.  Thus, by adding these two vectors to $4D_4$, we obtain an overlattice $\overline{4D_4}$ of $4D_4$ which has determinant $2^4$.  
On the other hand, it is known that $4D_4\otimes \bfQ \cap \NS(X)$ 
has determinant $2^4$ (see Matsumoto \cite[Table 1]{Matsumoto}).  
Thus $\overline{4D_4}=4D_4\otimes \bfQ \cap \NS(X)$ is primitive in 
$\NS(X)$.  Moreover, $2T_1 - H$ 
is a linear combination of $E_i, E_{ij}$ with integral coefficients.  
By adding $T_1$ to $\la 4 \ra \oplus \overline{4D_4}$, 
we get an overlattice with the determinant $2^4$ which is nothing but $L_{20}$ (Lemma \ref{(1,16)}).  Now the primitiveness of $L_{20}$ in $\NS(X)$ follows 
from that of $\overline{4D_4}$ in $\NS(X)$.
\end{proof}

\begin{remark}
Let $A$ be an abelian surface in any characteristic $p\geq 0$ (in case $p=2$, we assume its $p$-rank is 2 or 1).
The argument in the above proof can be generalized as follows.
Let $\pi : A\to A/\la \iota\ra$ be the quotient morphism by the inversion $\iota$ and
let $\phi : \Km(A) \to A/\la \iota \ra$ be the minimal resolution.
Then $A/\la \iota \ra$ has sixteen $A_1$-singularities (resp. four $D_4$- or two $D_8$-) when $p\ne 2$ (resp. $p=2$, $A$ has $p$-rank 2 or $p$-rank 1).  
Let $N$ be a sublattice of $\NS(\Km(A))$ generated by the exceptional curves which is isomorphic to $16A_1$ (resp. $4D_4$ or $2D_8$). 
Let $\overline{N} = N\otimes \bfQ \cap \NS(\Km(A))$.  Then 
$[\overline{N}: N]=2^5$ (resp. $2^2$ or $2$) when $p\ne 2$ (resp. $p=2$, 
$A$ has $p$-rank 2 or $p$-rank 1), and the root sublattice of 
$\overline{N}$ is $N$ (Matsumoto \cite{Matsumoto}). The sublattice 
$\phi^*\pi_*(\NS(A))$ is isomorphic to $\NS(A)(2)$, the lattice 
obtained from $\NS(A)$ by multiplying the intersection form by 2.
Both sublattices $\NS(A)(2)$ and $\overline{N}$ are primitive in 
$\NS(\Km(A))$ and are mutually orthogonal.  
Thus $\NS(\Km(A))$ is obtained from
$\NS(A)(2)\oplus \overline{N}$ as its overlattice (see Nikulin \cite[Proposition 1.5.1]{Nikulin}).
\end{remark}

Next we will discuss a finite polyhedron in $\NS(X)_\bfR$ by restricting 
Conway's fundamental domain $\calC$.  
Let $R = D_4\oplus D_5$ be a root sublattice of $\II_{1,25}$ generated by the following Leech roots
\[\alpha_0, \ \alpha_1, \ \alpha_2,\ \alpha_3,\ \alpha_4 =(1,1,[y=x]), \ \alpha_5= (1,1,[\infty_1]),\] 
\[\alpha_6=(1,1,[L_\infty]),\ \alpha_7 = (1,1,[\infty_{\bar\omega}]),\ \alpha_8 =(1,1,[\infty_\omega]),\]
where $\alpha_0, \alpha_1, \alpha_2, \alpha_3$ are as in (\ref{xyz})
(see Figure \ref{D4D5} and \cite[Figure 23.20]{CS}):

\begin{figure}[h]
\begin{center}
\scalebox{1.0}{\xy 
(-50,25)*{};
(70,-10)*{};
@={(0,0),(15,0),(30,0),(45,0),(-15,0),(-30,0),(-45,0),(-30,15),(30,15)}@@{*{\bullet}};
(0,0)*{};(15,00)*{}**\dir{-};(30,0)*{};(30,15)*{}**\dir{-};
(15,00)*{};(30,0)*{}**\dir{-};(30,0)*{};(45,0)*{}**\dir{-};
(-15,0)*{};(-30,0)*{}**\dir{-};(-30,0)*{};(-45,0)*{}**\dir{-};
(-30,0)*{};(-30,15)*{}**\dir{-};
(-45,-4)*{\alpha_1};(-30,-4)*{\alpha_0};(-15,-4)*{\alpha_3};(0,-4)*{\alpha_4};(15,-4)*{\alpha_5};(30,-4)*{\alpha_6};(45,-4)*{\alpha_7};(-34,15)*{\alpha_2};
(26,15)*{\alpha_8};
\endxy}
\end{center}
\caption{Root lattice of type $D_4\oplus D_5$}\label{D4D5}
\end{figure}

We will show that the orthogonal complement $R^{\perp}$ of $R$ is isomorphic to $\NS(X)$.  Let $\calC$ be Conway's fundamental domain
and let $\calC(R^{\perp})$ be the restriction of $\calC$ to
$R^{\perp}\otimes \bfR$ under the embedding $R^{\perp} \subset \II_{1,25}$.
There are 24 Leech roots perpendicular to $R_0 \oplus R_1$ of type $D_4 \oplus D_4$ which are listed in Lemma \ref{l3}.
Among those, 8 lines are perpendicular to $R = R \cup \{\alpha_6: [y=x]\}$.
A Leech root $(1,1,[(a, b)])$ is perpendicular to $R$ if and only if $a \ne b$.
There are twelve such points and we have 
20 (= 8+12) Leech roots.
Obviously we can identify eight $(1,1,[L])$ with $T_i, E_i$ and
twelve Leech roots $(1,1,[P])$ with twelve $(-2)$-curves $E_{ij}$
in Figure \ref{20curves} as in Figure \ref{20curves2}.

\begin{figure}[htbp]
\begin{center}
\scalebox{0.8}{
\xy
(70,65)*{};
(-15,-13)*{};
(0,0)*{};(0,30)*{}**\dir{-};
(00,00)*{};(30,00)*{}**\dir{-};
(0,30)*{};(30,30)*{}**\dir{-};
(30,0)*{};(30,30)*{}**\dir{-};
(40,10)*{};(30,0)*{}**\dir{-};
(40,10)*{};(50,20)*{}**\dir{-};
(00,00)*{};(20,20)*{}**\dir{-};
(20,20)*{};(50,20)*{}**\dir{-};
(20,20)*{};(20,50)*{}**\dir{-};
(50,20)*{};(50,50)*{}**\dir{-};
(30,30)*{};(50,50)*{}**\dir{-};
(00,30)*{};(20,50)*{}**\dir{-};
(50,50)*{};(20,50)*{}**\dir{-};
@={(0,0),(0,15),(0,30),(15,0),(15,30),(30,0),(30,30),(30,15),(40,10),(50,20),
(40,40),(50,50),(50,35),(10,40),(20,50),(35,50),(10,10),(20,20),(35,20),(20,35)}@@{*{\bullet}};
(-6,0)*{y=1};
(-6,30)*{x=0};(56,20)*{y=0};(34,-3)*{x=\bar\omega};
(53,53)*{x=1};(36,28)*{y=\omega};(25,23)*{x=\omega};
(12,52)*{y=\bar\omega};
\endxy 
}
\caption{20 Leech roots}
\label{20curves2}
\end{center}
\end{figure}

Since both $\NS(X)$ and $R^{\perp}$ contain 
twenty $(-2)$-vectors as their generators (Lemma \ref{(1,16)}), 
we have $\NS(X)\cong R^{\perp}$.  In particular, we have the following Lemma.
\begin{lemma}\label{PicardLattice}
$\NS(X)\cong U\oplus E_8 \oplus D_4\oplus A_3$.
\end{lemma}
\begin{proof}
Note that $\II_{1,25} \cong U\oplus 3E_8$ and
$D_4^{\perp}$ in $E_8$ (resp. $D_5^{\perp}$ in $E_8$) is isomorphic to $D_4$ (resp. $A_3$).
The assertion now follows.
\end{proof}

In the following we identify twenty $(-2)$-curves and twenty Leech roots, and use the notation in Figure
\ref{20curves}.  Let $\calC(X) =\calC(R^{\perp})$.
These twenty $(-2)$-curves define 20 faces of $\calC(X)$.  
The remaining faces are defined by the projections of Leech roots $r$ such that $r$ and $R$ generate a negative definite lattice, that is, a root lattice.  
There are four possibilities of $\la r, R\ra$: 
\[(a)\ D_4\oplus D_5 \oplus A_1,\ \ (b) \ D_4 \oplus D_6,\ \ (c)\ D_5\oplus D_5,\ \ (d)\ D_4\oplus E_6.\] 

\begin{lemma}\label{boundaries}
There are exactly $20$, $4$, $6$ or $8$ Leech roots of type $(a)$, $(b)$, $(c)$ or $(d)$, respectively.
The projection of a Leech root into $\NS(X)\otimes \bfQ$ has the norm $-2, -1, -1, -3/4$ respectively.
\end{lemma}
\begin{proof}
The case $(a)$ is nothing but the above 20 Leech roots. 
In case $(b)$, the desired Leech roots are
$(1,1,[P])$ such that $P$ lies on the line $y=x$ 
except $\infty$, and there are four such points:
$(0,0), (1, 1), (\omega, \omega), (\bar\omega, \bar\omega)$.
In case $(c)$, the desired Leech roots are $r=(1,1,[C_\sigma])$  with an odd permutation $\sigma$  with two fixed points, by Lemma~\ref{l3-2}.
Therefore, 
we have $6$ such Leech roots.

Finally in case $(d)$, the desired Leech roots are $(1,1,[L])$ where a line $L$
with slope $\omega, \bar\omega$.
There are eight such lines: $y=\omega x+ a$  and  $y=\bar\omega x+ a$  with $a \in \bfF_4$.
\end{proof}

Let $\delta$ be the projection of a Leech root of type $(a), (b), (c), (d)$.
Then $\delta \in \NS(X)$ for the case $(a)$, $2\delta \in \NS(X)$
for the cases $(b), (c)$ and $4\delta \in \NS(X)$ for the case $(d)$.
For each $\delta$, the hyperplane perpendicular to $\delta$ is a face of $\calC(X)$, and
hence $\calC(X)$ has $38=20+4+6+8$ faces by Lemma \ref{boundaries}.  We can easily determine
the incidence relation between $\delta$ and twenty $(-2)$-curves, by using the list of Leech roots in
the proof of Lemma \ref{boundaries}, as follows.

In case $(b)$, $\delta$ meets exactly two $(-2)$-curves among 20 $(-2)$-curves as follows:
\[\{T_i, E_i\}\ (i=1,2,3,4).\]
For example, if $\delta \cdot T_1=\delta\cdot E_1=1$, then 
\[2\delta = H_4 -4E_1-2(E_{21}+E_{31}+E_{41}),\]
where $H_4$ is the total transform of the hyperplane section of $S$.

In case $(c)$, $\delta$ meets exactly two $(-2)$-curves among 20 $(-2)$-curves as follows:
\[\{ E_{ij}, E_{ji}\} \ (1\leq i < j\leq 4).\]
For example, if $\delta \cdot E_{12}=\delta\cdot E_{21} =1$, then 
\[2\delta = H_4 -(2E_1+E_{21}+E_{31}+E_{41}) - (2E_2 + E_{12}+E_{32}+E_{42}) - E_{12}-E_{21}.\]

In case $(d)$, $\delta$ meets exactly three $(-2)$-curves among 20 $(-2)$-curves as follows:
\[\{E_{12}, E_{23}, E_{31}\},\ \{E_{12}, E_{24}, E_{41}\},\ \{E_{13}, E_{32}, E_{21}\},\ 
\{E_{13},E_{34}, E_{41}\}\]
\[\{E_{14}, E_{42}, E_{21}\},\ \{E_{14}, E_{43}, E_{31}\},\ \{E_{23}, E_{34}, E_{42}\}, \ 
\{E_{24}, E_{43}, E_{32} \}.\]
For example, if $\delta \cdot E_{12}=\delta\cdot E_{23}=\delta\cdot E_{31}=1$, then 
\[4\delta = 3H_4 -2((2E_2 + E_{12}+E_{32}+E_{42}) + (2E_3+E_{13}+E_{23}+E_{43})\]
\[+ (2E_1+E_{21}+E_{31}+E_{41}) + E_{12}+E_{23}+E_{31}).\]

Each $\delta$ of type $(b)$, $(c)$ and $(d)$ is not effective, but 
we show that it corresponds to an involution of $X$ by using 
Theorem \ref{M-W}.  More precisely, we take a suitable elliptic fibration on $X$ with a section whose Mordell-Weil group has rank 1 and is generated by $\delta$ over $\bfQ$.  Then the inversion is the desired 
involution.  We don't care whether its action on the hyperplane $\delta^\perp$ trivially or non-trivially.  The important thing is that it interchanges
the half-spaces adjacent along $\delta^\perp$ (see \S \ref{Mordell}).  This is the same as for the involutions considered in the remainder of
this paper.

\begin{lemma}\label{proj}
For each $\delta$ of type $(b)$, there exists an involution $\iota_{b,\delta}$ of $X$ which acts on $\NS(X)$ as the reflection associated with $\delta$.
In particular, $\iota_{b,\delta}^*(w) = w + a\delta \ (a > 0)$, where $w$ is
the projection of the Weyl vector $\rho$.
\end{lemma}
\begin{proof}
Let $\delta$ be as above: $\delta\cdot T_1=\delta\cdot E_1=1$.
Consider an elliptic fibration $\pi_1$ defined by a complete linear system
\[|E_{13}+E_{23} + E_{41}+E_{42}+2(E_3+E_{43}+T_4)|=|E_{14}+E_{24}+E_{31} + E_{32}+2(E_{4}+E_{34}+T_{3})|\]
with two singular fibers of type $\tilde{D}_{6}$.  We choose $E_2$ as a zero section $(O)$.
Note that $E_{12}, E_{21}$ are components of some 
reducible singular fibers of $\pi_1$.  
Since $\rank (\NS(X))=17$, 
the possible singular fibers are two of type $\tilde{A}_1$ or one of
type $\tilde{A}_3$.  If it is of type $\tilde{A}_3$, then the 
Shioda-Tate formula implies that the Mordell-Weil group is torsion.
However $\pi_1$ has four sections $E_1, E_2, T_1, T_2$ and hence
\[
16 =\det (\NS(X)) = 4^3/4^2= 4,
\]
which is a contradiction.  Therefore $\pi_1$ has two singular fibers of type $\tilde{A}_{1}$.    
Note that $\delta$ meets exactly two curves $T_1, E_1$ among 20 $(-2)$-curves and is perpendicular to a fiber and the section $E_2$.
Again the Shioda-Tate formula implies that the Mordell-Weil group of $\pi_1$ has 
rank at most 1.  
Together with the fact $\delta^2=-1$, $2\delta$ can not be represented as a linear combination of $(-2)$-curves in fibers of $\pi_1$. 
Thus the Mordell-Weil rank is equal to 1 and  
$\delta$ is perpendicular to the trivial lattice of $\pi_1$. 
Let $\iota_{b,\delta}$ be the involution in $\Aut(X)$ induced from the inversion of the Mordell-Weil group.  
Recall that the inversion acts trivially on the fibers of type $\tilde{A}_1$, 
$\tilde{D}_{2n}$, $\tilde{E}_7$ and $\tilde{E}_8$ 
(cf. Ujikawa \cite[\S 2]{Ujikawa}).  
The assertion now follows from  Theorem \ref{M-W}. 
\end{proof}

\begin{remark}
For each $\delta$ of type $(b)$, 
a projection from a double point of the Kummer quartic surface $S$ induces an involution 
$\iota$ of $X$ such that $\iota^*(\delta)=-\delta$.
\end{remark}

\begin{lemma}\label{transl-c}
For $\delta$ of type $(c)$, 
there exists an involution $\iota_{c,\delta}$ of $X$
which acts on $\NS(X)$ as the reflection associated with $\delta$.
In particular, $\iota_{c,\delta}^*(w) = w + a\delta \ (a > 0)$. 
\end{lemma}
\begin{proof}
Let $\delta$ be as above: $\delta \cdot E_{12} =\delta\cdot E_{21}=1$.
Consider an elliptic fibration $\pi_2$ defined by 
a complete linear system
\[|T_{1}+E_{13}+E_3+E_{23} + T_{2}+E_{24}+E_4+E_{14}|=|E_{1}+E_{31}+T_3+E_{32} + E_{2}+E_{42}+T_4+E_{41}|\]
with two singular fibers of type $\tilde{A}_{7}$ and choose $E_{34}$ as a zero section $(O)$. 
As in the proof of Lemma \ref{proj}, 
 $\pi_2$ has Mordell-Weil rank 1 and 
$\delta$ is perpendicular to the trivial lattice of $\pi_2$.
Let $\iota_{c,\delta}$ be the involution in $\Aut(X)$ induced from the inversion of the Mordell-Weil group.  The assertion now follows from  the same proof as that of Lemma \ref{proj}. 
\end{proof}

\begin{lemma}\label{transl-d}
For $\delta$ of type $(d)$, 
there exists an involution $\iota_{d,\delta}$ of $X$
satisfying $\iota_{d,\delta}^*(w) = w + a\delta \ (a > 0)$. 
\end{lemma}
\begin{proof}
Let $\delta$ be as above: $\delta \cdot E_{12} =\delta\cdot E_{23}= \delta\cdot E_{31} = 1$.  Consider an elliptic fibration $\pi_3$ defined by a complete linear system
\[|T_{1}+E_{13}+E_3+E_{43} + T_{4}+E_{42}+E_2+E_{32}+ T_{3}+E_{34}+E_4+E_{14}|\]
with a singular fiber of type $\tilde{A}_{11}$.
Note that $E_1, E_{21}, T_2$ are contained in a reducible fiber.  Since
$\rank (\NS(X))=17$, this fiber is of type $\tilde{A}_3$.  
Consider two sections $E_{24}$ and $E_{41}$ of the fibration.  Let $\iota_{d,\delta}$ be the involution in $\Aut(X)$ induced from the inversion
with respect to the bi-section $E_{24}+E_{41}$, which sends $\delta$ to 
$-\delta$.  Then $\iota_{d,\delta}$ preserves $17$ $(-2)$-curves in Figure \ref{20curves} perpendicular to $\delta$, and hence it fixes the projection of $w$ into $\delta^{\perp}$ in $\NS(X)$. The assertion now follows.
\end{proof}

\begin{remark}  
In the proof of Theorem \ref{MainCurve}, (2), we see that the involution $\iota_{d,\delta}$ is the covering involution of a double covering 
$\Phi_{|\Sigma_\delta|}: X \to \bfP^2$.
\end{remark}

\begin{lemma}\label{numtrivial}
Let $\Aut(\Gamma_{20})$ be the symmetry group of the dual graph 
$\Gamma_{20}$.  Then the natural map 
$\Aut(X)\to \Aut(\Gamma_{20})$ is injective.  In particular, the natural map $\Aut(X)\to O(\NS(X))$ is so.
\end{lemma}
\begin{proof}
Let $g \in \Ker(\Aut(X)\to \Aut(\Gamma_{20}))$.  Note that $g$ preserves
each $(-2)$-curve in $\Gamma_{20}$ and any elliptic fibration defined 
by some $(-2)$-curves in $\Gamma_{20}$.
Consider an elliptic fibration $\pi$ defined by the linear system
\[|2(T_1+E_{12}+E_2)+E_{13}+E_{14}+E_{32}+E_{42}| = |2(T_2+E_{21}+E_1)+E_{23}+E_{24}+E_{31}+E_{41}|.\]
It has four sections $T_3, T_4, E_3, E_4$ and has two reducible singular fibers of type 
$\tilde{D}_6$.  Moreover $E_{34}, E_{43}$ are components of some reducible singular fibers, 
and hence it has at least three reducible fibers.  Therefore $g$ acts trivially on the base of $\pi$
and acts on a general fiber of $\pi$ as an automorphism.  Since $g$ fixes the four points on the general
fiber which are the intersection with the four sections, $g$ acts trivially on the general fiber.  Therefore
$g$ is the identity.  The last assertion follows from Lemma \ref{NeronSeveri}.
\end{proof}

Now we study automorphisms preserving the polyhedron $\calC(X)$. 
The symmetry group $\Aut(\calC(X))$ coincides with $\Aut(\Gamma_{20})$.
 Obviously $\Aut(\Gamma_{20})$ is isomorphic to $\bfZ/2\bfZ\rtimes \mathfrak{S}_4$ where $\bfZ/2\bfZ$ is generated 
by the involution interchanging antipodal vertices and 
$\mathfrak{S}_4$ is the hexahedral group.
It is easy to see that the birational transformations of 
the Kummer quartic surface $S$ given in (\ref{CremonaOrdinary}),
(\ref{translation}) generate $(\bfZ/2\bfZ)^3$.
In the equation (\ref{KummerQuarticEq}) of $S$, if $a=b=c$, then involutions given in (\ref{transformation}) with the above 
$(\bfZ/2\bfZ)^3$ generate $\Aut(\Gamma_{20})$.  
It now follows from Lemma \ref{numtrivial} that
any element in $\mathfrak{S}_4$ preserves
the class of $H_4$, that is, if it is realized as an automorphism of $X$,
then it is induced from a projective transformation of $S$.
Thus if $a, b, c$ are pairwise different, 
$\Aut(X)\cap  \Aut(\Gamma_{20}) = (\bfZ/2\bfZ)^3$.
This is the case when $X$ is Picard general:

\begin{lemma}
If $X$ is Picard general, then $a, b, c$ are pairwise distinct.
\end{lemma}
\begin{proof}
If $a=b$, then the plane section $x+y = 0$ is defined by the equation 
$c(x^2+zw)^2 + x^2zw = 0$ and hence a union of two conics passing 
through $(0, 0,1, 0)$ and $(0, 0, 0, 1)$.
Their strict transform on $X$ is the union of two disjoint $\bfP^1$'s.
Note that the situation is the same as that of \cite{MO} as far as considering them in the Picard lattice.
We denote the sum of their divisor classes by  $F_{12}$.
Then the computation in  \cite[\S1]{MO} works the same and we have the following intersection numbers:
\[
(T_k+ E_k)\cdot F_{12} = 0, \quad (E_{kl}+E_{lk}) \cdot F_{12} = \left\{\begin{alignedat}{2} 4 & \ \{k, l\} = \{1, 2\} \\ 0 & {\rm \ otherwise.} \end{alignedat}\right.
\]
We have that $F_{12}$  is linearly independent from $L_{20}$ by considering invariant by Cremona involution.
( In the notation of \cite{MO},  $T_k+ E_k$ and $F_{12}$ belongs to 10A and 6B, which are of rank 7 and 3, respectively, and linearly independent.)\end{proof}

Now combining this lemma 
with Lemmas \ref{numtrivial}, \ref{proj}, \ref{transl-c}, \ref{transl-d},  
we have the following theorem.

\begin{theorem}\label{MainCurve}
Let $X$ be a Picard general Kummer surface associated with an ordinary curve of genus 2.  
Let $G$ be the (infinite) group generated by six involutions $\iota_{b,\delta}$, four involutions $\iota_{c,\delta}$ 
and eight involutions $\iota_{d,\delta}$.  Then we have
\begin{enumerate}
\item   $\Aut(X) \cong G \rtimes (\bfZ/2\bfZ)^3$, and
\item the eighteen involutions $\iota_{b,\delta}$, $\iota_{c,\delta}$ and  $\iota_{d,\delta}$ depend only on $\delta$ $($not on the choice of elliptic pencils$)$.
\end{enumerate}
\end{theorem}
\begin{proof}
(1) The proof is similar to the one for Lemma 7.3 in \cite{Kondo2} or Theorem 12.28 in \cite{Kondo3}.  

(2) Among the 20 divisors, we consider those with intersection number 0 with $\delta$.

In the case (c) there are 18 such divisors and their dual graph contains exactly two octagons, which are disjoint (Figure~\ref{fig:two octagons}).
Hence the choice of an elliptic pencil in the proof of Lemma~\ref{transl-c}
of $2\tilde{A}_{7}$-type is unique.
So we have nothing to prove.
\begin{figure}[htbp]
\begin{center}
\scalebox{0.8}{
\xy
(70,65)*{};
(-15,-13)*{};
(0,0)*{};(0,30)*{}**\dir{-};
(00,00)*{};(30,00)*{}**\dir{-};
(0,30)*{};(30,30)*{}**\dir{-};
(40,10)*{};(30,0)*{}**\dir{-};
(40,10)*{};(50,20)*{}**\dir{-};
(00,00)*{};(20,20)*{}**\dir{-};
(20,20)*{};(50,20)*{}**\dir{-};
(50,20)*{};(50,50)*{}**\dir{-};
(30,30)*{};(50,50)*{}**\dir{-};
(00,30)*{};(20,50)*{}**\dir{-};
(50,50)*{};(20,50)*{}**\dir{-};
@={(0,0),(0,15),(0,30),(15,0),(15,30),(30,0),(30,30),(40,10),(50,20),
(40,40),(50,50),(50,35),(10,40),(20,50),(35,50),(10,10),(20,20),(35,20)}@@{*{\bullet}}
\endxy 
}
\caption{Two octagons in the case (c)}
\label{fig:two octagons}
\end{center}
\end{figure}

In the case (d), there are 17 such divisors.
Their dual graph has two terminal vertices which are joined by three $A_5$-chains (Figure~\ref{fig:two ends+three chains}).
Hence their sum $\Sigma_\delta$ is nef and has self intersection number 2.
The complete linear system gives a double covering 
$\Phi_{|\Sigma_\delta|}: X \to \bfP^2$,
and the branch sextic has three points $P, Q, R$ which are simple singularities of type $a_5$.
The three points $P, Q, R$ are collinear and 
$\Sigma_\delta$ is the pull-back of the line passing through them.
There are three elliptic pencils with singular fibers of type $\tilde{A}_{11}, \tilde{A}_3$  which are used to prove Lemma~\ref{transl-d}.
But they are all pull-back from one of the three linear pencils passing through an $a_5$-point.    
Therefore, the covering involution of  $X \to \bar \bfP^2$  induces the inversion on all elliptic pencils corresponding to  $\delta$.
(Note that we choose a $0$-section among two sections perpendicular to 
$\delta$.
Then another section is a 2-torsion section 
and hence the inversion is independent on the choice of $0$-sections.)

\begin{figure}[h]
\begin{center}
\scalebox{0.8}{\xy 
(90,15)*{};
(-70,-15)*{};
@={(0,0),(15,0),(30,0),(45,0),(-15,0),(-30,0),(-45,0)}@@{*{\bullet}};
@={(0,10),(15,10),(30,10),(-15,10),(-30,10)}@@{*{\bullet}}; 
@={(0,-10),(15,-10),(30,-10),(-15,-10),(-30,-10)}@@{*{\bullet}};
(-30,10)*{};(30,10)*{}**\dir{-};(-30,-10)*{};(30,-10)*{}**\dir{-};
(-45,0)*{};(45,0)*{}**\dir{-};(-45,0)*{};(-30,10)*{}**\dir{-};
(45,0)*{};(30,10)*{}**\dir{-};
(45,0)*{};(30,-10)*{}**\dir{-};(-45,0)*{};(-30,-10)*{}**\dir{-};
\endxy}
\end{center}
\caption{Two ends joined by three chains in the case (d)}\label{fig:two ends+three chains}
\end{figure}

In the case (b) there are 18 such divisors and their dual graph is the dodecagon $\tilde A_{11}$ with equi-distributed 6 legs (Figure~\ref{fig:dodecagon}).
Instead of their sum we consider the sum $\Sigma_\delta$ taken with weight 2 on $\tilde A_{11}$  and 1 on 6 legs.
Then $\Sigma_\delta$ is nef and has self intersection number 12.
The morphism $\Phi_{|\Sigma_\delta|}: X \to \bfP^6$, given by the complete linear system, is the double cover of a sextic del Pezzo surface $dP_6$ in  
$\bfP^6$.  An elliptic pencil of $2\tilde D_6$ to prove Lemma~\ref{proj} 
is the pull-back of one of three conic pencils of $dP_6$.
(In fact, $\Sigma_\delta$ is the sum of three such elliptic divisors.)
Hence the covering involution of  $X \to dP_6$  induces the inversion on all elliptic pencils corresponding to  $\delta$.
(This involution also coincides with the covering involution of the projection from a $D_4$-singular point.)

\begin{figure}[ht]
\begin{center}
\unitlength1.5pt
\def\r{\color{red}\circle*{3}}
\def\b{\circle*{3}}
\begin{picture}(60,60)(-20,0)
\put(0,60)\b  
\put(28,52)\b  
\put(51,31)\b  
\put(60,0)\b 
\put(51,-31)\b 
\put(28,-52)\b 
\put(0,-60)\b  
\put(-28,-52)\b  
\put(-51,-31)\b  
\put(-60,0)\b  
\put(-51,31)\b 
 \put(-28,52)\b 

\put(0,32)\b 
\put(26,16)\b
\put(26,-16)\b  
\put(0,-30)\b  
\put(-26,-16)\b  
\put(-26,16)\b  

\qbezier(0,60)(55,55)(60,0) 
\qbezier(0,-60)(55,-55)(60,0) 
\qbezier(0,-60)(-55,-55)(-60,0) 
\qbezier(0,60)(-55,55)(-60,0)

\qbezier(0,60)(0,45)(0,30)
\qbezier(51,31)(38.5,23.5)(26,16) 
\qbezier(51,-31)(38.5,-23.5)(26,-16) 
\qbezier(0,-60)(0,-45)(0,-30) 
\qbezier(-51,-31)(-38.5,-23.5)(-26,-16) 
\qbezier(-51,31)(-38.5,23.5)(-26,16) 
\end{picture}
\end{center}
\vskip90pt
\caption{Dodecagon with 6 legs in the case (b)}
\label{fig:dodecagon}
\end{figure}
\end{proof}

\begin{remark}\label{hyperelliptic}
As we saw in the proof of (2) in the theorem, the involutions  $\iota_{d,\delta}$ and  
$\iota_{b,\delta}$  are the bi-elliptic and hyperelliptic involutions of the complete linear systems  $|\Sigma_\delta|$, determined canonically from $\delta$, respectively.
This holds true for $\iota_{c,\delta}$ also by taking the sum of 18 divisors in Figure~\ref{fig:two octagons} as $\Sigma_\delta$.
More precisely, $\Phi_{|\Sigma_\delta|}: X \to \bfP^3$  is a double cover of a quadric $Q$,
and  $\iota_{c,\delta}$  is the covering involution of  $X \to Q$.
Therefore,  in all cases, required involutions are determined only by  $\delta$.
\end{remark}

\begin{remark}\label{}
Elliptic fibrations $\pi_1, \pi_2, \pi_3$ considered in Lemmas
\ref{proj}, \ref{transl-c}, \ref{transl-d} have the Mordell-Weil rank 1.  This implies that $\Aut(X)$ is infinite.
\end{remark}

\begin{remark}\label{16-6}
Over the complex numbers, the Kummer surface contains two sets 
$\{N_\alpha\}, \{T_\beta \}$, called $(16_6)$-configuration, such that
$\{ N_\alpha\}$ and $\{T_\beta\}$ consist of disjoint sixteen $(-2)$-curves and each member of one set
meets exactly six member of other set.  The set $\{N_\alpha\}$ consists of the exceptional curves over
sixteen nodes of the Kummer quartic surface and $\{T_\beta \}$ the proper transforms of sixteen tropes.
The finite polyhedron $\calC(X)$ consists of 
32 $(-2)$-faces, 32 $(-1)$-faces, 60 $(-1)$-faces and 192 $(-3/4)$-faces.
First 32 faces are defined by $\{N_\alpha\}, \{T_\beta \}$, the second ones by
$H_4 -2N_\alpha$ and their images under the dual map, the third ones by $H_4 - (N_{\alpha_1}+N_{\alpha_2}+N_{\alpha_3}+
N_{\alpha_4})$ where $\{n_{\alpha_1}, n_{\alpha_2}, n_{\alpha_3}, n_{\alpha_4}\}$ is a G\"opel tetrad, and
the fourth ones by $3H_4 - 2\sum_{\alpha \in W} N_\alpha$ where $W$ is a Weber hexad.  For G\"opel tetrads
and Weber hexads, we refer the reader to Hudson \cite[\S 50--52]{H}, Kond\=o \cite[\S 12.3]{Kondo3}.
In our case, there are no disjoint sixteen $(-2)$-curves on $X$ (Shepherd-Barron \cite[Corollary 3.2]{SB}), in particular, no $(16_6)$-configurations. 
However if we allow reducible effective divisors, we have a $(16_6)$-configuration as follows:
\[ \{E_{ij}, \ 2E_i +\sum_{j\ne i} E_{ji}\}, \ \{T_i, \ T_i + E_{ij_1}+E_{j_1}+E_{ij_2}+E_{j_2}, j_1, j_2\ne i\}.\]
Then the case (b) $H_4 - 2(2E_i+\sum_{j\ne i} E_{ji})$, 
the case (c) $H_4 - ((2E_{i_1}+\sum_{j\ne i_1} E_{ji_1}) + (2E_{i_2}+\sum_{j\ne i_2} E_{ji_2}) + E_{i_1i_2}+E_{i_2i_1})$, the case (d)
$3H_4 - 2((2E_{i_1}+\sum_{j\ne i_1} E_{ji_1})+(2E_{i_2}+\sum_{j\ne i_2} E_{ji_2}) + (2E_{i_3}+\sum_{j\ne i_3} E_{ji_3}) + E_{i_1i_2}+E_{i_2i_3}+E_{i_3i_1})$
correspond to a node, a G\"opel  tetrad, a Weber hexad, respectively, as in the case over the complex numbers mentioned above..
\end{remark}

\section{The Kummer surface associated with a curve whose Jacobian has $p$-rank 1}\label{sec4}

Let $C_1$ be a general curve of genus two whose Jacobian $J(C_1)$ has 
$p$-rank 1.
In this case, the 2-torsion group of $J(C_1)$ is $\bfZ/2\bfZ$ and the Kummer quartic surface
$S_1$ has two rational double points of type $D_8$ and two tropes.

Let $X_1$ be the minimal resolution of singularities of $S_1$.
By taking the resolution of singularities,
we can see that $X_1$ contains 18 $(-2)$-curves corresponding to 18 black vertices except $D_2$, $D_2'$
in Figure \ref{p-rank1dualgraph}.

\begin{figure}[h]
\begin{center}
\scalebox{0.8}{\xy 
(80,55)*{};
(-55,-55)*{};
@={(40,0),(0,-40),(-40,0),(0,40),(27,27),(27,-27),(-27,27),(-27,-27),(15,35),(15,-35),(-15,35),(-15,-35),(35,15),(35,-15),(-35,15),(-35,-15),(-25,0),(25,0),(0,30),(0,-30)}@@{*{\bullet}};
(15,15)*{{\circ}};(-15,-15)*{{\circ}};
(40,0)*{};(35,15)*{}**\dir{-};(35,15)*{};(27,27)*{}**\dir{-};(27,27)*{};(15,35)*{}**\dir{-};
(15,35)*{},(0,40)*{}**\dir{-};
(0,40)*{};(-15,35)*{}**\dir{-};(-15,35)*{};(-27,27)*{}**\dir{-};
(-27,27)*{};(-35,15)*{}**\dir{-};(-35,15)*{};(-40,0)*{}**\dir{-};
(-40,0)*{};(-35,-15)*{}**\dir{-};(-35,-15)*{};(-27,-27)*{}**\dir{-};(-27,-27)*{};(-15,-35)*{}**\dir{-};
(-15,-35)*{},(0,-40)*{}**\dir{-};
(0,-40)*{};(15,-35)*{}**\dir{-};(15,-35)*{};(27,-27)*{}**\dir{-};
(27,-27)*{};(35,-15)*{}**\dir{-};(35,-15)*{};(40,0)*{}**\dir{-};
(40,0)*{};(25,0)*{}**\dir{-};(-40,0)*{};(-25,0)*{}**\dir{-};
(0,40)*{};(15.0,15.5)*{}**\dir{-};(0,-40)*{};(-15.0,-15.5)*{}**\dir{-};
(14.5,14.5)*{};(-14.5,-14.5)*{}**\dir{-};
(15,35)*{};(0,30)*{}**\dir{-};(-15,35)*{};(0,30)*{}**\dir{-};
(15,-35)*{};(0,-30)*{}**\dir{-};(-15,-35)*{};(0,-30)*{}**\dir{-};
(0,30)*{};(15.0,15.5)*{}**\dir{--};(0,-30)*{};(-15.0,-15.5)*{}**\dir{--};
(0,30)*{};(0,-30)*{}**\dir{--};
(44,0)*{E_{6}'};(-44,0)*{E_{6}};(0,-44)*{E_2'};(0,44)*{E_2};
(30,30)*{T_{1}};(30,-30)*{E_{4}'};(-30,30)*{E_4};(-30,-30)*{T_{2}};
(18,38)*{E_{1}};(18,-38)*{E_{3}'};(-18,38)*{E_{3}};(-18,-38)*{E_{1}'};
(38,18)*{E_{7}'};(38,-18)*{E_{5}'};(-38,18)*{E_{5}};(-38,-18)*{E_7};
(20,15)*{\delta_1};(-20,-15)*{\delta_2};
(25,4)*{E_{8}'};(-25,4)*{E_{8}};
(-5,27)*{D_2};(5,-28)*{D_2'};
(3,-10)*{4};(5,20)*{-1};(-5,-20)*{-1};
\endxy}
\end{center}
\caption{The dual graph of $(-2)$-curves ($p$-rank $1$ case)}
\label{p-rank1dualgraph}
\end{figure}

\noindent
Here we may assume that $T_1$ and $T_2$ are the proper transforms of two tropes.  Then 
$E_1, \ldots, E_8$ and $E_1', \ldots, E_8'$ 
are exceptional curves over two singular points of type $D_8$.
We will explain the symbols $D_2, D_2'$ and $\delta_1$, $\delta_2$ below.

In this section, we study the automorphism of Kummer surfaces with $p$-rank 1 under the following assumption.
\begin{assumption}\label{assum2}  
$\NS(J(C_1))$ has rank 1 (hence $\NS(X_1)$ has rank 17). 
We call such $X_1$ {\it Picard general}.
\end{assumption}

Proposition~\ref{Koizumi-Mori} in the previous section has the following analogue.
\begin{proposition}
We consider 
the locus ${\mathcal M}_{\leq 1}$ of curves of genus 2 and $p$-rank $\le 1$ in ${\mathcal M}_2$.
Then the generic curve  $C$  of genus 2 and  $p$-rank $\le 1$ over the algebraic closure of the function field of 
an irreducible component ${\mathcal N}$ of ${\mathcal M}_{\leq 1}$ is endo-general.
\end{proposition}
\begin{proof}
There exists a simple abelian surface  $A_0$  with $p$-rank 1 over a finite field by Honda~\cite[\S3.1]{Hon} (see also \cite{LO}).
Replacing by isogeny, there exists such an $A_0$ with a principal polarization (Mumford \cite[\S 23, Cor.1]{M1}).
Since $A_0$ is simple, $A_0 = J(C_0)$ for a curve  $C_0$ of genus 2.
Choose an irreducible component  ${\mathcal N}$  which contains the moduli point of this curve $C_0$.
Then by the same argument as \cite[Theorem 3.2]{Koi}, we have our assertion.
\end{proof}

By the assumption \ref{assum2}, $\NS(X_1)$ has signature $(1,16)$.
\begin{lemma}\label{PicardLatticeP-rank1} 
$\NS(X_1)$ is generated by the $18$ $(-2)$-curves in 
Figure \ref{p-rank1dualgraph} and is 
isomorphic to $U\oplus E_8\oplus D_7$.
\end{lemma}
\begin{proof}
Consider the linear system 
\[|E_3'+2E_2'+3E_1'+4T_2+5E_{7}+6E_6+3E_8+4E_5+2E_4|\]
which defines an elliptic fibration with singular fiber of type $\tilde{E}_8$ with a section $E_4'$.  Since
$\NS(X_1)$ has rank 17 and seven $(-2)$-curves $E_5', E_6', E_7', E_8', T_1, E_1, E_2$ forming a Dynkin diagram of type $D_7$ are contained in a singular fiber, this singular fiber is of type $\tilde{D}_7$, that is, there exists a $(-2)$-curve $D_2$ such that this fiber is given by
\[D_2 + E_2+ E_5'+E_8'+2(E_6'+E_7'+T_1+E_1)\]
(see Figure \ref{p-rank1dualgraph}).
Thus this fibration has two singular fibers of type $\tilde{E}_8$ and of type 
$\tilde{D}_7$ and a section $E_4'$. 
This implies that $\NS(X_1)$ contains $U\oplus E_{8}\oplus D_{7}$
as a sublattice of finite index.  Since there are no even unimodular lattices of signature $(1,16)$, we have
$\NS(X_1)\cong U \oplus E_{8}\oplus D_{7}$.  
\end{proof}

The $(-2)$-curve $D_2'$ in Figure \ref{p-rank1dualgraph} is obtained similarly. Since $D_2'$
meets the fiber of type $\tilde{E}_8$ in the first fibration with multiplicity $4$, we have $D_2\cdot D_2'=4$.

\begin{lemma}\label{numtrivialp-rank1}
The natural map $\Aut(X_1)\to O(\NS(X_1))$ is injective.
\end{lemma}
\begin{proof}
Let $g \in \Ker(\Aut(X_1)\to O(\NS(X_1)))$.
Consider an elliptic fibration $\pi$ defined by the linear system
\[|E_3+2E_4+3E_5+4E_6+2E_8+3E_7+2T_2+E_1'|\]
which has two singular fibers of type $\tilde{E}_7$ and four sections $E_2, E_2', D_2, D_2'$.
If $g$ acts trivially on the base of $\pi$, it induces an automorphism of a general fiber $F$ which fixes
four intersection points with the four sections.  
Thus we have $g=1$.  
If $g$ acts non-trivially on the base of $\pi$, 
then $g$ fixes three intersection points of $D_2$ 
with $E_1, E_3, D_2'$, which implies $g=1$.
\end{proof}

As in the previous section, we will first discuss a finite polyhedron 
$\calC(X_1)$ in 
$\NS(X_1)_\bfR$ by restricting the Conway's fundamental domain $\calC$.
Let $D_9$ be a root sublattice of $L$ generated by the following Leech roots
\[\alpha_0, \ \alpha_1, \ \alpha_2,\ \alpha_3,\ \alpha_4 =(1,1,[Q_0]),\ \alpha_5= (1,1,[\infty]),\]
\[ \alpha_6=(1,1,[\infty 0]),\ \alpha_7 = (1,1,[\infty 0.14.23]),\ \alpha_8 =(1,1,[14]),\]
where $\alpha_0, \alpha_1, \alpha_2, \alpha_3$ are as in (\ref{xyz}), and 
$[Q_0]$ has 2 at $\II$, $\III$ and at six points of the oval $Q_0$ and 0's elsewhere.
(see Figure \ref{D9} and \cite[Figure 23.20]{CS}).  We also denote by $D_8$ (resp. $D_7$) a root sublattice
generated by $\alpha_0, \alpha_1,\ldots, \alpha_7$ (resp. by $\alpha_0, \alpha_1,\ldots, \alpha_6$).

\begin{figure}[h]
\begin{center}
\scalebox{0.8}{\xy 
(-55,25)*{};
(80,-15)*{};
@={(0,0),(15,0),(30,0),(45,0),(60,0),(-15,0),(-30,0),(-45,0),(-30,15),}@@{*{\bullet}};
(0,0)*{};(15,00)*{}**\dir{-};(15,00)*{};(30,0)*{}**\dir{-};(30,0)*{};(45,0)*{}**\dir{-};
(0,0)*{},(-15,0)*{}**\dir{-};
(-15,0)*{};(-30,0)*{}**\dir{-};(-30,0)*{};(-45,0)*{}**\dir{-};
(-30,0)*{};(-30,15)*{}**\dir{-};(45,0)*{};(60,0)*{}**\dir{-};
(-45,-4)*{\alpha_{3}};(-30,-4)*{\alpha_0};(-15,-4)*{\alpha_{1}};(0,-4)*{\alpha_4};(15,-4)*{\alpha_5};(30,-4)*{\alpha_6};(45,-4)*{\alpha_7};(-36,15)*{\alpha_2};(60,-4)*{\alpha_8};
\endxy}
\end{center}
\caption{Root lattice of type $D_9$}\label{D9}
\end{figure}

The orthogonal complement $D_9^{\perp}$ of $D_9$ in $L$ is isomorphic to 
$U\oplus E_8\oplus D_7$ and the Leech roots perpendicular to $D_9$ are exactly 
18 vertices in Figure \ref{Dndualgraph} except $[14]$, $[23]$, $[\infty 2.03.14]$, 
$[\infty 3.02.14]$.  The Leech roots perpendicular to $D_8$ (resp. $D_7$) are exactly 
20 vertices (resp. all 22 vertices) in Figure \ref{Dndualgraph} except $[14]$, $[23]$.  
 The latter two cases will be used in the case of the Kummer surface 
 associated with the product of two non-isogenous ordinary elliptic curves in \S \ref{subs5-1} 
and the one with the self-product of an ordinary elliptic curve in \S \ref{subs5-2}.

\begin{figure}[h]
\begin{center}
\scalebox{0.8}{\xy 
(80,55)*{};
(-55,-55)*{};
@={(40,0),(0,-40),(-40,0),(0,40),(27,27),(27,-27),(-27,27),(-27,-27),(15,35),(15,-35),(-15,35),(-15,-35),(35,15),(35,-15),(-35,15),(-35,-15),(-25,0),(25,0)}@@{*{\bullet}};
(0,25)*{{\bullet}};(0,-25)*{{\bullet}};
(-5,-5)*{{\bullet}};(5,5)*{{\bullet}};
(40,0)*{};(35,15)*{}**\dir{-};(35,15)*{};(27,27)*{}**\dir{-};(27,27)*{};(15,35)*{}**\dir{-};
(15,35)*{},(0,40)*{}**\dir{-};
(0,40)*{};(-15,35)*{}**\dir{-};(-15,35)*{};(-27,27)*{}**\dir{-};
(-27,27)*{};(-35,15)*{}**\dir{-};(-35,15)*{};(-40,0)*{}**\dir{-};
(-40,0)*{};(-35,-15)*{}**\dir{-};(-35,-15)*{};(-27,-27)*{}**\dir{-};(-27,-27)*{};(-15,-35)*{}**\dir{-};
(-15,-35)*{},(0,-40)*{}**\dir{-};
(0,-40)*{};(15,-35)*{}**\dir{-};(15,-35)*{};(27,-27)*{}**\dir{-};
(27,-27)*{};(35,-15)*{}**\dir{-};(35,-15)*{};(40,0)*{}**\dir{-};
(40,0)*{};(25,0)*{}**\dir{-};(-40,0)*{};(-25,0)*{}**\dir{-};
(0,40)*{};(0,25)*{}**\dir{-};(0,-40)*{};(0,-25)*{}**\dir{-};
(5,6)*{};(0,25)*{}**\dir{-};(5,4)*{};(0,-25)*{}**\dir{-};
(-4,-5)*{};(25,0)*{}**\dir{-};(-6,-5)*{};(-25,0)*{}**\dir{-};
(8,8)*{[14]};(-8,-8)*{[23]};
(44,0)*{[01]};(-44,0)*{[04]};(0,-44)*{[02]};(0,44)*{[03]};
(30,30)*{[24]};(30,-30)*{[34]};(-30,30)*{[12]};(-30,-30)*{[13]};
(24,38)*{[\infty 1.03.24]};(24,-38)*{[\infty 1.02.34]};(-24,38)*{[\infty 4.03.12]};
(-24,-38)*{[\infty 4.02.13]};
(44,18)*{[\infty 3.01.24]};(44,-18)*{[\infty 2.01.34]};(-44,18)*{[\infty 3.04.12]};
(-44,-18)*{[\infty 2.04.13]};
(10,25)*{[\infty 2.03.14]};(10,-25)*{[\infty 3.02.14]};(25,4)*{[\infty 4.01.23]};(-25,4)*{[\infty 1.04.23]};
\endxy}
\end{center}
\caption{Leech roots: $D_9^{\perp}, D_8^{\perp}, D_7^{\perp}$}
\label{Dndualgraph}
\end{figure}

Now we identify 18 $(-2)$-curves except $D_2, D_2'$ in Fig.\ref{p-rank1dualgraph} with
18 vertices in Fig.\ref{Dndualgraph} except $[14]$, $[23]$, 
$[\infty 2.03.14]$, $[\infty 3.02.14]$.
The remaining faces of the finite polyhedron 
by restricting $\calC$ are defined by the projections of
Leech roots $r$ satisfying $\langle r, D_9\rangle$ being negative definite, that is, 
$\langle r, D_9\rangle$ is a root lattice of type $D_{10}$.
These Leech roots are given by
\[(1,1,[\infty 2.03.14]), \ (1,1,[\infty 3.02.14])\]
(see \cite[Figure 23.20]{CS}), whose projections 
$\delta_1, \delta_2$ in 
$D_9^{\perp}\otimes \bfQ$ have the intersection number $(-1)$.  In Figure \ref{p-rank1dualgraph}, we denote these
projections by two vertices (circle).  
The incidence relation between two $(-1)$-classes and 
18 $(-2)$-classes are given in Figure \ref{p-rank1dualgraph}.   

Now we show that each $(-1)$-class defines an involution of $X_1$ 
as in Lemma \ref{proj}.
Consider an elliptic fibration $\pi_1$ 
defined by a complete linear system
\[|E_{3}+2E_{4}+3E_5+4E_{6} + 3E_{7}+2T_{2}+E_1'+2E_{8}|\]
with two singular fibers of type $\tilde{E}_{7}$ and choose $E_{2}'$ as a zero section $(O)$.
As in the proof of Lemma \ref{proj},
$\pi_1$ has Mordell-Weil rank 1 and 
$\delta_1$ is perpendicular to the trivial lattice of $\pi_1$.
Let $\iota_1$ be the inversion of $\pi_1$.
Then it follows from Theorem \ref{M-W} and the same proof as that of Lemma \ref{proj} that 
which $\iota_1$ acts on $\NS(X)$ as the reflection associated with $\delta_1$.
In particular, $\iota_{b,\delta}^*(w) = w + a\delta_1 \ (a > 0)$, where $w$ is
the projection of the Weyl vector $\rho$.
Similarly we have an involution 
$\iota_2$ satisfying the same property.
Thus we have the following theorem
by the same proof as that of Theorem \ref{MainCurve}.

\begin{theorem}\label{MainCurvep-rank1}
Let $X_1$ be the Kummer surface associated with a curve of genus $2$ whose Jacobian has 
$p$-rank $1$.
Assume that $X_1$ is Picard general.
Let $G_1$ be the group generated by two involutions $\iota_1, \iota_2$ and let $(\bfZ/2\bfZ)^2$ 
be generated by two involutions
$\sigma, \varphi$ given by {\rm (\ref{Cremonaprank1}), (\ref{translprank1})}.  Then 
\[\Aut(X_1) \cong G_1 \rtimes (\bfZ/2\bfZ)^2.\]
\end{theorem}
\begin{proof}
Note that the symmetry group $\Aut(\calC(X_1))$ is isomorphic to 
$(\bfZ/2\bfZ)^2$
which are
induced from involutions $\sigma$, $\varphi$ given in (\ref{Cremonaprank1}), (\ref{translprank1}).
The assertion follows from Lemma \ref{numtrivialp-rank1}.
\end{proof}

Finally we consider the two projections of the Kummer quartic surface $S_1 \subset \bfP^3$  from a singular point.
It is easily verified that these give two double coverings  $X_1 \to \bfP^2$.
Moreover one covering involution interchanges $\{E_2, D_2\}$ and preserves other $(-2)$-curves in Figure~\ref{p-rank1dualgraph}.
Hence it coincides with  $\iota_1$.
The other interchanges  $\{E_2', D_2'\}$ and coincide with  $\iota_2$.
Hence  $\iota_1, \iota_2$  are cubic Cremona transformations
\begin{equation}\label{2 Cremonas/p=2}
pr_1 : [x, y, z, w]  \dashrightarrow \left[x, y, z, w + \frac{x^2z}{y^2}\right], \quad  
pr_2 : [x, y, z, w]  \dashrightarrow \left[x, y + \frac{xz^2}{w^2}, z, w\right].
\end{equation}

Hence we have the following corollary.

\begin{corollary}
$\Aut(X_1)$  is generated by  three cubic Cremona involutions 
\[\sigma :  [x,y,z,w]\dashrightarrow [xz^2, yz^2, b x^2z, b x^2w],\]
$pr_1, pr_2$  and the linear one $\varphi : [x, y, z, w] \mapsto [z, w, b x, b y]$.
\end{corollary}

\section{Kummer surfaces associated with the product of two ordinary elliptic curves}\label{sec5}

In this section, we 
consider Kummer surfaces associated with the product of two ordinary elliptic curves.  
Let $E, F$ be two ordinary elliptic curves and let $a, b$ be the unique non-zero 2-torsion points of $E, F$, respectively.  
Let $\iota_E, \iota_F$ be the inversions of $E, F$ respectively.  
Then $\iota = \iota_E\times \iota_F$ the inversion of the abelian surface $E\times F$.
Let $X' = \Km(E\times F)$ be the minimal resolution of 
the quotient surface $(E\times F)/\la \iota\ra$ which is also called a Kummer surface associated with  $E\times F$.  
If $E, F$  are defined by the Artin-Schreier coverings
\[
y^2 + y = x + \frac {c}x, \quad y^2 + y = x + \frac {d}x
\]
of  $\bfP^1_x$ for nonzero constants  $c, d$, then the quotient surface $(E\times F)/\la \iota\ra$  is defined by the Artin-Schreier covering
\begin{equation}\label{4D4-eq in char 2}
z^2 + z = x + \frac {c}x + y + \frac {d}y,
\end{equation}
or equivalently,
\[
w^2 + xyw = xy(x^2y + cy + xy^2 + dx),
\]
of  $\bfP^1_x \times \bfP^1_y$.
It has four rational double points of type $D_4$ at  $(0, 0), (0, \infty), (\infty, 0), (\infty, \infty)$
(Shioda \cite{Shioda}).
 Its minimal resolution $X'$ has 20 $(-2)$-curves corresponding to 
20 black vertices in Figure \ref{ExFdualgraph}.  

\begin{figure}[h]
\begin{center}
\scalebox{0.8}{\xy 
(80,55)*{};
(-55,-55)*{};
@={(40,0),(0,-40),(-40,0),(0,40),(27,27),(27,-27),(-27,27),(-27,-27),(15,35),(15,-35),(-15,35),(-15,-35),(35,15),(35,-15),(-35,15),(-35,-15),(0,25),(0,-25),(-25,0),(25,0)}@@{*{\bullet}};
(-5,-5)*{{\circ}};(5,5)*{{\circ}};
(40,0)*{};(35,15)*{}**\dir{-};(35,15)*{};(27,27)*{}**\dir{-};(27,27)*{};(15,35)*{}**\dir{-};
(15,35)*{},(0,40)*{}**\dir{-};
(0,40)*{};(-15,35)*{}**\dir{-};(-15,35)*{};(-27,27)*{}**\dir{-};
(-27,27)*{};(-35,15)*{}**\dir{-};(-35,15)*{};(-40,0)*{}**\dir{-};
(-40,0)*{};(-35,-15)*{}**\dir{-};(-35,-15)*{};(-27,-27)*{}**\dir{-};(-27,-27)*{};(-15,-35)*{}**\dir{-};
(-15,-35)*{},(0,-40)*{}**\dir{-};
(0,-40)*{};(15,-35)*{}**\dir{-};(15,-35)*{};(27,-27)*{}**\dir{-};
(27,-27)*{};(35,-15)*{}**\dir{-};(35,-15)*{};(40,0)*{}**\dir{-};
(40,0)*{};(25,0)*{}**\dir{-};(-40,0)*{};(-25,0)*{}**\dir{-};
(0,40)*{};(0,25)*{}**\dir{-};(0,-40)*{};(0,-25)*{}**\dir{-};
(4,4)*{};(-4,-4)*{}**\dir{--};
(5,6)*{};(0,25)*{}**\dir2{-};(5,4)*{};(0,-25)*{}**\dir2{-};
(-4,-5)*{};(25,0)*{}**\dir2{-};(-6,-5)*{};(-25,0)*{}**\dir2{-};
(8,8)*{\delta_2};(-8,-8)*{\delta_{1}};
(44,0)*{E_2^0};(-44,0)*{E_4^0};(0,-44)*{E_3^0};(0,44)*{E_1^0};
(30,30)*{\overline{F}'};(30,-30)*{\overline{E}'};(-30,30)*{\overline{E}};(-30,-30)*{\overline{F}};
(18,38)*{E_1^3};(18,-38)*{E_3^2};(-18,38)*{E_1^2};(-18,-38)*{E_3^3};
(38,18)*{E_2^2};(38,-18)*{E_2^3};(-38,18)*{E_4^3};(-38,-18)*{E_4^2};
(4,25)*{E_1^1};(4,-25)*{E_3^1};(25,4)*{E_2^1};(-25,4)*{E_4^1};
(-2,2)*{4};
\endxy}
\end{center}
\caption{The dual graph:  
$\Km(E\times F)$ (ordinary case)}
\label{ExFdualgraph}
\end{figure}

\noindent
Here we denote by $\overline{E}$, $\overline{E}'$, $\overline{F}$, $\overline{F}'$ the images of $E\times \{0\}, E\times \{b\}, \{0\}\times F, \{a\}\times F$ on $X'$, respectively. 
The remaining 16 $(-2)$-curves are the exceptional curves over the
four rational double points of type $D_4$ which are the images of four 2-torsion points on $E\times F$.

\subsection{The case $E$ is not isogenous to $F$.}\label{subs5-1}

First we assume that $E$ is not isogenous to $F$, that is, the Picard lattice of $X'$ has rank 18.

\begin{lemma}\label{PicardLatticeExF} 
$\NS(X')$ is generated by $20$ $(-2)$-curves in Figure \ref{ExFdualgraph}
and is isomorphic to $U\oplus E_8\oplus D_8$.
\end{lemma}
\begin{proof}
The elliptic fibration defined by a linear system
\[|2\overline{E} + 4E_4^3 + 6E_4^0+3E_4^1+ 5E_4^2+4\overline{F}+3E_3^3+2E_3^0+E_3^2|\]
has a singular fiber of type $\tilde{E}_8$ with a section $\overline{E}'$.
Since the Picard number is 18 and 8 $(-2)$-curves
\[E_1^1, E_1^0, E_1^3, \overline{F}', E_2^2, E_2^0, E_2^3, E_2^1\]
form a dual graph of type $D_8$, the fibration has a singular fiber of type $\tilde{D}_8$.
This implies that $\NS(X')$ contains $U\oplus E_8\oplus D_8$ as a sublattice.  
If $\NS(X') \ne U\oplus E_8\oplus D_8$, then $\NS(X')\cong U\oplus E_8\oplus E_8$.  In this case,
it is known (Vinberg \cite{V}) that the dual graph of $(-2)$-curves on $X'$ is given in Fig. \ref{ExFssdualgraph} (see Lemma \ref{PLatticeExF}), which is impossible.
\end{proof}

\begin{remark}
The lattice is the even part of the odd unimodular lattice $I_{1,17}$ of signature $(1,17)$ (\S~\ref{odd unimodular}).  
Thus Figure \ref{ExFdualgraph} is also a fundamental domain of
the reflection group of $I_{1,17}$ in Vinberg \cite[Table 4, A, p.347]{V}.
\end{remark}

\begin{lemma}\label{numtrivialproduct1}
The kernel of the natural map $\Aut(X')\to O(\NS(X'))$ is generated by the covering involution $\sigma : z \mapsto z+1$ of \eqref{4D4-eq in char 2}.
\end{lemma}
\begin{proof}
First note that 
the involution $\sigma$ acts trivially on the dual graph  
and thus is contained in the kernel since twenty $(-2)$-curves generate $\NS(X')$.
Conversely, let $g \in \Ker(\Aut(X')\to O(\NS(X')))$.   
Then $g$ acts trivially
on the above dual graph, and hence preserves each of 20 $(-2)$-curves in Figure \ref{ExFdualgraph}.
Consider the elliptic fibrations $\pi_1: X' \to \bfP^1_x, \pi_2: X' \to \bfP^1_y$ induced from the projections 
$E\times F \to E$, $E\times F \to F$.    
Both fibrations have a section and two singular fibers of type $\tilde{D}_8$ all of which are preserved by $g$.
Hence  $g$  preserves the double covering  $(\pi_1, \pi_2) : X' \to \bfP^1_x \times \bfP^1_y$
defined by  \eqref{4D4-eq in char 2}.
Therefore, $g$  is a lift of the automorphism  $(x, y) \mapsto (kx, ly)$  of  $\bfP^1_x \times \bfP^1_y$  for constants $k, l \ne 0$.
Assume that 
\[
 \left(kx + \frac {c}{kx} + ly + \frac {d}{ly}\right) -  \left(x + \frac {c}x + y + \frac {d}y\right) = f(x, y)^2 + f(x, y)
\]
holds for a rational function $f(x,y)$  on  $\bfP^1_x \times \bfP^1_y$.
Then it is easy to deduce that  $k=l=1$  and  $f(x, y) \in \bfF_2$  since  $f(x,y)$ cannot have a pole along any divisor of   $\bfA^1_x \times \bfA^1_y$.
Hence $g$  is either the identity or the involution $\sigma$ according as $f(x, y) \equiv 0$ or 1.
\end{proof}

\begin{remark}\label{(1,-1)}
The covering involution $\sigma$ in the lemma is induced from the involution $(1_E, \iota_F)$ of $E\times F$. 
\end{remark}

As in the previous section, 
we will first discuss a finite polyhedron $\calC(X')$ in $\NS(X')_\bfR$ 
by restricting the Conway's fundamental domain $\calC$.
Let $R = D_8$ be a root sublattice of $L$ generated by the Leech roots
$\alpha_0,\ \alpha_1,\ \alpha_2,\ \alpha_3,\alpha_4,\ \alpha_5,\ \alpha_6,\ \alpha_7$
as in Figure \ref{D9}.
Note that the orthogonal complement $R^{\perp}$ of $R$ in $L$ is isomorphic to the Picard lattice 
$\NS(X')$.  
Among these 42 Leech roots perpendicular to $\alpha_0, \alpha_1, \alpha_2, \alpha_3$, 
there are exactly 20 Leech roots perpendicular to $D_8$ which form
the graph given in Figure \ref{ExFdualgraph} (see also Figure \ref{Dndualgraph}). 
Under the identification between $R^{\perp}$ and $\NS(X')$, we can consider
these 20 Leech roots as $(-2)$-curves on $X'$ in Figure \ref{ExFdualgraph}.
Thus these 20 $(-2)$-curves define 20 faces of the restriction $\calC(X')$ of 
the Conway's fundamental domain $\calC$.  
The remaining faces are the restrictions of Leech roots $r$ such that $r$ and $R$ generate a negative definite lattice, that is, a root lattice.  There are two possibilities, that is, 
$\la r, R\ra = D_8 \oplus A_1$ or $D_9$.  In the first case, $r$ is nothing but one of the above 20 $(-2)$-curves.

\begin{lemma}\label{boundaries2}
There are exactly two Leech roots $r$ satisfying $\la r, R\ra = D_{9}$.
Their projections into $\NS(X')\otimes \bfQ$ have the norm $-1$.
\end{lemma}
\begin{proof}
The desired ones are $(1,1,[14]), (1,1,[23])$ (see \cite[Figure 23.20]{CS}).
\end{proof}

Let $\delta_1, \delta_2$ be the projection of $2(1,1,[14])$, $2(1,1,[23])$ in $\NS(X')$.  Then 
$\delta_1^2=\delta_2^2 = -4$, $\delta_1\cdot \delta_2= 4$, and the incidence relation of $\delta_1, \delta_2$
with twenty $(-2)$-curves are given in Figure \ref{ExFdualgraph}.

\begin{lemma}\label{translation2}
For each vertex $\delta_i$ $(i=1,2)$, there exists an involution $\iota_i$ of $X'$ such that $\iota_i^*(\delta_i)=-\delta_i$ and 
$\iota_i^*(w)= w+ a\delta_i$ \ $(a>0)$, where $w$ is
the projection of the Weyl vector $\rho$. 
\end{lemma}
\begin{proof}
Consider an elliptic fibration $\pi_1$
defined by a complete linear system of divisor of sixteen $(-2)$-curves forming the circle in Figure
\ref{ExFdualgraph} with a singular fiber of type $\tilde{A}_{15}$ and choose 
$E_{1}^1$ as a zero section $(O)$.   
As in the proof of Lemma \ref{proj}, $\pi_1$ has Mordell-Weil rank 1 and
$\delta_1$ is perpendicular to the trivial lattice of $\pi_1$.  
Let $\iota_1$ be the inversion of $\pi_1$.  In this case, $\iota_1$ preserves
the trivial lattice $\Triv(\pi_1)$, and hence $\iota_1^*(w)= w+ a\delta_i$ \ $(a>0)$ holds.  
Similarly we have an involution 
$\iota_2$ satisfying the same property. 
Thus we have the assertion.
\end{proof}

\begin{theorem}\label{MainProd1}
Let $X'$ be the Kummer surface associated with the product of two ordinary elliptic curves $E, F$.
We assume that $E$ is not isogenous to $F$.
Let $G'$ be the group generated by two involutions $\iota_1, \iota_2$ and let $(\bfZ/2\bfZ)^3$ be generated by two involutions induced from the translations by $2$-torsion points of $E\times F$ and 
the involution induced from $\iota_E \times 1_F$.
Then 
\[\Aut(X') \cong G' \rtimes (\bfZ/2\bfZ)^3.\]
\end{theorem}
\begin{proof}
Note that the symmetry group $\Aut(\calC(X'))$ is isomorphic to the dihedral group $D_8$, but only the subgroup $(\bfZ/2\bfZ)^2$ can be  realized as automorphisms because $E$ is not isogenous to $F$.
The assertion now follows from Lemma \ref{numtrivialproduct1}.
\end{proof}

\subsection{The case $E=F$.}\label{subs5-2}
Next we assume that $E = F$ and $X''$ is Picard general, i.e., the Picard number of $X''$ is 19.
An ordinary elliptic curve $E$  and  the Kummer surface $X''=\Km(E\times E)$  are the Artin-Schreier coverings
$y^2 + y = x + c/x$  and 
\begin{equation}\label{type I in char 2}
z^2 + z = x + \frac {c}x + y + \frac {c}y,
\end{equation}
the special case $c=d$  in \eqref{4D4-eq in char 2}, respectively.

Consider the following curves on $E\times E$:
\[\Delta =\{ (x,x)\ : \ x \in E\},\ \Delta' = \{ (x,\iota_E(x))\ : \ x \in E\},\]
\[\Delta''= \{ (x,x+a)\ : \ x \in E\}, \Delta'''= \{ (x,\iota_E(x)+a)\ : \ x \in E\},\]
where $a$ is the unique non-zero 2-torsion of $E$.
All of them are invariant under the inversion $\iota$ and hence their images on 
$X''=\Km(E\times E)$ are
non-singular rational curves $D, D', D'', D'''$ satisfying  
$D\cdot D'= D\cdot D'' = D'\cdot D'''= D''\cdot D'''=0$,
$D\cdot D'''=D'\cdot D''=2$.  See the following Figure \ref{ExEdualgraph}.

\begin{figure}[h]
\begin{center}
\scalebox{0.8}{\xy 
(80,55)*{};
(-55,-55)*{};
@={(40,0),(0,-40),(-40,0),(0,40),(27,27),(27,-27),(-27,27),(-27,-27),(15,35),(15,-35),(-15,35),
(-15,-35),(35,15),(35,-15),(-35,15),(-35,-15),(0,25),(0,-25),(-25,0),(25,0),(-10,0),(10,0),(0,10),(0,-10)}@@{*{\bullet}};
(40,0)*{};(35,15)*{}**\dir{-};(35,15)*{};(27,27)*{}**\dir{-};(27,27)*{};(15,35)*{}**\dir{-};
(15,35)*{},(0,40)*{}**\dir{-};
(0,40)*{};(-15,35)*{}**\dir{-};(-15,35)*{};(-27,27)*{}**\dir{-};
(-27,27)*{};(-35,15)*{}**\dir{-};(-35,15)*{};(-40,0)*{}**\dir{-};
(-40,0)*{};(-35,-15)*{}**\dir{-};(-35,-15)*{};(-27,-27)*{}**\dir{-};(-27,-27)*{};(-15,-35)*{}**\dir{-};
(-15,-35)*{},(0,-40)*{}**\dir{-};
(0,-40)*{};(15,-35)*{}**\dir{-};(15,-35)*{};(27,-27)*{}**\dir{-};
(27,-27)*{};(35,-15)*{}**\dir{-};(35,-15)*{};(40,0)*{}**\dir{-};
(40,0)*{};(25,0)*{}**\dir{-};(-40,0)*{};(-25,0)*{}**\dir{-};
(0,40)*{};(0,25)*{}**\dir{-};(0,-40)*{};(0,-25)*{}**\dir{-};
(25,0)*{};(0,10)*{}**\dir{-};(25,0)*{};(0,-10)*{}**\dir{-};
(-25,0)*{};(0,10)*{}**\dir{-};(-25,0)*{};(0,-10)*{}**\dir{-};
(0,25)*{};(10,0)*{}**\dir{-};(0,25)*{};(-10,0)*{}**\dir{-};
(0,-25)*{};(10,0)*{}**\dir{-};(0,-25)*{};(-10,0)*{}**\dir{-};
(-10,0)*{};(0,10)*{}**\dir2{-};(10,0)*{};(0,-10)*{}**\dir2{-};
(44,0)*{E_2^0};(-44,0)*{E_4^0};(0,-44)*{E_3^0};(0,44)*{E_1^0};
(30,30)*{\overline{F}'};(30,-30)*{\overline{E}'};(-30,30)*{\overline{E}};(-30,-30)*{\overline{F}};
(18,38)*{E_1^3};(18,-38)*{E_3^2};(-18,38)*{E_1^2};(-18,-38)*{E_3^3};
(38,18)*{E_2^2};(38,-18)*{E_2^3};(-38,18)*{E_4^3};(-38,-18)*{E_4^2};
(4,25)*{E_1^1};(4,-25)*{E_3^1};(25,4)*{E_2^1};(-25,4)*{E_4^1};
(0,14)*{D};(0,-14)*{D'};(-14,0)*{D'''};(14,0)*{D''};
\endxy}
\end{center}
\caption{The dual graph: 
$\Km(E\times E)$ (ordinary case)}
\label{ExEdualgraph}
\end{figure}

\begin{lemma}\label{PicardLatticeExE} 
$\NS(X'')$ is generated by $24$ $(-2)$-curves in Figure {\rm \ref{ExEdualgraph}}
and is isomorphic to $U\oplus E_8\oplus D_9$.
\end{lemma}
\begin{proof}
The elliptic fibration defined by a linear system
\[|E_2^3+2E_2^0+3E_2^2+4\overline{F}'+5E_1^3+6E_1^0+3E_1^1+4E_1^2+2\overline{E}|\]
\[=|E_3^1+E_3^2+D+D'+2(E_3^0+E_3^3+\overline{F}+E_4^2+E_4^0+E_4^1)|\]
has singular fibers of type $\tilde{E}_8$, $\tilde{D}_9$ and has a section $\overline{E}'$.
This implies that $\NS(X'')$ contains $U\oplus E_8\oplus D_9$ as a sublattice of finite index.  
Since there is no even unimodular lattice of signature $(1,18)$, we have
$\NS(X'')= U\oplus E_8\oplus D_9$.
\end{proof}

The covering involution 
\begin{equation}\label{key involution}
\sigma: (x,y,z) \mapsto (x, y, z+1)
\end{equation}
 of \eqref{type I in char 2}, similarly defined as Lemma~\ref{numtrivialproduct1},  is no more Picard trivial.
More generally, we have

\begin{lemma}\label{numericaltrivialExE}
The natural map $\Aut(X'')\to O(\NS(X''))$ is injective.
\end{lemma}
\begin{proof}
The proof is similar to Lemma \ref{numtrivial} by using an elliptic fibration with three singular fibers
of type $\tilde{A}_{15}, \tilde{A}_1, \tilde{A}_1$ and four sections $E^1_{i}$ $(i=1,2,3,4)$.
\end{proof}

Let $R$ be a root sublattice $D_7$ generated by Leech roots
$\alpha_0, \alpha_1, \alpha_2, \alpha_3,\alpha_4, \alpha_5, \alpha_6$
given  in Figure \ref{D9}.
Note that the orthogonal complement $R^{\perp}$ of $R$ in $L$ is isomorphic to the Picard lattice 
$\NS(X'')$.  
We can easily check that there are exactly $22$ Leech roots perpendicular to $R$ which 
form the graph given in Figure \ref{Dndualgraph}.  This graph coincides with 
Figure \ref{ExEdualgraph} except for $D'$, $D'''$.
Under the identification between $R^{\perp}$ and $\NS(X'')$, we can consider
these 22 Leech roots as $(-2)$-curves on $X''$ in Figure \ref{ExEdualgraph} except for $D'$, $D'''$.  
Thus these 22 $(-2)$-curves define 22 faces of the restriction $\calC(X'')$ of 
the Conway's fundamental domain $\calC$.  

The symmetry group of Figure \ref{Dndualgraph} is a symmetric group 
$\mathfrak{S}_4$ of degree 4.
The symmetry group $\mathfrak{S}_3$ acting on $E\times \{0\}, \{0\}\times E, E\times E$ and
the involution $(1_E, t)$ induce automorphisms of $X''$, where $t$ is the translation of $E$ by the non-zero 2-torsion point.  The symmetry group 
$\mathfrak{S}_4$ can be realized by these automorphisms of $X''$. 

The remaining faces of $\calC(X'')$ are the restrictions of Leech roots $r$ such that $r$ and $R$ generate a negative definite lattice, that is, a root lattice.  There are two possibilities, that is, 
$\la r, R\ra = D_8$ or $E_8$.

\begin{lemma}\label{boundaries3}
There are exactly three Leech roots $r$ satisfying $\la r, R\ra = D_8$.
Their projections into $\NS(X'')\otimes \bfQ$ have the norm $-1$.
There are exactly twelve Leech roots with $\la r, R\ra = E_8$ whose projections have the norm $-1/4$.
\end{lemma}
\begin{proof}
The Leech roots $r= (1,1,\infty 0.14.23)$, $(1,1,\infty 0.12.34)$, $(1,1,\infty 0.13.24)$ are the ones with
$\la r, R\ra = D_8$.  
By \S\ref{Leech}, for $\la r, R\ra = E_8$, we have  $r = (1,1,[Q])$  for an oval  $Q$  which satisfies the following:
\[
|Q \cap Q_0| = 3, \quad Q \not\ni  \infty \quad {\rm and}  \quad Q \ni 0, \ |Q\cap \infty 0|=2,
\]
which are equivalent to  $\la r, \alpha_4\ra = \la r, \alpha_5\ra =\la r, \alpha_6\ra =0$.
Hence  $Q$  contains (exactly) two from $\{1, 2, 3, 4\}$.
By Lemma \ref{Ovals} below, we have $12 = \binom{4}{2}\times 2$ such ovals  $Q$,  and hence $12$ $(-1/4)$-Leech roots.
\end{proof}

\begin{lemma}\label{Ovals}
For a $3$-subset $T$ of the oval  $Q_0$,  there exist exactly two ovals $Q$  with  $Q \cap Q_0 = T$.
The corresponding octad is obtained by appending two romans $\I, \III$  for the one,
and $\I, \II$  for the other.
\end{lemma}
\begin{proof}
We may assume  $T= \{0, 1, 2\}$, and hence  $Q_0 \setminus T = \{\infty, 3, 4\}$.
Consider the pencil  $P$  of plane cubics generated by two triangles
$M_0 := \triangle 012$ and  $M_1 = \triangle \infty34$ corresponding to $T$ and  $Q_0 \setminus T$,
where a \emph{triangle} means a union of three lines.
Then  $P$  is a Hesse pencil and it contains two other (Maclaurin) triangles  $M_2 = \triangle ***$  and  
$M_3 = \triangle \bullet\bullet\bullet$,
from which the required ovals are  obtained by appending  $T$.
More explicitly, the octads corresponding to our two ovals $Q$  are
\begin{center}
 \begin{tabular}{|c|c|c|}
   \hline  
$*$ \ \ $\cdot$  & $\cdot$ \ \ $\cdot$ & $\cdot$ \ \ $\cdot$ \\ 
I \ \ $\cdot$  & $\cdot$ \ \ $\cdot$  &  $\cdot$ \ \ $\cdot$ \\    \hline 
\ $\cdot$\ \ $\cdot$  & $*$ \ \ $\cdot$ & 1 \ \ $\cdot$ \\ 
III  \ 0 &  $\cdot$ \ \ $*$ & $\cdot$ \ \ 2 \\
\hline
\end{tabular}
\quad {\rm and } \quad
  \begin{tabular}{|c|c|c|}
   \hline  
$\cdot$ \ \ $\bullet$  & $\cdot$ \ \ $\cdot$ & $\cdot$ \ \ $\bullet$ \\ 
I \ \ $\cdot$  & $\cdot$ \ \ $\cdot$  &  $\bullet$ \ \ $\cdot$ \\    \hline 
II \ \ $\cdot$  & $\cdot$ \ \ $\cdot$ & 1 \ \ $\cdot$ \\ 
$\cdot$ \ \ 0  & $\cdot$ \ \ $\cdot$ & $\cdot$ \ \ 2 \\
\hline
\end{tabular}
\end{center}
by Curtis's table of the standard sextets \cite[Figure 11.17]{CS}.
This also follows from the fact that the base points of the pencil are nine points `$\bullet$' in the following figure:
\begin{center}
 \begin{tabular}{|c|c|c|}
   \hline  
$\cdot$ \ \ $\cdot$  & $\bullet$ \ \ $\bullet$ & $\bullet$ \ \ $\cdot$ \\ 
I \ \ $\bullet$  & $\bullet$ \ \ $\bullet$  &  $\cdot$ \ \ $\bullet$ \\    \hline 
II\ \ $\infty$  & $\cdot$ \ \ $\bullet$ & 1 \ \ 3 \\ 
III  \ 0 &  $\bullet$ \ \ $\cdot$ & 4 \ \ 2 \\
\hline
\end{tabular}
\end{center}
\end{proof}

Denote the projection of $(1,1,\infty 0.14.23)$, $(1,1,\infty 0.12.34)$, 
$(1,1,\infty 0.13.24)$ by 
$\delta_{D,1}$,
$\delta_{D,2}$, $\delta_{D,3}$, respectively.
The incident relation of $\delta_{D,1}$ with 24 $(-2)$-curves is as in the following dual graph.

\begin{figure}[h]
\begin{center}
\scalebox{0.8}{\xy 
(80,55)*{};
(-55,-55)*{};
@={(40,0),(0,-40),(-40,0),(0,40),(27,27),(27,-27),(-27,27),(-27,-27),(15,35),(15,-35),(-15,35),
(-15,-35),(35,15),(35,-15),(-35,15),(-35,-15),(0,25),(0,-25),(-25,0),(25,0),(-10,0),(10,0),(0,10),(0,-10)}@@{*{\bullet}};
(15,15)*{{\circ}};
(40,0)*{};(35,15)*{}**\dir{-};(35,15)*{};(27,27)*{}**\dir{-};(27,27)*{};(15,35)*{}**\dir{-};
(15,35)*{},(0,40)*{}**\dir{-};
(0,40)*{};(-15,35)*{}**\dir{-};(-15,35)*{};(-27,27)*{}**\dir{-};
(-27,27)*{};(-35,15)*{}**\dir{-};(-35,15)*{};(-40,0)*{}**\dir{-};
(-40,0)*{};(-35,-15)*{}**\dir{-};(-35,-15)*{};(-27,-27)*{}**\dir{-};(-27,-27)*{};(-15,-35)*{}**\dir{-};
(-15,-35)*{},(0,-40)*{}**\dir{-};
(0,-40)*{};(15,-35)*{}**\dir{-};(15,-35)*{};(27,-27)*{}**\dir{-};
(27,-27)*{};(35,-15)*{}**\dir{-};(35,-15)*{};(40,0)*{}**\dir{-};
(40,0)*{};(25,0)*{}**\dir{-};(-40,0)*{};(-25,0)*{}**\dir{-};
(0,40)*{};(0,25)*{}**\dir{-};(0,-40)*{};(0,-25)*{}**\dir{-};
(25,0)*{};(0,10)*{}**\dir{-};(25,0)*{};(0,-10)*{}**\dir{-};
(-25,0)*{};(0,10)*{}**\dir{-};(-25,0)*{};(0,-10)*{}**\dir{-};
(0,25)*{};(10,0)*{}**\dir{-};(0,25)*{};(-10,0)*{}**\dir{-};
(0,-25)*{};(10,0)*{}**\dir{-};(0,-25)*{};(-10,0)*{}**\dir{-};
(-10,0)*{};(0,10)*{}**\dir2{-};(10,0)*{};(0,-10)*{}**\dir2{-};
(14,15)*{};(0,10)*{}**\dir{-};
(15,14)*{};(10,0)*{}**\dir{-};
(14,15)*{};(-10,0)*{}**\dir{--};
(15,14)*{};(0,-10)*{}**\dir{--};
(44,0)*{E_2^0};(-44,0)*{E_4^0};(0,-44)*{E_3^0};(0,44)*{E_1^0};
(30,30)*{\overline{F}'};(30,-30)*{\overline{E}'};(-30,30)*{\overline{E}};(-30,-30)*{\overline{F}};
(18,38)*{E_1^3};(18,-38)*{E_3^2};(-18,38)*{E_1^2};(-18,-38)*{E_3^3};
(38,18)*{E_2^2};(38,-18)*{E_2^3};(-38,18)*{E_4^3};(-38,-18)*{E_4^2};
(4,25)*{E_1^1};(4,-26)*{E_3^1};(26,4)*{E_2^1};(-26,4)*{E_4^1};
(0,14)*{D};(0,-14)*{D'};(-14,0)*{D'''};(14,0)*{D''};
(18,18)*{\delta_{D,1}};
(0,4)*{-1};(1,-3)*{-1};
\endxy}
\end{center}
\caption{The incidence relation of $\delta_{D,1}$}
\label{ExEdualgraph2}
\end{figure}

It is immediately seen that the involution  $\sigma$  in \eqref{key involution}  acts on the Picard group by reflection with respect to  $\delta_{D,1}$.
In fact, $\sigma$  switches $D$ and $D'''$, $D'$ and $D''$,  and leaves the other 20 $\bfP^1$'s invariant  (see Remark~\ref{(1,-1)}).
Moving $\sigma$  by the action of  $\mathfrak{S}_4 \subset \Aut(X'')$, we have three involutions $\sigma=\sigma_1, \sigma_2, \sigma_3$, and the following:

\begin{lemma}\label{translation3}
The involution $\sigma_i$  acts on the Picard group by reflection with respect to  $\delta_{D,i}$ for every $i=1, 2, 3$.
\end{lemma}

\begin{lemma}\label{translation4}
Let $\delta_{E,i}$ $(i=1,\ldots, 12)$ be the projections of twelve Leech roots satisfying $\la r, R\ra = E_8$.
For each vertex $\delta_{E,i}$ $(i=1,\ldots, 12)$, there exists an involution $\iota_{E,i}$  which acts on the Picard group by reflection with respect to  $\delta_{E,i}$. 
\end{lemma}
\begin{proof}  
We assume that $\delta_E= \delta_{E,i}$ is the projection of the Leech root $r=(1,1,Q)$ where
$Q$ is the oval given in the left hand side MOG in the proof of Lemma \ref{Ovals}, that is,
\[
Q= \{0,\ 1,\ 2,\ \infty 0.14.23,\ \infty 2.04.13,\ \infty 1.03.24\}.
\]
It follows that $\delta_E$ meets exactly three vertices $[34]$, $[\infty 1.03.24]$, $[\infty 2.04.13]$
in Figure \ref{Dndualgraph}. Thus the incidence relation of $\delta_E$ with 22 $(-2)$-curves is as in Figure
\ref{ExEdualgraph3}:

\begin{figure}[h]
\begin{center}
\scalebox{0.7}{\xy 
(-55,65)*{};
(80,-65)*{};
@={(40,0),(0,-40),(-40,0),(0,40),(27,27),(27,-27),(-27,27),(-27,-27),(15,35),(15,-35),(-15,35),
(-15,-35),(35,15),(35,-15),(-35,15),(-35,-15),(0,25),(0,-25),(-25,0),(25,0),
(-10,0),
(0,10),
}@@{*{\bullet}};
(10,-10)*{{\circ}};
(40,0)*{};(35,15)*{}**\dir{-};(35,15)*{};(27,27)*{}**\dir{-};(27,27)*{};(15,35)*{}**\dir{-};
(15,35)*{},(0,40)*{}**\dir{-};
(0,40)*{};(-15,35)*{}**\dir{-};(-15,35)*{};(-27,27)*{}**\dir{-};
(-27,27)*{};(-35,15)*{}**\dir{-};(-35,15)*{};(-40,0)*{}**\dir{-};
(-40,0)*{};(-35,-15)*{}**\dir{-};(-35,-15)*{};(-27,-27)*{}**\dir{-};(-27,-27)*{};(-15,-35)*{}**\dir{-};
(-15,-35)*{},(0,-40)*{}**\dir{-};
(0,-40)*{};(15,-35)*{}**\dir{-};(15,-35)*{};(27,-27)*{}**\dir{-};
(27,-27)*{};(35,-15)*{}**\dir{-};(35,-15)*{};(40,0)*{}**\dir{-};
(40,0)*{};(25,0)*{}**\dir{-};(-40,0)*{};(-25,0)*{}**\dir{-};
(0,40)*{};(0,25)*{}**\dir{-};(0,-40)*{};(0,-25)*{}**\dir{-};
(25,0)*{};(0,10)*{}**\dir{-};
(-25,0)*{};(0,10)*{}**\dir{-};
(0,25)*{};(-10,0)*{}**\dir{-};
(0,-25)*{};(-10,0)*{}**\dir{-};
(27,-27)*{};(10.5,-10.5)*{}**\dir{-};
(15,35)*{};(10.5,-9.5)*{}**\dir{-};
(-35,-15)*{};(9.5,-10.5)*{}**\dir{-};
(44,0)*{E_2^0};(-44,0)*{E_4^0};(0,-44)*{E_3^0};(0,44)*{E_1^0};
(30,30)*{\overline{F}'};(30,-30)*{\overline{E}'};(-30,30)*{\overline{E}};(-30,-30)*{\overline{F}};
(18,38)*{E_1^3};(18,-38)*{E_3^2};(-18,38)*{E_1^2};(-18,-38)*{E_3^3};
(38,18)*{E_2^2};(38,-18)*{E_2^3};(-38,18)*{E_4^3};(-38,-18)*{E_4^2};
(4,25)*{E_1^1};(4,-26)*{E_3^1};(26,4)*{E_2^1};(-26,4)*{E_4^1};
(0,14)*{D};
(-14,0)*{D''};
(15,-10)*{\delta_{E}};
\endxy}
\end{center}
\caption{The incidence relation of $\delta_{E}$}
\label{ExEdualgraph3}
\end{figure}
Now consider an elliptic fibration $\pi$
defined by a complete linear system
\[|E_{4}^3+2E_{4}^0+3E_4^1+4D+ 5E_{2}^1+6E_2^0+3E_2^3+4E_2^2+2\overline{F}'|\]
with two singular fibers of type $\tilde{E}_{8}$ and a section 
$\overline{E}$.
As in the proof of Lemma \ref{proj},
$\pi$ has Mordell-Weil rank 1 and 
$\delta_E$ is perpendicular to the trivial lattice of $\pi$.
Let $\iota_{E}$ be the inversion of $\pi$.
Then it follows from Theorem \ref{M-W}
that $\iota_{E}^*(\delta_{E})=-\delta_{E}$.
It is easy to see that  $\iota_{E}^*$  is a reflection since it leave 19 $\bfP^1$'s, which generate a sublattice of rank 17, invariant.
\end{proof}

Let $G$ be a subgroup of $\Aut(X'')$ generated by $\sigma_i, \iota_{E,j}$ $(i=1,2,3; j=1,\ldots, 12)$.  
By combining Lemma \ref{numericaltrivialExE}, we obtain the following theorem.

\begin{theorem}\label{MainProd2}
Let $X''$ be the Kummer surface associated with 
an ordinary elliptic curve $E$ in characteristic $2$.  Assume that
$X''$ is Picard general, i.e., the Picard number of $X''$ is $19$.
Then
$\Aut(X'')$ is isomorphic to $G\rtimes \mathfrak{S}_4$.
\end{theorem}

\subsection{Remark on Enriques surfaces covered by $\Km(E\times F)$}\label{subs5-3}
We recall that an Enriques surface of {\it Lieberman type} over  $\bfC$ is the quotient of a Kummer surface  $\Km(E \times F)$  by the composite $\varepsilon$ of two induced involutions: one is induced from $1_E \times \iota_F \in \Aut_0(E \times F)$ and the other from the translation by a diagonal 2-torsion element (\cite[Example 1]{MN}, see also \cite[Appendix A]{Mu0}).
In our situation of characteristic 2, this type of free involution $\varepsilon$  is given by  $x \mapsto c/x, y \mapsto d/y, z \mapsto z+1$  of \eqref{4D4-eq in char 2}.
The quotient Enriques surface $\Km(E \times F)/\varepsilon$  has a configuration of 10 $(-2)$-curves whose dual graph is  ${\rm Figure \ref{ExFdualgraph}}/\varepsilon$ given
in Figure \ref{Barth-Peters}.

\begin{figure}[h]
\begin{center}
\scalebox{0.9}{\xy 
@={(0,0),(10,10),(20,20),(30,10),(40,0),(30,-10), (20,-20),(10,-10),(10,0),(30,0)}@@{*{\bullet}};
(0,0)*{};(20,20)*{}**\dir{-};(20,20)*{};(40,0)*{}**\dir{-};(40,0)*{};(20,-20)*{}**\dir{-};(20,-20)*{};(20,-20),(0,0)*{}**\dir{-};
(0,0)*{};(10,0)*{}**\dir{-};(30,0)*{};(40,0)*{}**\dir{-};
\endxy}
\end{center}
\caption{}
\label{Barth-Peters}
\end{figure}

\noindent
This is the same as that of Enriques surfaces studied in Barth-Peters~\cite{BP} and 
this change of configuration is due to that of quotient singularity from  $16A_1$  to  $4D_4$.

In the special case $E=F$, two curves  $x=y$  and $xy=c$  of bidegree $(1,1)$  on $\bfP^1_x \times \bfP^1_y$  split over the double covering.
Hence we have four more $(-2)$-curves and the dual graph of $(-2)$-curves in Figure \ref{ExFdualgraph} extends to that in Figure \ref{ExEdualgraph}.
The 1-parameter family of $\Km(E\times E)/\varepsilon$  is now the Enriques surfaces of type I with only finite automorphisms in characteristic 2 constructed in Katsura-Kond\=o~\cite[Theorem 4.1]{KK1} (see also \cite[\S8.10]{DK1}).

\section{Kummer surfaces associated with the product of an ordinary elliptic curve and a supersingular elliptic curve}\label{sec6}

In this section we will consider Kummer surfaces associated with the product of an ordinary elliptic
curve and a supersingular elliptic curve.  

Let $E$ be an ordinary elliptic curve, $F$ a supersingular elliptic curve  and 
$\iota$ the inversion of the abelian surface $E\times F$.
It is known that the quotient surface $(E\times F)/\la \iota\ra$ is a $K3$ surface with two rational double points of type $D_8$
(Shioda \cite{Shioda}).
Let $X'''$ be the Kummer surface associated with $E\times F$, that is, the minimal resolution of 
the quotient surface. 
An automorphism of  $E \times F$ which commutes with $\iota$ descends to 
$(E\times F)/\la \iota\ra$.
Hence we have a natural injective homomorphism
\begin{equation}\label{natural hom}
E_{(2)} \times \Aut_0(F) \longrightarrow \Aut(X'''),
\end{equation}
where  $E_{(2)} \cong \bfZ/2\bfZ$  is the group of 2-torsion points.
The automorphism group  $\Aut_0(F)$  (as group variety) is the binary tetrahedral group  $SL_2(\bfF_3)$
(cf. \cite[Chap. 3, Remark (6.3)]{Huse}).

\begin{theorem}\label{MainProd3}
Let $E, F$ and  $X'''$ be as above. 
Then \eqref{natural hom} is an isomorphism and  $\Aut(X''')\cong \bfZ/2\bfZ\times SL_2(\bfF_3)$.
\end{theorem}

The Kummer surface  $X'''$  has 19 $(-2)$-curves forming the dual graph in Figure \ref{ExFssdualgraph}.  
\begin{figure}[h]
\begin{center}
\scalebox{0.8}{\xy 
(90,15)*{};
(-70,-15)*{};
@={(0,0),(7.5,0),(15,0),(22.5,0),(30,0),(37.5,0),(45,0),(52.5,0),(60,0),(-7.5,0),(-15,0),(-22.5,0),(-30,0),(-37.5,0),(-45,0),(-45,7.5),(-52.5,0),(-60,0),(45,7.5)}@@{*{\bullet}};
(0,0)*{};(15,00)*{}**\dir{-};(15,00)*{};(30,0)*{}**\dir{-};(30,0)*{};(45,0)*{}**\dir{-};
(0,0)*{},(-15,0)*{}**\dir{-};
(-15,0)*{};(-30,0)*{}**\dir{-};(-30,0)*{};(-60,0)*{}**\dir{-};
(-45,0)*{};(-45,7.5)*{}**\dir{-};(45,0)*{};(60,0)*{}**\dir{-};
(45,7.5)*{};(45,0)*{}**\dir{-};
(-60,-4)*{\alpha_{1}};(60,-4)*{\alpha_2};(0,-4)*{\alpha_{3}};
\endxy}
\end{center}
\caption{The dual graph: a product type (ordinary and supersingular case)}\label{ExFssdualgraph}
\end{figure}

\noindent
Here we denote by $\alpha_1, \alpha_2, \alpha_3$ 
the images of $\{0\}\times F, \{a\}\times F, E\times \{0\}$ on $X'''$ respectively, 
where $a$ is the non-zero 2-torsion point of $E$.
The remaining 16 $(-2)$-curves, forming $2D_8$ diagram, come from the resolution $X''' \to (E\times F)/\la \iota\ra$.
The first projection $E\times F \to E$ induces an elliptic fibration 
$\pi_1$ on $X'''$ with exactly two singular fibers of type 
$\tilde{E}_8$ and a section $\alpha_3$, and the second one gives an elliptic fibration $\pi_2$ with a unique singular fiber of type 
$\tilde{D}_{16}$ and two sections $\alpha_1, \alpha_2$.

\begin{lemma}\label{PLatticeExF} 
$\NS(X''')$ is generated by $19$ $(-2)$-curves 
in Figure {\rm \ref{ExFssdualgraph}} and is 
isomorphic to $U\oplus E_8\oplus E_8$.
\end{lemma}
\begin{proof}
Obviously $E\times F$ has Picard number 2 and hence $\rank(\NS(X'''))= 18$.  The existence of an elliptic fibration
with two singular fibers of type $\tilde{E}_8$ and a section implies that 19 $(-2)$-curves in Figure \ref{ExFssdualgraph}
generate $U\oplus E_8\oplus E_8$.  Since this lattice is unimodular, the assertion follows.
\end{proof}

It follows from Vinberg \cite{V} that the subgroup $W(X''')$ of 
$O(\NS(X'''))$ generated by
19 $(-2)$-classes in Figure \ref{ExFssdualgraph} is of finite index and the quotient group
$O(\NS(X'''))/\{\pm 1\}\cdot W(X''')$ is isomorphic to $\bfZ/2\bfZ$, the symmetry group 
of the dual graph which is induced from the translation by the non-zero 
2-torsion point of $E$.
On the other hand, we have the following lemma.

\begin{lemma}\label{numtrivialproduct3}
Let $r: \Aut(X''')\to O(\NS(X'''))$.  Then, 
$\Ker(r) \cong \Aut_0(F)$. 
\end{lemma}
\begin{proof} 
Consider 
the action of $g \in \{ 1_E\} \times \Aut_0(F)$ on $E\times F$.    
Since $g$ fixes both $\alpha_1$ and 
$\alpha_2$, it acts trivially on the dual
graph (Figure \ref{ExFssdualgraph}).  Thus  
we have $\Aut_0(F) \subset \Ker (r)$.

Now let $g\in \Ker(r)$, $g\ne 1$. 
First consider the elliptic fibration $\pi_1$ with
two singular fibers of type $\tilde{E}_8$ and a section $\alpha_3$.  

If $g$ acts trivially on the base of $\pi_1$,
then it induces an automorphism of a general fiber.  
Since a supersingular elliptic curve is unique, 
a general fiber is isomorphic to $F$ and hence $g \in \Aut_0(F)$.   

If $g$ acts on the base of $\pi_1$ non-trivially, 
then $g$ has two fixed points on the base which are the images of two
singular fibers.  Therefore $|g|$ is odd.  Next we consider the
elliptic fibration $\pi_2$ with a unique singular fiber of type 
$\tilde{D}_{16}$ with sections $\alpha_1, \alpha_2$. 
If $g$ acts on the base of $\pi_2$ trivially,
then $g$ acts on any fiber as an automorphism. Since $|g|$ is odd,  
this contradicts the fact $\Aut_0(E)\cong \bfZ/2\bfZ$.  
Finally assume that $g$ acts on the base of $\pi_2$ non-trivially.
Since $|g|$ is odd, $g$ fixes two points on the base one of which is
the image of the unique singular fiber, and 
the other is the image of a smooth
fiber $E_0$.  Thus $g$ acts on $E_0$ as an automorphism.  
If $g$ acts trivially on $E_0$, then 
$g$ also acts trivially on the base of $\pi_1$, 
because $E_0$ dominates the base of $\pi_1$ under $\pi_1$.  This is a contradiction.  Thus
$g$ acts on $E_0$ non-trivially which again contradicts the fact $\Aut_0(E_0)\cong \bfZ/2\bfZ$. 
\end{proof}

Theorem~\ref{MainProd3}  immediately follows from this lemma.

\section{Quartics with $4D_4$ singularities: complex version of Theorem~\ref{MainCurve}}\label{sec7}
First, recall that Kummer \cite[p.253]{Kum} gave the following equation of the 
Kummer quartic surface:
\begin{equation}\label{KummerOriginal}
\phi^2=16Kpqrs,
\end{equation}
where
$\phi = p^2+q^2+r^2+s^2 +2a(qr+ps)+2b(rp+qs)+2c(pq+rs)$
and the coefficients $a,b,c,K$ satisfy the equation $K=a^2+b^2+c^2-2abc -1$.
In this subsection we study the quartic surfaces with four rational double points of type $D_4$ over the complex number field $\bfC$ related to the Kummer's equation as above.
These surfaces appears naturally as {\it sub-quotients} in the 4-dimensional family of quartic surfaces
\begin{equation}\label{Heisenberg}
S: A(x^2y^2+z^2w^2) + B(x^2z^2+y^2w^2)+ C(x^2w^2 + y^2z^2)+ Dxyzw+E(x^4+y^4+z^4+w^4)=0
\end{equation}
with  $A, B, C, D, E \in \bfC$, which are invariant under the action of the Heisenberg group
\[
1 \to \mu_2 \to H \to \bar H \to 0, \quad \bar H \simeq(\bfZ/2\bfZ)^4
\]
(see e.g., \cite[\S 10.2]{DoC}). 

The quotient group $\bar H$  is a 4-dimensional vector space over the binary field $\bfF_2$ and carries a symplectic form  $\bar H \times \bar H \to \mu_2$  coming from the commutator relation.
$\bar H$ has 15 involutions and 15 G\"opel subgroups.
The fixed point locus of an involution is the union of two skew lines.
So we have 30 special lines as fixed loci of involutions.
If  $S$ has once a node outside these lines then it has 16 nodes and becomes a Kummer surface, which is characterized by  the Segre cubic equation 
\[
\Delta_0(A, B, C, D, E)=(4A^2+4B^2+4C^2-D^2-16E^2)E - 4ABC =0
\] 
(see e.g., \cite[\S 10.3]{DoC}).

A G\"opel subgroup $F \subset \bar H$ contains three involutions and their common fixed locus consists of four points.
The quartic surfaces with 4 nodes at these four points form a codimension one linear subfamily of \eqref{Heisenberg}.
For example, in the case $F$ is generated by $\text{diag}[1,-1,1,-1]$ and $\text{diag}[1,1,-1,-1]$, then the common fixed points are the coordinate points $(1,0,0,0), \ldots, (0,0,0,1)$,
and the linear subfamily 
\begin{equation}\label{Heisenbergwith E=0}
S_{E=0}: A(x^2y^2+z^2w^2) + B(x^2z^2+y^2w^2)+ C(x^2w^2 + y^2z^2)+ Dxyzw=0
\end{equation}
is defined by $E=0$.

Now we consider the quotient of  $S$  in \eqref{Heisenberg} by the action of  a G\"opel $F$, say 
\[
F= \la \text{diag}[1,-1,1,-1], \text{diag}[1,1,-1,-1]\ra.
\]
Then the quotient surface is given by
\begin{equation}\label{Heisenberg/Gopel}
S/F: \{A(xy+zw) + B(xz+yw)+ C(xw + yz) +E(x^2+y^2+z^2+w^2)\}^2 - D^2xyzw=0.
\end{equation}

\begin{remark}
(1) \ It is remarkable that 
in the Kummer case, that is, if $\Delta_0(A, B, C, D, E)=0$, the quotient $S/F$  is the
same as the one given 
by (\ref{KummerOriginal}) of the  Kummer quartic surface.  

(2)\ Assume the base field is of characteristic 2, then amusing things happen.  In this case, note that $\Delta_0= D^2E$.  If $D=0$, then
$S$ and $S/F$ are a (double) quadric. If $E=0$, then $S$ and $S/F$ are nothing but the
Kummer quartic surface in characteristic 2 given in (\ref{KummerQuarticEq}).
Moreover the quotient map $S\to S/F$ is the Frobenius map.
\end{remark}

If  $S$  belongs to the subfamily \eqref{Heisenbergwith E=0}, i.e., $E=0$, then  
$S_{E=0}/F$   has  $4D_4$ singularities at the coordinate points.
Let  $X \to S_{E=0}/F$ be the minimal resolution of $4D_4$ singularities.
Then $X$  has 20 $(-2)$-curves with the configuration given in Figure \ref{20curves}.
The standard Cremona transformation \eqref{CremonaOrdinary} exists characteristic freely,
and defines an involution of  $X$.
Moreover, since the involutions $\iota_{b,\delta}$, $\iota_{c,\delta}$ and $\iota_{d,\delta}$ are constructed in terms of elliptic pencils, they exist in our case also.
Thus the proof of Theorem~\ref{MainCurve} works in this situation and we have 

\begin{theorem}\label{MainCurveOverC}
Assume that  $X$ is Picard general, that is, the Picard number of $X$ is
$17$.
Let $G$ be the group generated by six involutions $\iota_{b,\delta}$, four involutions $\iota_{c,\delta}$ 
and eight involutions $\iota_{d,\delta}$.  Then 
\[
\Aut(X) \cong G \rtimes (\bfZ/2\bfZ)^3.
\]
Moreover, the eighteen involutions $\iota_{b,\delta}$, $\iota_{c,\delta}$ and  $\iota_{d,\delta}$ depends only on $\delta$.
\end{theorem}

\begin{remark}
The quotient of a quartic surface in \eqref{Heisenbergwith E=0} by the standard Cremona transformation is an Enriques surface of type  $E_7$ (\cite[\S8]{Mu2}), which has a numerically reflective involution (\cite{Mu1}).
\end{remark}

\section{Quartics with $2D_8$ singularities: Complex version of Theorem~\ref{MainCurvep-rank1} and a 3-parameter families in characteristic two}\label{Qw2D8}

Recall that a Picard general Kummer surface $X_1$ in Section~\ref{sec4} has 
Picard number 17.  So it is natural to search a 3-dimensional family of $K3$ surfaces with the same Picard lattice as $X_1$.
Our answer over $\bfC$ is the family of quartics
\begin{equation}\label{lift of Duquesne}
\{q(x,z)+yw\}^2 - \gamma^2x^2z^2 + x^2zw + xyz^2  = 0  \subset \bfP^3_{xyzw},
\end{equation}
where $q(x, z)=\beta x^2 + \alpha xz + z^2$  is a binary quadratic form and 
$\alpha, \beta (\ne 0), \gamma$ are constants.
In this section we first study this family of quartic surfaces.
The family \eqref{lift of Duquesne}
is equivalent to that of 
\begin{equation}\label{lift of Duquesne 2}
S_{a, b, c} : \{q_1(x, z) + yw\}\{q_2(x, z) + yw\}   + x^2zw + xyz^2  = 0  \subset \bfP^3_{xyzw}
\end{equation}
for pairs of quadratic forms $q_1(x, z) = bx^2 + axz + z^2$, 
$q_2(x, z) = bx^2 + (a+c)xz + z^2$ $(a, b, c \in \bfC, b\ne 0)$  whose difference is a constant multiple of  $xz$.

We also study modulo 2 reduction of this family, namely
\begin{equation}\label{Gen'dDuquesne}
S_{a, b, c} : b^2x^4+a^2x^2z^2  + z^4 + y^2w^2 + cxz(bx^2+axz+z^2+yw) + x^2zw + xyz^2  = 0 \subset \bfP^3
\end{equation} 
in characteristic 2,
which contains Duquesne's  hexanomial quartics \eqref{KummerEq1} 
at $c=0$.

\subsection{3-parameter family over $\bfC$}
Quartic surface  $S_{a, b, c}$ has $D_8$ singularity at $[0, 1, 0, 0]$  and $[0, 0, 0, 1]$.
We denote its minimal resolution by  $X_{a, b, c}$.
Since  $S_{a, b, c}$  has also two tropes cut out by  $x=0$  and $w=0$, we have
\begin{lemma}
$X_{a, b, c}$  contains a configuration of $18$ $(-2)$ curves with the same dual graph as $\Km(J(C))$ of Jacobian $J(C)$  with  $p$-rank $1$.
\end{lemma}
The dual graph of the configuration in this lemma is a hexadecagon with two legs and depicted by 18  {\it bullets} $\bullet$ in Figure~\ref{16gon2legs}.

\begin{figure}[h]
\begin{center}
\xy
(40,45)*{};
(-50,-5)*{};
@={(0,40),(10,40),(20,40),(30,40),(40,40),
(0,30),(10,30),(40,30),
(0,20),(40,20),
(0,10),(30,10),(40,10),
(0,0),(10,0),(20,0),(30,0),(40,0)
}@@{*{\bullet}};
(40,0)*{};(40,40)*{}**\dir{-};
(0,40)*{};(40,40)*{}**\dir{-};
(0,0)*{};(0,40)*{}**\dir{-};
(0,0)*{};(40,0)*{}**\dir{-};
(0,40)*{};(10,30)*{}**\dir{-};
(30,10)*{};(40,0)*{}**\dir{-};
(30,10)*{};(40,0)*{}**\dir{-};
(10,10)*{{\circ}};(30,30)*{{\circ}};
(0,0)*{};(9.5,9.5)*{}**\dir{=};
(30.5,30.5)*{};(40,40)*{}**\dir{=};
(10.5,10.5)*{};(29.5,29.5)*{}**\dir{--};
(8,14)*{-4};(28,34)*{-4};(20,24)*{4}
\endxy
\end{center}
\caption{Hexadecagon, two legs and two $(-4)$'s}
\label{16gon2legs}
\end{figure}
The surface 
$S_{a,b,c}$ is preserved by the linear transformation $[x, y, z, w] \mapsto [z, w, b x, b y]$ given in \eqref{translprank1}, which induces an involution $\varphi$  of  $X_{a, b, c}$.
The involution $\varphi$ is nothing but the translation of the elliptic fibration 
of  $X_{a, b, c}$ with $I_{16}$ of the hexadecagon by the 2-torsion section 
corresponding to two legs.

\begin{lemma}
The Picard lattice of  $X_{a, b, c}$  contains  $I_{1,16}^{ev}$ as primitive sublattice,
where  $I_{1,16}$ is the odd unimodular lattice of signature  $(1, 16)$  and  $I_{1,16}^{ev}$  is its even part {\rm (\S\ref{odd unimodular})}.
\end{lemma}
\begin{proof}
$\Pic (X_{a, b, c})$  contains  $U \oplus A_{15}$  
by virtue of the hexadecagonal elliptic fibration.
Since the fibration has a 2-torsion section,  $\Pic (X_{a, b, c})$ contains an index 2 overlattice  $\Lambda$  of $U \oplus A_{15}$.  Then 
$\Lambda$  is isomorphic to $I_{1,16}^{ev}$.
Since $I_{1,16}^{ev}$ has no nontrivial even overlattice, $\Lambda$  is primitive in  
$\Pic (X_{a, b, c})$.
\end{proof}

\begin{remark}
Both $I_{1,16}$ and $I_{1,16}^{ev}$ are reflective with Figure~\ref{16gon2legs} as their common Coxeter graph (Proposition~\ref{reflectivity inherits}).
\end{remark}

Now we assume that  $X_{a, b, c}$ {\it Picard general}, i.e., of Picard rank 17.
The projections of $S_{a,b,c}$ from two singular points induce two double covers $X_{a, b, c} \to \bfP^2$.
Hence we have two (covering) involutions  $pr_1, pr_2$  of  $X_{a,b,c}$.
 These involutions act on the Picard lattice ($\simeq I_{1,16}^{ev}$)  as reflection by 
 $(-4)$-roots, 
 depicted by white vertices $\circ$  in  Figure~\ref{16gon2legs}.
 These involutions generate an infinite dihedral group  $G_1$ in  $\Aut(X_{a, b, c})$. 
Furthermore, the Cremona transformation given by the equation 
 (\ref{Cremonaprank1}) 
 preserves  $S_{a,b,c}$  and induces an involution  $\sigma$  of  $X_{a,b,c}$.
The same argument as \S4.2  works for  $X_{a,b,c}$, and we have the following theorem.

\begin{theorem}\label{MainGen'dDuquesne/C}
Assume that  $X_{a,b,c}$  is Picard general.
Then  $\Aut(X_{a, b, c})$  is the semi-direct product 
$\la pr_1, pr_2 \ra \rtimes \la \sigma, \varphi \ra \cong G_1 \rtimes (\bfZ/2\bfZ)^2$.
\end{theorem}

\subsection{3-parameter family in characteristic 2}\label{Qw2D8-2}
In characteristic 2, two involutions, called {\it projections}, are given by
 
\begin{equation}\label{2 Cremonas 2}
\begin{split}
pr_1 : [x, y, z, w]  \dashrightarrow \left[x, y, z, w + \frac{x^2z}{y^2} + \frac{cxz}y\right], \\ 
pr_2 : [x, y, z, w]  \dashrightarrow \left[x, y + \frac{xz^2}{w^2} + \frac{cxz}w, z, w\right].
\end{split}
\end{equation}
For general $a, b, c$, the Picard number of  the minimal resolution $X_{a,b,c}$  of  
$S_{a, b, c}$ in \eqref{Gen'dDuquesne} is 17,
since the same holds for Kummer surfaces $S_{a,b,0}$ ($S_{a,b,0}$ is characterized 
in the family $S_{a,b,c}$ as the locus of Kummer quartic surfaces).
Therefore, the same argument as the previous subsection works for  $X_{a,b,c}$.
\begin{theorem}\label{MainGen'dDuquesne/p=2}
Assume that  $X_{a,b,c}$  is Picard general.
Then $\Aut(X_{a, b, c})$  is 
the semi-direct product 
$\la pr_1, pr_2 \ra \rtimes \la \sigma, \varphi \ra \cong G_1 \rtimes (\bfZ/2\bfZ)^2$.
\end{theorem}

\section{Double $\bfP^1 \times \bfP^1$ with $4D_4$ singularities: Complex version of Theorem~\ref{MainProd1}}\label{DQ4D4}

The complex analogy of Theorem~\ref{MainProd1}  is easy in this case.
Consider the double cover of $\bfP^1 \times \bfP^1$  defined by
\begin{equation}\label{4D4-eq}
w^2 =  \left(x + \frac{a}x\right) + b + \left(y + \frac{c}y\right)
\end{equation}
for general  $a, b, c \in \bfC$, where $(x, y)$  is an inhomogeneous coordinate.
The branch is the union of $x = 0, \infty$, $y = 0, \infty$  and the (smooth) curve 
\begin{equation}
(x^2+a)y + bxy + x(y^2+c) = 0 
\end{equation}
of bidegree $(2, 2)$ passing through the coordinate points $(0, 0), (0, \infty), (\infty, 0), (\infty, \infty)$.
Hence the double cover has four $D_4$-singular points over the coordinate points.
Let  $S$  be the minimal resolution.
Then $S$  contains 20 $(-2)$-curves whose dual graph is Figure~\ref{16gon4legs}.
\begin{figure}[h]
\begin{center}
\xy
(40,45)*{};
(-50,-5)*{};
@={(0,40),(10,40),(20,40),(30,40),(40,40),
(0,30),(10,30),(40,30),
(0,20),(40,20),
(0,10),(30,10),(40,10),
(0,0),(10,0),(20,0),(30,0),(40,0)
}@@{*{\bullet}};
(40,0)*{};(40,40)*{}**\dir{-};
(0,40)*{};(40,40)*{}**\dir{-};
(0,0)*{};(0,40)*{}**\dir{-};
(0,0)*{};(40,0)*{}**\dir{-};
(0,40)*{};(10,30)*{}**\dir{-};
(30,10)*{};(40,0)*{}**\dir{-};
(30,10)*{};(40,0)*{}**\dir{-};
(10,10)*{{\bullet}};(30,30)*{{\bullet}};
(0,0)*{};(10,10)*{}**\dir{-};
(30,30)*{};(40,40)*{}**\dir{-};
\endxy
\end{center}
\caption{Hexadecagon with four legs}
\label{16gon4legs}
\end{figure}
Obviously  $S$ has an action of  $(\bfZ/2\bfZ)^3$  
generated by three involutions 
\begin{equation}\label{3 involutions}
w\to -w, \quad  x\to a/x  \quad {\rm and} \quad y \to c/y.
\end{equation}
The crucial elliptic fibration to determine  $\Aut(S)$  is the one with the hexadecagon of  $\bullet$'s as $I_{16}$-type fiber,
which we denote by  $f_{16}: S \to  \bfP^1$.
Corresponding to two extra roots  $\delta_1, \delta_2$, we find two inversions $\iota_1, \iota_2$ with respect to two sections  $f_{16}$.

\begin{theorem}\label{CpxMainProd1}
Let $S$ and the action of $(\bfZ/2\bfZ)^3$ be as above
and assume that  $S$  is Picard general, i.e., of Picard rank $18$.
Let $G$ be the infinite dihedral group generated by two involutions $\iota_1, \iota_2$ preserving the fibration  $f_{16}$
Then 
\[\Aut(S) \cong G \rtimes (\bfZ/2\bfZ)^3.\]
\end{theorem}

\begin{remark} 
In characteristic 2,  the Artin-Schreier version \eqref{4D4-eq in char 2}
of \eqref{4D4-eq}  is the standard equation of Kummer surfaces of ordinary and product type.
\end{remark}

\begin{remark}\label{BPEnriques} 
The product of three involutions of \eqref{3 involutions}, namely  $\varepsilon : (x, y, w) \mapsto (a/x, c/y, -w)$, has no fixed points.
The quotient  $S/\varepsilon$  is  
the Enriques surfaces studied by Barth-Peters (cf. \cite[Example 2]{MN}).
\end{remark}

\section{Quartics with $6A_3$ singularities: complex version of Theorem~\ref{MainProd2}}\label{Q6A3}

The $K3$ surfaces studied in \S~\ref{subs5-2}, i.e., the minimal resolution of $\Km(E \times E)$, have a configuration of 22 $\bfP^1$'s whose dual graph is Figure~\ref{tetrahedron}.
Here we study a $K3$ surface  $X$  with such a configuration in more general setting.

\begin{figure}[h]
\begin{center}
\xy
(-110,0)*{}
@={(-40,40),(-70,-10),(0,-10),(-40,-30),
(-40,22.5), (-40,5),(-40,-12.5),
(-62.5,2.5), (-55,15),(-47.5,27.5),
(-10,2.5), (-20,15),(-30,27.5),
(-54.2,-10),(-37.5,-10),(-18.7,-10),
(-62.5,-15),(-55,-20),(-47.5,-25),
(-30,-25),(-20,-20),(-10,-15)
}@@{*{\bullet}};
(0,-10)*{};(-40,40)*{}**\dir{-};
(0,-10)*{};(-40,-30)*{}**\dir{-};
(0,-10)*{};(-70,-10)*{}**\dir{-};
(-40,-30)*{};(-40,40)*{}**\dir{-};
(-40,-30)*{};(-70,-10)*{}**\dir{-};
(-40,40)*{};(-70,-10)*{}**\dir{-};
(-44,41)*{4};
(-43,23.5)*{3};(-51,29)*{3};(-26,29)*{3};
(-44,5)*{2};(-59,16)*{2};(-16,16)*{2};
(-43,-14)*{1};(-66,4)*{1};(-7,4)*{1};
\endxy
\end{center}
\caption{Tetrahedral graph on 22 vertices and a member in $|D|$}
\label{tetrahedron}
\end{figure}

\begin{lemma}
If a $K3$ surface  $X$  (over an arbitrary field) has 22 $\bfP^1$'s on it and if their dual graph is Figure{\rm ~\ref{tetrahedron}}, 
then  $X$  is the minimal resolution of a quartic surface 
\begin{equation}\label{mirror}
(x+y+z+w)^4 = \lambda xyzw  \subset \bfP^3_{(x:y:z:w)}
\end{equation}
\end{lemma}
\begin{proof}
We birationally embed the $K3$ surface in $\bfP^3$ by the complete linear system  $|D|$, where  $D$  is the sum of 22 $\bfP^1$'s in the configuration.
Obviously we have  $(D^2)=4$.
We start our proof with the key observation that $X$  has four elliptic pencils $\Lambda_1, \ldots, \Lambda_4$ of type  $\tilde E_6 + \tilde A_{11}$.
Let  $f_i$, $i = 1, 2, 3, 4$, be their {\it elliptic parameters}, that is, $f_i$  is a rational function on  $S$ which has simple poles along $\tilde A_{11} \in \Lambda_i$  and  zeros along $\tilde E_6 \in \Lambda_i$ with prescribed multiplicities.
Then  $D_i := D +(f_i)$  belongs to  $|D|$ (see Figure~\ref{tetrahedron} for the coefficients of  $D_i$).
Especially,  $|D|$  has no base points, and the four divisors $D_1, \ldots, D_4 \in |D|$ determine a morphism $\varphi : X \to \bfP^3_{(x:y:z:w)}$,
for which we observe the following:
\begin{itemize}
\item Four plane sections $x=0, \ldots, w=0$  of  $\varphi(X)$  are lines $l_i$, $i = 1, 2, 3, 4$, counted with multiplicity 4.
\item The plane section of  $\varphi(X)$  corresponding to  $D$  is $\sum_1^4 l_i$.
\item Six $A_3$'s, i.e., chains of three $\bfP^1$'s, on  six edges of the tetrahedron are contracted by  $\varphi$.
\end{itemize}
Hence $\varphi$  is the birational morphism onto a quartic $(ax+by+cz+dw)^4 = \lambda xyzw$ for a suitable constant $\lambda$, where  $ax+by+cz+dw = 0$  is the plane cutting out  $\varphi(D)$ from $\varphi(X)$.
Taking $ax, \ldots, dw$ as new coordinates $x, \ldots, w$, we have our lemma.
\end{proof}

A  quartic surface in \eqref{mirror} is characterized by the presence of
\begin{itemize}
\item a plane section which is a complete quadrilateral, and 
\item six rational double points of type  $A_3$ at the intersection points.
\end{itemize}

Now we consider  $S$  in the pencil \eqref{mirror}  over  $\bfC$.
It is the quotient of the quartic surface in the {\it Dwork pencil}
\begin{equation}\label{mirror2}
x^4+y^4+z^4+w^4 = \lambda xyzw  \subset \bfP^3_{(x:y:z:w)}
\end{equation}
by the diagonal symplectic action of  $\mu_4^2$, and a mirror quartic surface in the sense of   \cite{Sei} (see also \cite{DoM}, \cite{NS}).
Let  $X$  be the minimal resolution of  $S$.
If  $\lambda$  is general, then the Picard lattice of  $X$  is generated by the 22 $\bfP^1$'s, namely
\begin{itemize}
\item the six edges of the quadrilateral  $xyzw = 0$ in the plane  $x+y+z+w = 0$, and 
\item the exceptional divisors over the six rational double points of type  $A_3$,
\end{itemize}
and isomorphic to  $U \oplus E_8 \oplus D_9$.
The proof of Theorem~\ref{MainProd2}  works in this situation and we have

\begin{theorem}
Let $X$ be the minimal resolution of a general member of  \eqref{mirror} and 
assume that $X$ is Picard general, i.e., the Picard number of $X$ is $19$.
Then
$\Aut(X)$ is isomorphic to $G\rtimes \mathfrak{S}_4$.
Here  $G$ is a subgroup of $\Aut(X)$ generated by involutions $\sigma_i$ $(i=1,2,3)$ (see below), inversions $\iota_{E,j}$ $(j=1,\ldots, 12)$ of suitable elliptic fibrations {\rm (cf. \S~\ref{subs5-2})}, and $\mathfrak{S}_4$  is the group of permutations of coordinates of the ambient projective space  $\bfP^3_{(x:y:z:w)}$. 
\end{theorem}

The following change of variables 
\begin{equation}
u = \frac{xy}{(x+y+z+w)^2}, \quad  v = \frac{yz}{(x+y+z+w)^2}, \quad   t = \frac{x+z}{x+y+z+w}
\end{equation}
transforms \eqref{mirror} into
\begin{equation}\label{special 4D4}
t(1-t) = u + \frac\lambda{u} + v+ \frac\lambda{v}
\end{equation}
(cf. \cite[(2.3)]{NS}).
Hence $X$  is the special case $a=c$ of \eqref{4D4-eq}.
The involution $\sigma_1$  in the theorem is  $(u,v,t) \mapsto (u, v, 1-t)$ in the new coordinate and induced from the cubic Cremona involution
 \begin{equation}\label{Cremonamirror}
 (x: y: z: w) \mapsto  (x(z+w)^2: y(z+w)^2: z(x+y)^2: w(x+y)^2) 
 \end{equation}
in the original coordinate of $\bfP^3_{(x:y:z:w)}$.
Other involutions $\sigma_2$ and $\sigma_3$  are obtained similarly by permuting coordinates.
There are three elliptic fibrations $g_i$ $(i=1,2,3)$ of type $2\tilde{D}_8$ (and with Mordell-Weil rank one).
Each $\sigma_i$  is the common inversion of two $g_j$  with $j \ne i$.

\begin{remark}
The theorem is essentially a consequence of the reflectivity of the odd unimodular lattice $I_{1, 18}$ and its even part  $I_{1, 18}^{ev}$ (\S\ref{odd unimodular}), which have  Figure~\ref{tetrahedron} as the essential part of their common Coxeter diagram.
\end{remark}

\begin{remark}
The Enriques surfaces of type I with only finite automorphisms in the sense \cite{Kondo1} is the quotient of our $K3$ surface \eqref{mirror}  by the cubic Cremona involution
 \begin{equation}
 (x: y: z: w) \mapsto  (y(z+w)^2: x(z+w)^2: w(x+y)^2: z(x+y)^2), 
  \end{equation}
  which is the composite of \eqref{Cremonamirror} and the permutation $(xy)(zw)$ of coordinates.
  It is also the quotient of  \eqref{special 4D4}  by  $\varepsilon : ( u, v, t) \mapsto (1/u, 1/v, 1-t)$ (see Remark~\ref{BPEnriques}).
\end{remark}

\begin{remark}\label{SI3}
Let  $f: E \to F$ be a degree 2 isogeny between elliptic curves and let 
$\Gamma \subset E \times F$  be its graph.
Then the divisor  $\Gamma - 2(0 \times F) - E \times 0$ is parpendicular to $U$-polarization and has self-intersection number $-4$.
Hence the abelian surface $A = E \times F$ has the same transcendental lattice with K3 surfaces $X$ in the theorem (cf. the last column of Table~\ref{table2}).
\end{remark}

\begin{remark}
The quartic surface in characteristic 2 defined by the equation (\ref{mirror}) (= 
 (\ref{mirror2})) is considered in
\cite{DK3} as an analogue of Desmic quartic surfaces over the complex numbers.
\end{remark}

\appendix

\section{{\bf The automorphism group of a generic Jacobian Kummer surface in characteristic different from two}}\label{OddChar}

The result of \cite[Theorem 7.4]{Kondo2} holds true in odd characteristic as well.
In this appendix, we prove it in the following form:

\begin{theorem}\label{OddCharThm}
Let $C$  be a curve of genus $2$ which is Picard general, i.e., $\rho(J(C)) = 1$.
Then the automorphism group of the Kummer surface $\Km(C)$ of the Jacobian variety  $J(C)$ is generated by the following$:$
\begin{enumerate}
\item  the $2$-elementary abelian group $2^5$  of order $32$ which is generated by translations $2^4$  and a switch involution,
\item  $16$ projections and $16$ correlations,
\item  $60$ G\"opel involutions, and
\item  $192$ Hutchinson-Weber involutions {\rm (\cite{Ohashi})}.
\end{enumerate}
\end{theorem}

We here version-up the original proof so that it works in odd characteristic as well by resolving problems related with Torelli type theorems and the primitive hull of 16 nodes.
They cause problems since the former needs analysis and the determination of the latter used topology.
Our solution here is this:
\begin{enumerate}
\item First we confirm that automorphisms in the theorem are constructed without Torelli type theorems in \S\ref{ConstInv}.
\item Next we show that the primitive hull of the lattice generated by 16 $(-2)$ nodal classes of a Kummer surface in its Picard lattice is the same as in characteristic 0 by a direct argument avoiding the use of topology in \S\ref{primitive hull}.
\item Torelli type theorems, including that of Nikulin, were used in the original proof (\cite[\S3]{Kondo2}).
 In \S\ref{Bolza-Igusa}, we replace them by the Bolza-Igusa classification of  the automorphism group $\Aut\ C$.
\end{enumerate}
The final subsection  \S\ref{JK-MOG} has nothing to do with positive characteristic problem.
The Jacobian Kummer lattice is embedded into the extended Leech lattice in the framework of \S\ref{Leech}.

\subsection{Construction of involutions}\label{ConstInv}
We recall here how automorphisms in the theorem  
are constructed and observe that the construction works well in positive odd characteristic.
\begin{enumerate}
\item 
The Kummer quartic surface  $\overline{\Km(C)} \subset \bfP^3$  is (projectively) self-dual.
More precisely, a switch (involution) is the composite of the Gauss map and a certain linear isomorphism  associated with theta characteristics $2D \sim K_C$ (\cite[\S5]{Ohashi}).

\item 
The Jacobian $J(C)$  is obtained from the symmetric product ${\rm Sym}^2 C$  of  $C$  by contracting the anti-diagonal  $A : = \{(x, \iota(x)) \,|\, x \in C$\}, being a $(-1)$  $\bfP^1$, where  $\iota$  is the hyperelliptic involution of  $C$.
The hyperelliptic covering $C \to \bfP^1$  induces the double covering  
$\Km(C) \to {\rm Sym}^2 \bfP^1 \simeq \bfP^2$,
whose covering involution is nothing but a projection.
The 16 correlations are well defined as the pullback of 16 projections 
by the Gauss map.

\item 
The diagonal $\Delta = \{(x, x) \,|\, x \in C\} \simeq C$  is mapped onto a conic  $Q \simeq \bfP^1$  in  $\bfP^2$.
The branch locus of  $\Km(C) \to {\rm Sym}^2 \bfP^1 \simeq \bfP^2$  is the union of 6 tangent lines  $l_1, \ldots, l_6$ of  $Q$ at the Weierstrass points $p_1, \ldots, p_6$.
A suitable quadratic Cremona transformation with center $l_a \cap l_b$,  $l_c \cap l_d$ and $l_e \cap l_f$  with  $\{a, \ldots, f\} = \{1, \ldots, 6\}$  lifts to an involutive automorphism of  $Km(C)$ (\cite[Proposition 5.1]{Mu1}).
These are nothing but the (Hutchinson-) G\"opel involutions corresponding to 15 G\"opel tetrads containing the origin.
Other 45 G\"opel involutions are obtained by taking conjugation by translations.

\item
The construction in \cite[\S7]{Ohashi} of Hutchinson-Weber involution from a Weber hexad works in odd characteristic as well.
\end{enumerate}

\subsection{The Kummer lattice $\Lambda_{Kum}$}\label{primitive hull}
Let  $A$  be an abelian surface and $\Km(A) \to \overline{\Km(A)}$  the minimal resolution of 16 nodes, corresponding to $A_{(2)}$, the group of 2-torsion points.
We put $I := A_{(2)}$, the index set of 16 nodes and exceptional divisors $e_i$ over them, and regard it as the 4-dimensional vector space over the binary field $\bfF_2$.

The Picard lattice of  $\Km(A)$ contains the negative definite sublattice  $\bigoplus _{i \in I} \bfZ e_i$  of rank 16.
The following was proved in \cite[Appendix to \S5]{P-SS} topologically in characteristic 0.

\begin{lemma}
The $16$ divisor class
\[
\sum_{i \in I} e_i \quad {\rm and} \quad \sum_{j \in J} e_j
\]
are divisible by $2$  in  $\Pic (\Km\ A)$, where  $J$ runs over all $3$-dimensional subspaces $J \subset I$.
\end{lemma}
\begin{proof}(1) The first one is obvious since the blow-up of  $A$ with center  $A_{(2)}$  is the double cover of  $\Km(A)$ with branch  $\sum_{i \in I} e_i$.

(2) To prove the divisibility of the rest, we consider the isogeny  $B_J \to A$  of degree 2 corresponding to  $J$,
that is, $B_J \to A$ is the dual of the quotient morphism of  $\hat A \to \hat A/a_J$ of the dual Abelian surface  $\hat A$ by the non-zero 
2-torsion point $a_J$ dual to $J$.
Then the induced map  $B_{J(2)} \to A_{(2)}$  on 2-torsions is a 2--1 surjection onto  $J$.
The induced morphism  $\overline\Km(B_J) \to \overline\Km(A)$  is ramified over 8 nodes indexed by $I \setminus J$.
Hence  $\sum_{j \not\in J} e_j$  is divisible by 2 in  $\Pic (\Km\ A)$.
So is the complement $\sum_{j \in J} e_j$  by (1).
\end{proof}

We define  $\Lambda_{Kum}$  as the overlattice obtained from $\bigoplus _{i \in I} \bfZ e_i$  by adding the halves of the 16 divisor classes in the lemma.
This lattice does not depend on $A$, and we denote it by $\Lambda_{Kum}$.
The discriminant group is the 2-elementary abelian group $2^6$ and canonically isomorphic to the 2nd exterior product $\bigwedge^2 I$  over $\bfF_2$.
Moreover, the discriminant form is the square map  
\begin{equation}\label{square map}
\bigwedge^2 I  \to \bigwedge^4 I \simeq \bfF_2, \quad w \mapsto w \wedge w.
\end{equation}

\begin{lemma}
The sublattice $\Lambda_{Kum} \subset \Pic (\Km\ A)$ is primitive.
\end{lemma}
\begin{proof}
Assume that $\Lambda_{Kum}$  is imprimitive.
Then there exists a nonzero isotropic element  $w$  in the discriminant group 
which vanishes in that of  $\Pic (\Km\ A)$.
By \eqref{square map}, $w$  is the Pl\"ucker coordinates of a 2-dimensional subspace  $\{0, a, b, c\}$  of  $I$.
This means that the sum $e_0+e_a+e_b+e_c$  is divisible by 2 in $\Pic (\Km\ A)$.
This contradicts the fact that the sum of disjoint $n$  $\bfP^1$'s on a $K3$ surface is divisible by 2 only when $n = 8, 16$ corresponding to $K3$ and abelian coverings.
This fact is due to Nikulin~\cite[\S1]{NKum} over the complex number field and the fact $n\ne 4$ remains to be true in odd characteristic (c.f. \cite[Cor. 3.3]{KaKoSch}).
\end{proof}

\begin{remark}
The above construction of  $\Lambda_{Kum}$  is nothing but Construction A in the sense of \cite[Chap. 5]{CS} from the first-order Reed-Muller code  [16, 5, 8].
Construction B from the same code yields the Barnes-Wall lattice of discriminant 256 (\cite[Chap. 4 \S10]{CS}).
\end{remark}

From now on we assume that  $A$  is the Jacobian $J(C)$ of a smooth curve  $C$ of genus 2.
Twice of the theta divisor descends to a polarization of degree 4 on the singular Kummer surface $\overline{\Km(C)}$.
We denote its pull-back on  $\Km(C)$  by $h$.
The Picard lattice of  $\Km(C)$  contains the (orthogonal) direct sum  $\bfZ h \oplus \Lambda_{Kum}$.
The Kummer quartic surface  $\overline{\Km(C)} \subset \bfP^3$  has a trope, that is, a double conic passing through 6 nodes.
Let $T \subset \Km(C)$  be the strict transform of its reduced part.
Then we have $2T + \sum_{i=1}^6 e_i \in  |h|$.
In particular, $h - \sum_{i=1}^6 e_i \in \bfZ h \oplus \Lambda_{Kum}$  is divisible by 2 in $\Pic (\Km\ C)$.
We denote the overlattice of  $\bfZ h \oplus \Lambda_{Kum}$  augmented by $(h - \sum_{i=1}^6 e_i)/2$  by  $\Lambda_{JK}$.

By the above lemma, we have

\begin{proposition}\label{JK-lattice}
The embedding $\Lambda_{JK} \subset \Pic (\Km\ C)$  is primitive.
\end{proposition}
In particular, we have the following:
\begin{corollary}
If  $C$  is Picard general, then  $\Pic (\Km\ C)$  is isomorphic to  $\Lambda_{JK}$  .
\end{corollary}

\subsection{Weaker form}\label{weaker form}
By this corollary the argument of \cite{Kondo2,Ohashi}  works in our case.
Namely, any automorphism $g$ of  $\Km(C)$ becomes an automorphism  $g'$  preserving the Conway-Borcherds cone after composing with involutions (2), (3) and (4) in Theorem~\ref{OddCharThm}.
Composing with a switch (involution) if necessary, $g'$  preserves the projective embedding  $\overline{\Km(C)} \hookrightarrow \bfP^3$  
as Kummer quartic.
For a suitable translation $t \in 2^4$, the composite $tg'$  preserves the node $N_0$  corresponding to the origin of $J(C)$.
Projecting $\overline{\Km(C)} \subset \bfP^3$ from $N_0$, we have the double covering  $\Km(C) \to \bfP^2$ described in  \S\ref{ConstInv}.
Hence $tg'$  is induced from a {\it reduced automorphism} of  $C$, that is, an automorphism of  the projective line  $\bfP^1$ preserving the 6 Weierstrass points  $p_1, \ldots,  p_6$.

The symmetry of the Conway-Borcherds cone is the semi-direct product of  the 2-elementary abelian group  $2^5$ by the symmetric group  $\frakS_6$  of degree 6, which permutes $p_1, \ldots,  p_6$ as monodromy.
Therefore, we have

\begin{theorem}\label{WeakerThm}
Let $C$  be a Picard general curve of genus $2$.
Then the automorphism group of the Jacobian Kummer surface $\Km(C)$ is generated by the involutions {\rm (2), (3), (4)} in Theorem~\ref{OddCharThm} and 

{\rm (1')}  the semi-direct product of  the $2$-elementary abelian group  $2^5$ by the reduced automorphism group $\Aut_0(C):= \Aut(\bfP^1; p_1, \ldots,  p_6)$ of  $C$.
\end{theorem}

\subsection{Proof of Theorem~\ref{OddCharThm}}\label{Bolza-Igusa}

It suffices to show that  $\Aut_0(C) \ne \{id\}$ implies $\rho(J(C)) > 1$.
We prove this using the classification of the reduced automorphism groups by Igusa~\cite[\S8]{Ig}, which algebraizes that of Bolza~\cite[Table on p. 70]{Bol} in characteristic 0.
There are six cases where  $\Aut_0(C) \ne \{id\}$ according to them:

{\small
  \begin{tabular}{c|c|c|c|c|c|c}\\
  Igusa \# & (1)  & (2)       & (3)       & (4)           & (5) & (6) \\\hline
  $\Aut_0$  & $C_2$  & $D_6$ & $C_2 \times C_2$ & $D_{12}$ & $\frakS_4$& $C_5$\\\hline
  Bolza form  &$x^6+\alpha x^4+\beta x^2+1$ & $x^6+ \alpha x^3+1$ & $x(x^4+\alpha x^2 +1)$ & $x^6+1$ &  $x(x^4+1)$ & $x^5+1$\\
in char. 0 &&&&&&\\
  \end{tabular}
  }
  \smallskip
  
  \noindent
with understanding that the dihedral group $D_{12}$ of order 12 in (4) does not appear in characteristic 3, and  $\Aut_0$ in (6) enlarges from the cyclic group $C_5$ of order 5 to $PGL(2, \bfF_5)$ of order 120 in characteristic 5.
 
In all cases except (6),  $C$  has a {\it bi-elliptic} involution $\sigma$, i.e.,  the quotient  $E = C/\sigma$ is elliptic.
For example, in the case $C: y^2 = x(x^4+\alpha x^2 +1)$  of  (4), such an involution $\sigma$  is given by  $(x, y) \mapsto (1/x, y/x^3)$.
The Jacobian $J(C)$  contains  $E$  as Abelian subvariety and hence it is not Picard general.
In the exceptional case (6), the endomorphism ring $\End^0(J(C))$  of  the Jacobian  contains the cyclotomic field  $\bfQ(\zeta)$ with  $\zeta^5=1$.
By the classification of endomorphism algebras of Abelian varieties (see e.g., \cite[Chap. IV]{M1}), 
the real quadratic subfield $\bfQ(\sqrt{5})$  belongs to the N\'eron--Severi group of  $J(C)$,
and hence $J(C)$  is not Picard general either.
Thus we have proved Theorem~\ref{OddCharThm}.

\begin{remark}
In the last case (6) the moduli point of  $J(C)$  as principally polarized abelian surface (ppas for short) belongs to the Humbert divisor of discriminant 5 (\cite[Chap. IX \S4]{G}).
The divisor consists of Comessatti surfaces, i.e., ppas's with real multiplication by  $\sqrt{5}$, which is crucial in studying Hutchinson-Weber involution (\cite{Ohashi2}).
\end{remark}

\subsection{The Jacobian Kummer lattice in the extended Leech lattice in the framework of \S\ref{Leech}}\label{JK-MOG}
In the original paper \cite[\S4]{Kondo2} a primitive embedding $\Lambda_{JK} \subset \II_{1,25} = U \oplus  \Lambda$  
was constructed using the fact that Golay code is the quadratic residue code modulo 23.
Here, taking advantage of our treatment in \S\ref{Leech}, we construct the embedding of $\Lambda_{JK}$ using MOG related with the projective plane  $\bfP^2(\bfF_4)$ over the quaternary field as explained there.

We find  9 Leech roots such that their orthogonal complement in  $U \oplus  \Lambda$  is $\Lambda_{JK}$.
More concretely, we find a standard set of 32 Leech roots of Kummer configuration  $16_6$ in their complement.
Our choice of 9 Leech roots is the following:
\[
(2, 1, [\widehat{\infty_\infty}]), \quad \emptyset : = (-1, 1, 0), \quad (2, 1, [\widehat{\infty_0}]).
\]
and
\[
(1,1, [\I]), (1,1, [\II]), (1,1, [\III]), (1,1, [\infty_1]), (1,1, [\infty_{\omega}]), (1,1, [\infty_{\bar\omega}]).
\]
The first three form $A_3$.
The rest is orthogonal to each other, and orthogonal to the first $A_3$.
Hence the total generates a root lattice  $R$ of type $A_3 + 6A_1$.

\begin{figure}[ht]
 \centering
\begin{tabular}{ccc}
$\emptyset$  && \\ 
 $/$ \ \ $\backslash$ && \\   \hline
$[\widehat{\infty_\infty}]$  $[\widehat{\infty_0}]$  & $\ [P] \ $ \ \ $\ \circ \ $ & $\ \circ \ $ \ \ $\ \circ \ $ \\ \hline
[I] \ \ \ \ \ $[\infty_1]$  & $\ \ \circ \ $ \ \ \ \ \ $\ \cdot \ $ & $\ \cdot \ $ \ \ $\ \cdot \ $ \\     \hline 
[II] \ \ \ \ $[\infty_\omega]$  & $\ \ \circ \ $ \ \ \ \ \ $\ \cdot \ $ & $\ \cdot \ $ \ \ $\ \cdot \ $ \\  \hline
[III] \ \ \ $[\infty_{\bar\omega}]$  & $\ \ \circ \ $ \ \ \ \ \ $\ \cdot \ $ & $\ \cdot \ $ \ \ $\ \cdot \ $  \\
\hline
\end{tabular}
\caption{$R =  A_3 + 6A_1$, $\bfA^2(\bfF_4)$ and Kummer configuration}\label{Kummer in Leech}
\end{figure}

The following is easy to verify:

\begin{lemma}
The Leech roots in the orthogonal complement of  $R$  in $U \oplus  \Lambda$ are the following $16+16:$ 
\begin{enumerate}
\item  one $16$ are  $(1,1, [P])$, where  $P$  runs over all points of the affine plane  $\bfA^2(\bfF_4)$, and
\item  the other $16$ are $(1, 1, [Q(P)])$, where $P$  is the same as above and $Q(P)$  is the symmetric difference of two lines passing through  $P$ parallel to $x$- and $y$-axes in $\bfP^2(\bfF_4)$, 
\end{enumerate}

\noindent
which form a configuration of Kummer  $16_6$.
\end{lemma}
Octads $Q(P)$'s  in the lemma are of the last type (0+8) in Proposition~\ref{MOGoctads}.
For example, $Q(P)$ for the point $P(0,0) \in \bfA^2(\bfF_4)$  is the union of six $\circ$ and the pair of points of infinity $\infty_\infty, \infty_0$ in the direction of  $x$- and $y$-axes as shown in Figure~\ref{Kummer in Leech}.
The vector $[Q(P)]$  has 2 on the eight points $Q(P)$ and 0's elsewhere in the notation of \S\ref{Leech}.

\bigskip

\bigskip

\section{{\bf KUMMER QUARTIC SURFACES IN CHARACTERISTIC TWO} }

\bigskip

\begin{center}
SHIGERU MUKAI
\end{center}

\bigskip

\subsection{A normal form}
Le $C$  be a curve of genus 2 over an algebraically closed field  $k$ of characteristic 2.
We assume that $C$ is ordinary, that is, the hyperelliptic double covering  $C \to \bfP^1$  has 3 Weierstra{\ss} points.
Taking the affine coordinate such that $x= 0, 1, \infty$ are its Weierstra{\ss} points, we have
\[
C: y^2+(x^2+x)y +f_6(x)=0
\]
for a sextic polynomial $f_6(x)$. 
We have a freedom of replacing the variable $y$ with $y + d(x)$  for a cubic polynomial $d(x)$.
By this change of variables, we can eliminate the term of $x^5, x^3, x$ from  $f_6(x)$, 
and obtain a simpler form
\[
C: y^2+(x^2+x)y + f_3(x)^2 = 0.
\]
Furthermore, by replacing  $y$ with $y + e(x^2+x)$ for a suitable constant $e \in k$, we can eliminate the quadratic term from  $f_3(x)$, and reach a normal form
\begin{equation}\label{perfect square}
C: y^2+(x^2+x)y + (ax^3+bx+c)^2 = 0
\end{equation}
for constants $a, b, c \in k$.
Since $C$  is smooth at $(x, y) = (0, c)$, $c$ is non-zero.
Similarly, we have $a \ne 0$  and $a+b+c \ne 0$.

\subsection{Theorem and Igusa form}
The purpose of this appendix is to show the following for this ordinary curve $C$ of genus 2.

\begin{theorem}\label{isogeny in char 2}
The Kummer surface of the Jacobian of $C$  is isomorphic to the minimal resolution of the quartic surface
\[
\{(xw+cyz) + (yw+\tilde bxz) + (zw+axy)\}^2 + xyzw= 0, 
\]
in  $\bfP^3_{(xyzw)}$, where we put $\tilde b = a+b+c$.
\end{theorem}

We apply the theorem to a curve of Igusa normal form:
\begin{equation}\label{Igusa reduction}
u^2+u = \frac\alpha{x} + \frac\beta{x+1} + \gamma x 
\end{equation}
or equivalently,
\[
y^2+(x^2+x)y + \gamma x^3(x+1)^2 + \alpha x(x+1)^2 + \beta x^2(x+1) = 0.
\]
Changing $y$ with $y+ \gamma x^3 + (\alpha+\beta+\gamma)x+\alpha$, we have
\[
y^2+(x^2+x)y + \gamma^2 x^6 + \gamma x^4+ \{(\alpha+\beta+\gamma)^2+\gamma\}x^2 + \alpha'^2 = 0.
\]
Further, changing $y$  with $y+\varepsilon(x^2+x)$ with $\varepsilon^2+\varepsilon=\gamma$, we reach a normal form
\[
y^2+(x^2+x)y + \{\gamma x^3 + (\alpha+\beta+\gamma)x + \alpha\}^2 = 0
\]
in the sense above.
By applying Theorem \ref{isogeny in char 2}, and 
finally, by replacing $x, y, z, w$ by $\sqrt{\alpha}x$, $\sqrt{\beta}y$, $\sqrt{\gamma}z$, $\sqrt{\alpha \beta \gamma}w$, respectively and dividing by $\alpha \beta \gamma$, we have

\begin{corollary}
For a curve of genus $2$ with Igusa normal form \eqref{Igusa reduction},
the Jacobian Kummer surface is
\[
\alpha(xw+yz)^2 + \beta(yw+xz)^2 + \gamma(zw+xy)^2 + xyzw= 0.
\]
\end{corollary}

\subsection{Proof of Theorem~\ref{isogeny in char 2}}
The function field of the self product  $C \times C$ is generated by two elements 
\[
u_1=\frac{y_1}{x_1^2+x_1}, \quad u_2=\frac{y_2}{x_2^2+x_2}
\]
with two relations
\[
u_i^2+u_i + \left(\frac{ax_i^3+bx_i+c}{x_i^2+x_i}\right)^2 = 0, \quad  i=1, 2
\]
over the purely transcendental extension $k(x_1, x_2)$ of  $k$.
Hence the function field of the Jacobian is generated by
\[
u : =\frac{y_1}{x_1^2+x_1} + \frac{y_2}{x_2^2+x_2}
\]
with relation
\[
u^2+u + \left(\frac{ax_1^3+bx_1+c}{x_1^2+x_1}\right)^2 + \left(\frac{ax_2^3+bx_2+c}{x_2^2+x_2}\right)^2 = 0
\]
over  $k(x_1+x_2, x_1x_2)$.
Putting  $t= (x_1^2+x_1)(x_2^2+x_2)u$, we have
\[
t^2 + (x_1^2+x_1)(x_2^2+x_2)t + (x_1+x_2)^2 \{ax_1x_2(x_1+x_2+x_1x_2)+bx_1x_2+c(x_1+x_2+1)\}^2.
\]
Consider the compactification of the affine 2-space  $\bfA^2$ with coordinate  $(x_1x_2, x_1+x_2)$ by the projective plane  $\bfP^2 = {\rm Sym}^2 \bfP^1$  with the homogeneous coordinate  $(x:y:z) = (x_1x_2: (x_1+1)(x_2+1): 1)$.
Then we have
\[
t^2 +xyzt+ (x+y+z)^2(axy + \tilde b xz + cyz)^2 = 0,
\]
the equation of (Artin-Schreier) double cover of  $\bfP^2$ with branch  $xyz=0$.
Putting $t = (x+y+z)^2w$, we have 
\[
(x+y+z)^2w^2 +xyzw+ (axy + \tilde b xz + cyz)^2 = 0,
\]
which is the quartic \eqref{isogeny in char 2} in the theorem.
Geometrically speaking, the above double plane is the projection from the double point  $(1,0,0,0)$.

\end{document}